\documentclass[11pt,fleqn]{amsart}

\usepackage{epsfig}
\usepackage{amsmath}
\usepackage{amssymb}
\usepackage{amscd}
\usepackage{latexsym}
\usepackage{mathtools}
\usepackage{tabularx}
\usepackage{blkarray}
\usepackage{a4wide}
\usepackage[usenames]{color}
\usepackage{mathdots}
\usepackage{mleftright}
\usepackage{subfigure}
\usepackage{url}
\usepackage{tikz,tikz-cd}
\usetikzlibrary{decorations.markings}
\usetikzlibrary{calc}
\usetikzlibrary{3d}
\usepackage[matrix,arrow,tips,curve]{xy}
\usepackage{pb-diagram,pb-xy}
\usepackage{verbatim}
\usepackage{mathrsfs}
\usepackage{color}
\usepackage{scalerel} 
\usepackage{mathabx}
\usepackage{enumitem}

\setlist[enumerate]{label=\textnormal{(\arabic*)}}

\newcommand{\cell}{\ssub[-2pt]{\mathcal Z}!}
\newcommand{\Gr}{\mathrm{Gr}}
\definecolor{cadmiumgreen}{rgb}{0.0, 0.42, 0.24}
\usepackage[
colorlinks, citecolor=cadmiumgreen,
pdfauthor={Omid Amini, Lucas Gierczak}, 
pdftitle={Combinatorial flag arrangements},
pdfstartview ={FitV},
]{hyperref}

\usepackage{todonotes} %% use disable option to disable

\newtheorem{thm}{Theorem}[section]

\newtheorem{lemma}[thm]{Lemma}

\newtheorem{prop}[thm]{Proposition}

\newtheorem{question}[thm]{Question}
\newtheorem*{question*}{Question}

\theoremstyle{definition}
\newtheorem{defii}[thm]{Definition}
\newenvironment{defi}
	{\pushQED{\qed}\defii}
	{\popQED\enddefii}
	\newtheorem{convention}[thm]{Convention}
\newenvironment{conv}
	{\pushQED{\qed}\convention}
	{\popQED\endconvention}
\newtheorem{remm}[thm]{Remark}
\newenvironment{remark}
	{\pushQED{\qed}\remm}
	{\popQED\endremm}
\newtheorem{exx}[thm]{Example}
\newenvironment{example}
	{\pushQED{\qed}\exx}
	{\popQED\endexx}

\newcommand{\Z}{\ssub{\mathbb{Z}}!}

\newcommand{\R}{\mathbb{R}}

\newcommand{\N}{\mathbb{N}}

\newcommand{\ord}{\ssub{\mathrm{ord}}!}

\newcommand{\Rat}{\mathrm{Rat}}

\renewcommand{\k}{\kappa}

\newcommand{\grdeg}{{m}}

\newcommand{\rest}[1]{\raisebox{-1pt}{$\vert$}_{#1}}

\newcommand{\st}{\bigm|} % such that in sets
\newcommand{\rkstandard}{\rk^{\scaleto{\mathrm{st}}{4.5pt}}}

\newcommand{\rquot}[2]{#1\big/#2} %%% Quotient space

\usepackage{manfnt}

\newcommand{\cube}{\mbox{\,\mancube}}

\makeatletter
\newsavebox\myboxA
\newsavebox\myboxB
\newlength\mylenA

\newcommand*\overbar[2][0.65]{%
	\sbox{\myboxA}{$\m@th#2$}%
	\setbox\myboxB\null% Phantom box
	\ht\myboxB=\ht\myboxA%
	\dp\myboxB=\dp\myboxA%
	\wd\myboxB=#1\wd\myboxA% Scale phantom
	\sbox\myboxB{$\m@th\overline{\copy\myboxB}$}% Overlined phantom
	\setlength\mylenA{\the\wd\myboxA}% calc width diff
	\addtolength\mylenA{-\the\wd\myboxB}%
	\ifdim\wd\myboxB<\wd\myboxA%
		\rlap{\hskip .9\mylenA\usebox\myboxB}{\usebox\myboxA}%
	\else
		\hskip -0.5\mylenA\rlap{\usebox\myboxA}{\hskip 0.5\mylenA\usebox\myboxB}%
	\fi}

\NewDocumentCommand{\ssub}{O{0pt} O{.9} m t! e{_^}}{
	#3%
	\IfValueT{#5}{
		\IfBooleanTF{#4}{\sb{\hspace{#1}\scaleobj{#2}{#5}}}{\sb{#5}}
	}
	\IfValueT{#6}{
		\IfBooleanTF{#4}{\sp{\hspace{#1}\scaleobj{#2}{#6}}}{\sp{#6}}
}
}

\NewDocumentCommand{\ssubb}{O{0pt} O{0pt} O{.9} m t! e{_^}}{
	#4%
	\IfValueT{#6}{
		\IfBooleanTF{#5}{\sb{\hspace{#1}\scaleobj{#3}{#6}}}{\sb{#6}}
	}
	\IfValueT{#7}{
		\IfBooleanTF{#5}{\sp{\hspace{#2}\scaleobj{#3}{#7}}}{\sp{#7}}
}
}

\ExplSyntaxOn
\NewDocumentCommand{\tossub}{o o m}{
	\expandafter\let\csname old\cs_to_str:N #3\endcsname#3
	\renewcommand#3%
	{\ssub[#1][#2]{\csname old\cs_to_str:N #3\endcsname}}
}
\ExplSyntaxOff

\tossub\omega
\tossub\Rat

\newcommand{\ssF}{\ssub{F}!}
\newcommand{\order}{\scaleto{\mathcal O}{4pt}}

\newcommand{\ssa}{\ssub{a}!}
\newcommand{\ssb}{\ssub{b}!}
\newcommand{\ssx}{\ssub{x}!}
\newcommand{\ssy}{\ssub{y}!}
\newcommand{\ssn}{\ssub{n}!}
\newcommand{\ssrho}{\ssub{\rk}!}
\newcommand{\ssi}{\ssub{i}!}
\newcommand{\ssI}{\ssub{I}!}
\newcommand{\ssJ}{\ssub{J}!}
\newcommand{\ssd}{\ssub{d}!}
\newcommand{\sse}{\ssub{e}!}
\newcommand{\ssvare}{\ssub{\varepsilon}!}

\newcommand{\ssv}{\ssub{v}!}
\newcommand{\ssc}{\ssub{c}!}

\newcommand{\ssY}{\ssub{Y}!}

\newcommand{\sst}{\ssub{t}!}

\newcommand{\ssfilt}{\ssubb[-1pt]{\filt}!}
\newcommand{\ssfiltg}{\ssubb[-.5pt]{\filtg}!}

\newcommand{\ssS}{\ssub{S}!}
\newcommand{\sss}{\ssub{s}!}
\newcommand{\ssp}{\ssub{p}!}

\newcommand{\ssR}{\ssub{R}!}
\newcommand{\ssr}{\ssub{\varrho}!}
\newcommand{\ssL}{\ssub{L}!}
\newcommand{\unda}{\ssub{\underline a}!}
\newcommand{\undb}{\ssub{\underline b}!}
\newcommand{\undc}{\ssub{\underline c}!}
\newcommand{\undd}{\ssub{\underline d}!}
\newcommand{\unde}{\ssub{\underline e}!}

\newcommand{\undx}{\ssub{\underline x}!}
\newcommand{\undy}{\ssub{\underline y}!}
\newcommand{\undz}{\ssub{\underline z}!}
\newcommand{\undw}{\ssub{\underline w}!}
\newcommand{\ssw}{\ssub{w}!}
\newcommand{\ssk}{\ssub{k}!}

\newcommand{\ssP}[1]{\ssub{P}!_{\hspace{-2pt}#1}}

\newcommand{\undp}{\ssub{\underline p}!}
\newcommand{\undq}{\ssub{\underline q}!}
\newcommand{\undu}{\ssub{\underline u}!}

\newcommand{\mat}{\ssub{\mathrm{M}}!}
\newcommand{\polymat}{\ssub{\mathrm{P}}!}

\newcommand{\rank}{\ssub{\mathrm{rank}}!}
\newcommand{\VS}{\mathrm{H}}
\newcommand{\filt}{\mathrm{F}}
\newcommand{\filtg}{\mathrm{G}}
\newcommand{\hcube}[2]{\ssub{\cube}_{#1}^{^{#2}}}

\newcommand{\rk}{\ssub{\mathbf{r}}!}
\newcommand{\undr}{\underline \varrho}

\newcommand{\codim}{\mathrm{codim}}

\newcommand{\matricube}{\ssub{\mathscr M}!}
\newcommand{\unimatricube}{\ssub{\mathscr U}!}

\newcommand{\Fl}{\mathscr F}
\newcommand{\Ind}{\mathscr I}
\newcommand{\Jind}{\mathcal J}

\newcommand{\dvs}{n} %%%% maximum width (dimension of the ambient vector space in representable case)

\newcommand{\subfaceeq}{\preceq}
\newcommand{\supfaceeq}{\succeq}
\newcommand{\subface}{\prec}
\newcommand{\ssubface}{\mathbin{\mathchoice
	{\subface\!\!\!\cdot}%
	{\subface\!\!\!\cdot}%
	{\subface\!\cdot}%
	{\subface\!\cdot}%
}} % covering elements
\newcommand{\supface}{\succ}
\newcommand{\ssupface}{\mathbin{\mathchoice
	{\cdot\!\!\!\supface}%
	{\cdot\!\!\!\supface}%
	{\cdot\!\supface}%
	{\cdot\!\supface}%
}}

\renewcommand{\gneq}{\supface}
\renewcommand{\lneq}{\subface}

\newcommand{\mysetminusD}{\hbox{\tikz{\draw[line width=0.6pt,line cap=round] (3pt,0) -- (0,6pt);}}}
\newcommand{\mysetminusT}{\mysetminusD}
\newcommand{\mysetminusS}{\hbox{\tikz{\draw[line width=0.45pt,line cap=round] (2pt,0) -- (0,4pt);}}}
\newcommand{\mysetminusSS}{\hbox{\tikz{\draw[line width=0.4pt,line cap=round] (1.5pt,0) -- (0,3pt);}}}

\newcommand{\mysetminus}{\mathbin{\mathchoice{\mysetminusD}{\mysetminusT}{\mysetminusS}{\mysetminusSS}}}
\renewcommand{\setminus}{\smallsetminus}

\newcommand{\mycontrD}{\hbox{\tikz{\draw[line width=0.6pt,line cap=round] (0,-6pt) -- (4pt,0);}}}
\newcommand{\mycontrT}{\mycontrD}
\newcommand{\mycontrS}{\hbox{\tikz{\draw[line width=0.45pt,line cap=round] (0,-4pt) -- (3pt,0);}}}
\newcommand{\mycontrSS}{\hbox{\tikz{\draw[line width=0.4pt,line cap=round] (0, 1.5pt) -- (3pt,0);}}}

\newcommand{\mycontr}{\mathbin{\mathchoice{\mycontrD}{\mycontrT}{\mycontrS}{\mycontrSS}}}

\newcommand{\ellone}[1]{\ssub{\mleft|#1\mright|}!_{\ell_1}}

\DeclareFontFamily{U}{mathx}{}
\DeclareFontShape{U}{mathx}{m}{n}{<-> mathx10}{}
\DeclareSymbolFont{mathx}{U}{mathx}{m}{n}
\DeclareMathAccent{\widecheck}{0}{mathx}{"71}
\newcommand{\Cir}{\mathscr C}
\newcommand{\CCir}{\widecheck{\mathscr C}}

\newcommand{\x}{\scaleto{X}{6pt}}
\newcommand{\y}{\scaleto{Y}{6pt}}
\newcommand{\ssT}{\ssub[-1.5pt]{T}!}

\numberwithin{equation}{section}

\title{Combinatorial flag arrangements}

\author{Omid Amini}
\address{CNRS - CMLS, \'Ecole polytechnique, Institut polytechnique de Paris}
\email{\href{omid.amini@polytechnique.edu}{omid.amini@polytechnique.edu}}

\author{Lucas Gierczak}
\address{CMLS, \'Ecole polytechnique, Institut polytechnique de Paris}
\email{\href{lucas.gierczak-galle@polytechnique.edu}{lucas.gierczak-galle@polytechnique.edu}}

\date{\today}

\begin{document}

\begin{abstract}
	We introduce combinatorial objects named \emph{matricubes} that provide a generalization of the theory of matroids. As matroids provide a combinatorial axiomatization of hyperplane arrangements, matricubes provide a combinatorial axiomatization of arrangements of initial flags in a vector space. We give cryptomorphic axiomatic systems in terms of rank function, flats, circuits, and independent sets, and formulate a duality concept. We also provide precise links between matricubes, permutation arrays and matroids, and raise several open questions.
\end{abstract}

\maketitle

\setcounter{tocdepth}{1}

\tableofcontents

\section{Introduction}
	
	Consider a vector space $\VS$ of finite dimension over a ground field $\k$ and a collection $\mathcal A$ of $\grdeg$ initial flags $\ssfilt_1^\bullet, \dots, \ssfilt_\grdeg^{\bullet}$. For $j = 1, \dots, \grdeg$, the flag $\ssfilt_j^{\bullet}$ consists of a chain of $\ssr_j + 1$ vector subspaces, $0 \leq \ssr_j \leq \dim_{\k}(\VS)$,
	\[ \VS = \ssfilt_j^0 \supseteq \ssfilt_j^1 \supseteq \dots \supseteq \ssfilt_j^{\ssr_j - 1} \supseteq \ssfilt_j^{\ssr_j} \supsetneq (0) \]
	where, for every $j \in \{1, \dots, \grdeg\}$ and $i \in \{1, \dots, \ssr_j\}$, $\ssfilt^i_j$ is a vector subspace of codimension $0$ or $1$ in $\ssfilt^{i - 1}_j$. We call this collection a \emph{flag arrangement}. In the case $\ssr_j = 1$ for all $j$, and codimension of $\ssfilt^1_j$ in $\VS$ equal to one, we obtain a hyperplane arrangement.
	
	The aim of this paper is to introduce mathematical structures called \emph{matricubes} that provide a combinatorial axiomatization for the intersection patterns of a finite collection of initial flags in a vector space (as the one above). The case $\ssr_1, \dots, \ssr_\grdeg = 1$ recovers the theory of matroids. Like matroids which in the representable case come from matrices, representable matricubes come from cubical matrices (i.e., three-dimensional matrices).
	
	Let us start with a few notations. For $\dvs$ a non-negative integer, we set $[\dvs] \coloneqq \{0, \dots, \dvs\}$.
	Let $\grdeg$ be a positive integer, and $\ssr_1, \dots, \ssr_\grdeg$ be non-negative integers. Let $\undr \coloneqq (\ssr_1, \dots, \ssr_\grdeg)$. The hypercuboid $\hcube{\undr}{}$ of width $\undr$ is the product $\prod_{j = 1}^{\grdeg}[\ssr_j]$. It is endowed with a natural partial order $\subfaceeq$ defined by declaring $\undx \subfaceeq \undy$ for elements $\undx = (\ssx_1, \dots, \ssx_\grdeg)$ and $\undy = (\ssy_1, \dots, \ssy_\grdeg)$ in $\hcube{\undr}{}$, if $\ssx_j \leq \ssy_j$ for all $j$. The minimum and maximum elements of $\hcube{\undr}{}$ are the points $\underline 0 \coloneqq (0, \dots, 0)$ and $\undr$. We define two operations $\vee$ and $\wedge$ of join and meet by taking the maximum and the minimum coordinate-wise, respectively:
	\[ \undx \vee \undy \coloneqq (\max(\ssx_1, \ssy_1), \dots, \max(\ssx_\grdeg, \ssy_\grdeg)), \qquad \undx\wedge \undy \coloneqq (\min(\ssx_1, \ssy_1), \dots, \min(\ssx_\grdeg, \ssy_\grdeg)), \]
	for any pair of elements $\undx, \undy \in \hcube{\undr}{}$. For $i \in \{1, \dots, \grdeg\}$ and $t \in [\ssr_i]$, we denote by $t \, \unde_i$ the point of $\hcube{\undr}{}$ whose $i$-th coordinate is $t$ and whose other coordinates are zero.
	
	We first give the definition of matricubes in terms of their rank functions, and then provide cryptomorphic axiomatic systems in terms of their flats, circuits and independent sets.
	
	\subsection{Definition in terms of rank function}
		
		A function $f\colon \hcube{\undr}{} \rightarrow \Z$ is called \emph{submodular} if for every two elements $\undx$ and $\undy$, we have
		\[ f(\undx) + f\mleft(\undy\mright) \geq f\mleft(\undx \vee \undy\mright) + f\mleft(\undx \wedge \undy\mright). \]
		
		A \emph{matricube} $\matricube$ with ground set $\hcube{\undr}{}$ is defined in terms of a function $\rk \colon \hcube{\undr}{} \rightarrow \Z_{\geq 0}$ called \emph{the rank function of $\matricube$} that satisfies the following conditions:
		
		\begin{enumerate}[label=(R\arabic*)]
			\item \label{def:axiom_rank_function_1_intro} $\rk(\underline 0) = 0$, and for every $1 \leq i \leq \grdeg $ and $1 \leq t \leq \ssr_i$, we have $\rk(t \, \unde_i) - \rk((t - 1) \, \unde_i) \leq 1$.
			
			\item \label{def:axiom_rank_function_2_intro} $\rk$ is non-decreasing, that is, if $\undx \subfaceeq \undy$, then $\rk(\undx) \leq \rk\mleft(\undy\mright)$.
			
			\item \label{def:axiom_rank_function_3_intro} $\rk$ is submodular.
		\end{enumerate}
		
		We call the quantity $r = r(\matricube) \coloneqq \rk\mleft(\undr\mright)$, the maximum value taken by the function $\rk$, the \emph{rank} of $\matricube$. In the case $\ssr_j = 1$ for all $j$, $\matricube$ gives a matroid with ground set $E = \{1, \dots, \grdeg\}$.
		
		Note that it follows from \ref{def:axiom_rank_function_1_intro} and \ref{def:axiom_rank_function_2_intro} that $\rk(t \, \unde_i) \leq t$. We say that $\matricube$ is \emph{simple} if the following stronger version of \ref{def:axiom_rank_function_1_intro} holds:
		
		\begin{enumerate}[label=(R\arabic*$^\ast$)]
			\item \label{def:axiom_rank_function_1nd_intro} $\ssr_i > 0$ and $\rk(t \, \unde_i) = t$ for all $i = 1, \dots, \grdeg$ and $t \in [\ssr_i]$.
		\end{enumerate}
		
		For a collection $\mathcal A$ of initial flags $\ssfilt_1^\bullet, \dots, \ssfilt_\grdeg^{\bullet}$ in a vector space of dimension $\dvs$, as above, the codimensions of the intersection patterns of their elements define a rank function. That is, the function $\rk \colon \hcube{\undr}{} \rightarrow \Z$ defined, for every $\undx = (\ssx_1, \dots, \ssx_\grdeg)$, by
		\[ \rk(\undx) \coloneqq \codim_\k \mleft(\ssfilt_{1}^{\ssx_1} \cap \dots \cap \ssfilt_{\grdeg}^{\ssx_\grdeg}\mright) = \dvs - \dim_\k \mleft(\ssfilt_{1}^{\ssx_1} \cap \dots \cap \ssfilt_{\grdeg}^{\ssx_\grdeg}\mright) \]
		is the rank function of a matricube that we denote by $\matricube_{\mathcal A}$. Note that $\matricube_{\mathcal A}$ is simple if all the inclusions in each flag are strict. Like for matroids, a matricube appearing in this way will be called \emph{representable} over the field $\k$. Note that by duality of vector spaces, a representable matricube can be described equivalently by a collection of initial increasing flags in the dual vector space. This point of view allows to associate a matricube to any three-dimensional matrix with entries in a given field. We refer to Section~\ref{sec:cubical_matrices} for more details.
		
		Abstracting the example given above of an arrangement of initial flags in a vector space, we show in Section~\ref{subsec:flag_matroids} that a finite collection of initial flag matroids, all defined on the same ground set, defines a matricube.
		
		\medskip
		
		In the next three sections, we present alternative axiomatic systems for matricubes, that will be discussed more thoroughly in the paper.
	
	\subsection{Flats of matricubes}
		
		Let $\matricube$ be a matricube of rank $r$ with ground set $\hcube{\undr}{}$ and rank function $\rk$. A point $\unda$ in $\hcube{\undr}{}$ is called a \emph{flat} of $\matricube$ if
		
		\begin{enumerate}
			\item[$(*)$] \label{defi:jump2} for every $i = 1, \dots, \grdeg$ such that $\unda + \unde_i$ belongs to $\hcube{\undr}{}$, we have $\rk(\unda + \unde_i) = \rk(\unda) + 1$.
		\end{enumerate}
	 	
		Note that in particular, $\undr$ is a flat of $\matricube$. We denote by $\Fl(\matricube) \subseteq \hcube{\undr}{}$ the set of flats of $\matricube$. In the case $\matricube$ is a matroid, $\Fl(\matricube)$ is the set of flats of that matroid.
		
		As in the case of matroids, a matricube can be defined in terms of its flats. The axiomatic system of flats of a matricube is~\ref{def:axiom_flats_1_intro}-\ref{def:axiom_flats_2_intro}-\ref{def:axiom_flats_3_intro}, provided below.
	 	
		Given a poset $(P, \subfaceeq)$ and two elements $x, y \in P$, we say that that $y$ \emph{covers} $x$, and write $y \ssupface x$, if $y \supface x$ in $P$ and there is no element $z \in P$ such that $y \supface z \supface x$. Let $\Fl$ be a subset of $\hcube{\undr}{}$. Endowed with the partial order $\subfaceeq$ of the hypercuboid $\hcube{\undr}{}$, $\Fl$ is a poset.
		
		We prove in Section~\ref{sec:flats} that $\Fl \subseteq \hcube{\undr}{}$ is the set of flats of a matricube with underling ground set $\hcube{\undr}{}$ if, and only if, the following properties hold.

		\begin{enumerate}[label=(F\arabic*)]
			\item \label{def:axiom_flats_1_intro} $\undr$ is in $\Fl$.
			
			\item \label{def:axiom_flats_2_intro} $\Fl$ is closed under meet.
			
			\item \label{def:axiom_flats_3_intro} If $\unda$ is an element of $\Fl$ and $i$ is such that $\unda + \unde_i \in \hcube{\undr}{}$, then there exists an element $\undb \supfaceeq \unda + \unde_i$ in $\Fl$ such that $\undb \ssupface \unda$ in $\Fl$.
		\end{enumerate}
		
		In other words, the axiomatic systems \ref{def:axiom_flats_1_intro}-\ref{def:axiom_flats_2_intro}-\ref{def:axiom_flats_3_intro} and \ref{def:axiom_rank_function_1_intro}-\ref{def:axiom_rank_function_2_intro}-\ref{def:axiom_rank_function_3_intro} are equivalent.
	
	\subsection{Duality, and circuits of matricubes}
		
		Again, let $\matricube$ be a matricube on the ground set $\hcube{\undr}{}$. In Section~\ref{subsec:duality}, we define the dual matricube $\matricube^*$ on the same ground set $\hcube{\undr}{}$. In terms of rank functions, the rank function $\rk^*$ of $\matricube^*$ is given by
		\[ \rk^*(\undx) \coloneqq \ellone{\undx} + \rk(\undx^c) - \rk(\matricube) \qquad \forall \, \undx \in \hcube{\undr}{}, \]
		where $\undx^c \coloneqq \undr - \undx$ is the \emph{complement} of $\undx$ in $\hcube{\undr}{}$, and $\rk$ denotes the rank function of $\matricube$.
		
		Denote by $\Fl(\matricube^*)$ the set of flats of the dual matricube, and consider
		\[
			\CCir \coloneqq \mleft\{\unda^c \, \st \, \unda \in \Fl(\matricube^*) \mright\} \subseteq \hcube{\undr}{}.
		\]
		
		A point $\undc$ in $\hcube{\undr}{}$ is called a \emph{circuit} of $\matricube$ if
		
		\begin{enumerate}
			\item[$(*)$] $\undc$ is an element of $\CCir$ which is not the join of any set of elements of $\CCir \setminus \{\undc\}$.
		\end{enumerate}
		
		\smallskip

		We denote by $\Cir(\matricube) \subseteq \hcube{\undr}{}$ the set of circuits of $\matricube$. This definition extends that of circuits in matroids. Moreover, as in the case of matroids, a matricube can be defined in terms of its circuits, via the following axiomatic system. We prove in Section~\ref{sec:circuits} that a subset $\Cir \subseteq \hcube{\undr}{}$ is the set of circuits of a matricube with underlying ground set $\hcube{\undr}{}$ if, and only if, the following properties hold.
		
		\begin{enumerate}[label=(C\arabic*)]
			\item \label{def:axiom_circuits_1_intro} $\underline 0$ is not in $\Fl$.
			
			\item \label{def:axiom_circuits_2_intro} All elements of $\Cir$ are join-irreducible in $\Cir$.
			
			\item \label{def:axiom_circuits_3_intro} If $\unda \in \CCir$ and $i \in \{1, \dots, \grdeg\}$ is such that $\unda - \unde_i \in \hcube{\undr}{}$, then there exists an element $\undb \subfaceeq \unda - \unde_i$ in $\CCir \cup \{\underline 0\}$ such that $\undb \ssubface \unda$ in $\CCir \cup \{\underline 0\}$.
		\end{enumerate}
		
		In other words, the axiomatic systems \ref{def:axiom_circuits_1_intro}-\ref{def:axiom_circuits_2_intro}-\ref{def:axiom_circuits_3_intro} and \ref{def:axiom_rank_function_1_intro}-\ref{def:axiom_rank_function_2_intro}-\ref{def:axiom_rank_function_3_intro} are equivalent.
	
	\subsection{Independents of matricubes}
		
		Let $\matricube$ be a matricube on the ground set $\hcube{\undr}{}$. We say that a point $\unda$ of $\hcube{\undr}{}$ is an \emph{independent of $\matricube$} if
		
		\begin{enumerate}
			\item[$(*)$] for every $i = 1, \dots, \grdeg$ such that $\unda - \unde_i \in \hcube{\undr}{}$, we have $\rk\mleft(\unda - \unde_i\mright) = \rk\mleft(\unda\mright) - 1$.
		\end{enumerate}
	 	We denote by $\Ind(\matricube) \subseteq \hcube{\undr}{}$ the set of independents of $\matricube$.
		
		The set of independents of a matricube is nonempty and closed under meet $\wedge$. Moreover, removing unit vectors from an independent reduces the rank in the following sense: for every independent $\unda \in \Ind(\matricube)$, and every distinct elements $\ssi_1, \dots, \ssi_k \in \{1, \dots, \grdeg\}$ with $\ssa_{\ssi_j} \neq 0$, $1 \leq j \leq k$, we have $\rk\mleft(\unda - \unde_{\ssi_1} - \dots - \unde_{\ssi_k}\mright) = \rk(\unda) - k$.

		We provide an axiomatic system for independent sets of a matricube. In order to do this, we need to define an operation of \emph{removal} in independents.
		
		Let $\Jind$ be a subset of $\hcube{\undr}{}$. Let $\unda$ be an element of $\Jind$ and $i \in \{1, \dots, \grdeg\}$ such that $\ssa_i \neq 0$. If there is at least one element $\undb \subface \unda$ in $\Jind$ which differs from $\unda$ only in the $i$-th component, we define $\unda \setminus i$ to be such an element in $\Jind$ with the largest $i$-th coordinate. In this case, we say that $\unda \setminus i$ is the \emph{removal of $i$ in $\unda$.}
		
		\begin{defi} \label{def:orderable_subset_intro}
			Let $\Jind$ be a subset of $\hcube{\undr}{}$.
				
			\begin{enumerate}[label=(\alph*)]
				\item We say that \emph{removals exist in $\Jind$} if for every $\unda \in \Jind$ and $i \in \{1, \dots, \grdeg\}$, if $\ssa_i \geq 1$, the removal $\unda \setminus i$ exists in $\Jind$.
			\end{enumerate}
			
			If removals exist in $\Jind$, then $\underline 0 \in \Jind$ and for every element $\unda$, there exists a sequence of removals that reduces $\unda$ to $\underline 0$.
			
			\begin{enumerate}[label=(\alph*),resume]	
				\item We say that $\Jind$ is \emph{orderable} if removals exist in $\Jind$ and for every $\unda \in \Jind$, all the sequences of removals that bring $\unda$ to $\underline 0$ have the same length.
			\end{enumerate}
			
			If $\Jind$ is orderable, we define the \emph{size} of $\unda$, denoted by $|\unda|$, as the number of removals needed to reduce $\unda$ to $\underline 0$.
		\end{defi}
		
		The axiomatic system of independents can be formulated as follows. Let $\Ind$ be a subset of $\hcube{\undr}{}$ that verifies the following property:
		
		\begin{enumerate}[label=(I\arabic*)]
			\item \label{def:axiom_independents_1_intro} Removals exist in $\Ind$ and the following holds. For all $\undp \in \Ind$ and removals $\undp \setminus i$ and $\undp \setminus j$, with $i, j \in \{1, \dots, \grdeg\}$, the meet $\undq \coloneqq \mleft(\undp \setminus i\mright) \wedge \mleft(\undp \setminus j\mright)$ belongs to $\Ind$ and, moreover, the two intervals $\mleft[\undq, \undp \setminus i \mright]$ and $\mleft[\undq, \undp \setminus j\mright]$ in $\Ind$ have the same size.
		\end{enumerate}
		
		Here, the interval $[\unda, \undb]$ in $\Ind$ means the set of all $\undc \in \Ind$ which verify $\unda \subfaceeq \undc \subfaceeq \undb$.
		
		We prove in Section~\ref{sec:independents} that \ref{def:axiom_independents_1_intro} is equivalent to the orderability of the set $\Ind$. In particular, if $\Ind$ verifies \ref{def:axiom_independents_1_intro}, we can define the size of $\unda$ as the number of removals needed to reduce $\unda$ to $\underline 0$. This enables to formulate the second property of interest. For $\unda, \undb \in \hcube{\undr}{}$, denote by $D(\unda, \undb)$ the set of elements in $\{1, \dots, \grdeg\}$ such that $\ssa_i < \ssb_i$. The following is understood as a matricube analogue of the augmentation property for independents of matroids.
		
		\begin{enumerate}[label=(I\arabic*),resume]
			\item \label{def:axiom_independents_2_intro} $|\cdot|$ is increasing on independents, i.e., for all $\unda, \undb \in \Ind$ such that $\unda \subface \undb$, we have $|\unda| < |\undb|$. Moreover, let $\unda, \undb$ be two elements of $\Ind$ such that $|\unda| < |\undb|$ and $D(\unda, \undb)$ contains at least two elements. Then, there exists $\undc \in \Ind$ that verifies:
			\begin{itemize}
				\item $\undc \subfaceeq \unda \vee \undb$,
				
				\item $|\undc| > |\unda|$.
				
				\item There exists $i \in D(\unda,\undb)$ such that $\ssc_i < \ssb_i$.
			\end{itemize}
		\end{enumerate}
		
		We prove in Theorem~\ref{thm:axiomatic_independents} that \ref{def:axiom_independents_1_intro}-\ref{def:axiom_independents_2_intro} are equivalent to \ref{def:axiom_rank_function_1_intro}-\ref{def:axiom_rank_function_2_intro}-\ref{def:axiom_rank_function_3_intro}.
	
	\subsection{Permutation arrays}
		
		A combinatorial approach to the study of intersection patterns of a configuration of complete flags was introduced by Eriksson--Linusson in the notion of permutation arrays~\cite{eriksson2000combinatorial, eriksson2000decomposition}. Our Theorem~\ref{thm:equiv_perm_supermodular} proved in Section~\ref{sec:permutation_arrays} shows that permutation arrays are in a one-to-one correspondence with matricubes of rank $r$ or $r + 1$ on the ground set $\hcube{r}{d} \coloneqq \hcube{\undr}{}$ with $\undr = (r, r, \dots, r)$, that is, with all $\ssr_j = r$.
	
	\subsection{Matricubes as coherent complexes of matroids}		
		
		As we explain in Section~\ref{sec:local_matroids}, a matricube locally gives rise to a collection of matroids. Local obstructions for the representability of a matricube can then be formulated in terms of matroid representability. In the case of permutation arrays, via our Theorem~\ref{thm:equiv_perm_supermodular}, this gives obstructions for representability that generalize the examples found in the work of Billey and Vakil~\cite{billey2008intersections}. We moreover go further by proving Theorem~\ref{thm:representability_matricube_matroid}, which shows that the representability of matricubes over infinite fields can be reduced to matroid representability, see Section~\ref{subsec:natural_matroid_representability}.
		
		In Section~\ref{sec:matroidal-characterization}, we provide a matroidal characterization of matricubes by establishing an equivalence between matricubes and coherent complexes of matroids labeled by the elements of a hypercuboid satisfying Properties~\ref{def:axiom_coherent_complex_1_intro} and~\ref{def:axiom_coherent_complex_2_intro} below. Namely, let $\mleft(\mat_{\unda}\mright)_{\unda \in \hcube{\undr}{}}$ be a family of matroids indexed by $\hcube{\undr}{}$, with $\mat_{\unda}$ a matroid on the set $\ssI_{\unda}$ consisting of all $j \in \{1, \dots, \grdeg\}$ with $\ssa_j < \ssr_j$. Denote the rank function of $\mat_{\unda}$ by $\rk_{\unda}$. We say that the collection $\mleft(\mat_{\unda}\mright)$ forms a \emph{coherent complex of matroids} if the following two conditions are met:
		
		\begin{enumerate}[label=(CC\arabic*)]
			\item \label{def:axiom_coherent_complex_1_intro} For all $i \in \{1, \dots, \grdeg\}$ and $0 \leq t < \ssr_i$, we have $\ssrho_{t \unde_i}(\unde_i) \leq 1$.
			
			\item \label{def:axiom_coherent_complex_2_intro} The matroids $\mat_{\unda}$ satisfy the following relation:
			\[ \mat_{\unda + \unde_i} =
			\begin{cases}
				\mat_{\unda} \mycontr i &\textrm{if } \ssa_i = \ssr_i - 1 \\
				\mat_{\unda} \mycontr i \sqcup \{i\} &\textrm{else}
			\end{cases}. \]
		\end{enumerate}
		
		Theorem~\ref{thm:matroidal_characterization} provides an equivalence between \ref{def:axiom_coherent_complex_1_intro}-\ref{def:axiom_coherent_complex_2_intro} and \ref{def:axiom_rank_function_1_intro}-\ref{def:axiom_rank_function_2_intro}-\ref{def:axiom_rank_function_3_intro}.
	
	\subsection{The natural matroid of a matricube and representability} \label{subsec:natural_matroid_representability}
		
		Remembering only the data of the subspaces in a flag arrangement results in a subspace arrangement, whose combinatorics is encoded in an integer polymatroid. In the same way, any matricube on the ground set $\hcube{\undr}{}$, $\undr = (\ssr_1, \dots, \ssr_{\grdeg})$, gives rise to an integer polymatroid on the ground set the disjoint union $[\ssr_1] \sqcup \dots \sqcup [\ssr_\grdeg]$.
		
		Bases and exchange properties for integer polymatroids have been studied by Herzog and Hibi~\cite{HH02}. Csirmaz~\cite{csirmaz2020cyclic} gives axiomatic systems for cyclic flats. In recent work, Bonin, Chun and Fife~\cite{BCF23} study bases, circuits, and cyclic flats in integer polymatroids, connecting them to a classical construction going back to McDiarmid~\cite{McD73}, Lov\'asz~\cite{Lov77}, and Helgason~\cite{Helg06}, which shows that the data of an integer polymatroid on a ground set $E$ is equivalent to the data of a matroid on a larger set $\widehat E$ obtained from $E$ by replacing each element $e$ of $E$ with $\rho(e)$ distinct copies of that element, called the \emph{natural matroid}. (This has gained recent interest in the work~\cite{OSW18}, in developing a decomposition theorem for 2-polymatroids, as well as in the works~\cite{CHLSW22, EL23} related to combinatorial Hodge theory.)
		
		In our setting, starting from a matricube, we can thus associate to it first an integer polymatroid and then use the above construction to replace the integer polymatroid by a matroid on a larger ground set. We will review this construction in Section~\ref{subsec:matricube_polymatroid}. As we explain there, this leads to a story different from the theory presented in this paper.
		
		This point of view is however useful for treating the question of the representability of matricubes. We show in Theorem~\ref{thm:representability_matricube_matroid} that a matricube is representable over an infinite field (or a field of large enough cardinality) if, and only if, the corresponding natural matroid is representable over the same field.
	
	\subsection{Further related work} \label{subsec:further_related_questions}
		
		Our original motivation for developing the theory exposed here comes from the problem of describing tropical degenerations of linear series on algebraic curves. In companion work~\cite{AG22}, matricubes are used as the combinatorial structure underlying a combinatorial theory of limit linear series on metric graphs (the geometric situation behind this theory is briefly discussed in Section~\ref{subsec:geometric_rank_functions}). While working on the degeneration problem for linear series, we gradually realized how similarly matricubes and matroids behave. Apart from bases, for which we do not provide a definition and an axiomatic system, the other relevant constructions in the theory of matroids have their matricube analogues.
		
		The recent work of Baker and Bowler~\cite{BB19} develops a theory of matroids over hyperstructures. The extension of this theory to flag matroids is given by Jarra and Lorscheid~\cite{JL24}, and a generalization to quiver matroids is the subject of a forthcoming work of Jarra, Lorscheid and Vital. The work by Baker and Lorscheid~\cite{BL21, BL20} studies the moduli space of matroids and deduce applications to representability questions for matroids. It seems plausible and interesting to generalize these results to the context of matricubes.
		
		In~\cite{bollen2018algebraic}, Bollen, Draisma and Pendavingh show that each representation of an algebraic matroid $\mat$ over a field of positive characteristic comes naturally with a valuation, that they name the \emph{Lindström valuation} of that representation. To this end, using the Frobenius map of the base field, they associate to any such representation what they call a \emph{matroid flock}, an infinite family of linear matroids of the same rank as $\mat$, indexed by $\Z^E$, where $E$ is the ground set of $\mat$. It is interesting to note that, although these notions arise in totally different contexts, the axiomatic systems of coherent complexes and matroid flocks are reminiscent of each other. There are however some major differences. Namely, the matroids appearing in a matroid flock all have the same rank, and there is an invariance property with respect to the direction $(1, \dots, 1)$. Besides, the boundary condition imposed on the coherent complexes does not appear in matroid flocks.
		
		Submodular functions on distributive lattices are a central topic in the study of a large class of combinatorial optimization problems. We refer to the books by Schrijver~\cite{schrijver2003combinatorial} and Fuji~\cite{fujishige2005submodular} for a discussion of these aspects.
		
		Murota~\cite{murota1998discrete} investigates a theory of convex analysis in the discrete setting that involves functions $f \colon \Z^n \to \Z$. Classical duality theorems about real convex functions are proved in the discrete setting. Discrete convexity in that setting is similar in spirit to the submodularity property studied in the present paper.
	
	\subsection{Organization of the text}
		
		In Section~\ref{sec:basic_properties}, we define matricubes using rank functions and give basic examples, including uniform and representable matricubes. We define operations of deletion and contraction on matricubes, and formulate a duality concept.
		
		In Sections~\ref{sec:flats}, \ref{sec:circuits} and~\ref{sec:independents}, we explore alternative axiomatic systems for matricubes, relying on flats, circuits and independents, respectively.
		
		In Section~\ref{sec:diamond_property}, we prove elementary combinatorial results useful throughout the paper, which provide a simpler way of checking whether a function on a hypercuboid is a rank function.
		
		In Section~\ref{sec:permutation_arrays}, we show that particular kinds of matricubes are in a natural one-to-one correspondence with permutation arrays.
		
		In Section~\ref{sec:local_matroids}, we provide the equivalence of matricubes with coherent complexes of matroids, and provide local obstructions for representability.
		
		Finally, in Section~\ref{sec:further_discussions}, we discuss further interesting features of matricubes and raise several open questions.
	
	\subsection{Acknowledgments}
		
		The authors thank Matt Baker, Alex Fink, St\'ephane Gaubert, June Huh, Oliver Lorscheid, Matthieu Piquerez and Farbod Shokrieh for discussions related to the content of this paper. Special thanks to Alex Fink for making us aware of references~\cite{bollen2018algebraic, BCF23, EL23}. O.\thinspace A. thanks Math+, the Berlin Mathematics Research Center and TU Berlin, where part of this research was carried out. L.\thinspace G. thanks the Mathematical Foundation Jacques Hadamard (FMJH) for financial support.

\section{Basic properties} \label{sec:basic_properties}
	
	Let $\dvs$ be a non-negative integer and $[\dvs] = \{0, 1, \dots, \dvs\}$. For elements $\ssr_1, \dots, \ssr_\grdeg \in[\dvs]$, the \emph{hypercuboid} $\hcube{\undr}{}$ of \emph{width} $\undr = (\ssr_1, \dots, \ssr_\grdeg)$ is the product $\prod_{j = 1}^\grdeg [\ssr_j]$. When $\ssr_1 = \cdots = \ssr_d = r$, we simply denote the hypercuboid by $\hcube{r}{\grdeg}$. We denote the elements of $\hcube{\undr}{}$ by vectors $\undx = (\ssx_1, \dots, \ssx_\grdeg)$, for $\ssx_1 \in [\ssr_1], \dots, \ssx_\grdeg \in [\ssr_\grdeg]$. In the hypercuboid, we define, for every $i = 1, \dots, \grdeg$ and $t \in [\ssr_i]$, the $t$-th layer in the direction $i$ as $\ssL^i_t \coloneqq \mleft\{\undx \in \hcube{\undr}{}, \ssx_i = t\mright\}$.
	
	We endow $\hcube{\undr}{}$ with the partial order $\subfaceeq$: For a pair of elements $\undx, \undy \in \hcube{\undr}{}$, we have $\undx \subfaceeq \undy$ provided that $\ssx_j \leq \ssy_j$ for all $j = 1, \dots, \grdeg$. The smallest and largest elements with respect to $\subfaceeq$ are $\underline 0$ and $\undr$, respectively. Moreover, there is a lattice structure on $\hcube{\undr}{}$, where the two operations of join $\vee$ and meet $\wedge$ correspond to taking the maximum and the minimum coordinate-wise, respectively.
	
	A function $f \colon \hcube{\undr}{} \rightarrow \Z$ is called \emph{submodular} if for every pair of elements $\undx$ and $\undy$, we have
	\[ f(\undx) + f\mleft(\undy\mright) \geq f\mleft(\undx \vee \undy\mright) + f\mleft(\undx \wedge \undy\mright). \]
	We will be interested in a special kind of submodular function on $\hcube{\undr}{}$. For each integer $i \in \{1, \dots, \grdeg\}$, we denote by $\unde_i$ the vector whose coordinates are all zero except the $i$-th coordinate, which is equal to one. For $0 \leq t \leq \ssr_i$, the vector $t \, \unde_i$ lies in $\hcube{\undr}{}$.
	
	\begin{defi}[Matricube] \label{def:matricube}
		A \emph{matricube} $\matricube$ with ground set $\hcube{\undr}{}$ is defined in terms of a function $\rk \colon \hcube{\undr}{} \rightarrow \Z$ called \emph{the rank function of $\matricube$} that satisfies the following conditions:
		
		\begin{enumerate}[label=(R\arabic*)]
			\item \label{def:axiom_rank_function_1} $\rk(\underline 0) = 0$, and for every $1 \leq i \leq \grdeg $ and $1 \leq t \leq \ssr_i$, we have $\rk(t \, \unde_i) - \rk((t - 1) \, \unde_i) \leq 1$.
			
			\item \label{def:axiom_rank_function_2} $\rk$ is non-decreasing with respect to $\subfaceeq$, that is, if $\unda \subfaceeq \undb$, then $\rk(\unda) \leq \rk(\undb)$.
			
			\item \label{def:axiom_rank_function_3} $\rk$ is submodular.
		\end{enumerate}
		
		We call $r = r(\matricube) \coloneqq \rk\mleft(\undr\mright)$, the maximum value taken by the function $\rk$, the \emph{rank} of $\matricube$.
	\end{defi}
	
	Note that \ref{def:axiom_rank_function_1} implies that $\rk(t \, \unde_i) \leq t$ for all $i = 1, \dots, \grdeg$ and $t \in [\ssr_i]$. We say that $\matricube$ is \emph{simple} if the following alternate form of \ref{def:axiom_rank_function_1} holds:
	\begin{enumerate}[label=(R\arabic*$^\ast$)]
		\item \label{def:axiom_rank_function_1nd} $\ssr_i > 0$ and $\rk(t \, \unde_i) = t$ for all $i = 1, \dots, \grdeg$ and $t \in [\ssr_i]$.
	\end{enumerate}
	
	\begin{remark} \label{rem:link_coordinate_rank}
		In $\matricube$ is simple, then the above properties imply that if $\undx \in \hcube{\undr}{}$ has rank $j$, then $\ssx_i \leq j$ for all $i = 1, \dots, \grdeg$. In particular, $\underline 0$ is the only element of rank $0$ in $\hcube{\undr}{}$.
	\end{remark}
	
	To be able to present examples of rank functions easily, we adopt the following convention.
	
	\begin{conv}[Cases $\grdeg = 1, 2, 3$]
		In this article, for $\grdeg = 1$, a function on $\hcube{r}{}$ is described by a tuple with $r + 1$ entries $(\sst_0, \dots, \sst_r)$, which means that the value of the function on the $i$-th entry of $\hcube{r}{}$ is $\sst_i$.
		
		In the same way, for $\grdeg = 2$, a function on $\hcube{(\ssr_1, \ssr_2)}{}$ will often be described by an array of size $(\ssr_1 + 1) \times (\ssr_2 + 1)$, $\ssub{(\sst_{ij})}!_{0 \leq i \leq \ssr_1, 0 \leq j \leq \ssr_2}$, which means that the function takes value $\sst_{ij}$ on $(i, j) \in \hcube{(\ssr_1, \ssr_2)}{}$. We choose the convention that the first direction is horizontal, the second direction is vertical, and the origin is the bottom left-hand corner.
		
		When $\grdeg = 3$ and $\undr = (\ssr_1, \ssr_2, \ssr_3)$, a function defined on $\hcube{\undr}{}$ will be specified by $\ssr_3 + 1$ arrays $\ssR_0, \dots, \ssR_{\ssr_3}$ of size $(\ssr_1 + 1) \times (\ssr_2 + 1)$, where $\ssR_k$ describes the values of the function on $\hcube{(\ssr_1, \ssr_2)}{} \times \{k\} \subseteq \hcube{\undr}{}$.
	\end{conv}
	
	Here are two examples of matricubes with $\undr = (4, 3)$.
	\[
		\begin{pmatrix}
			3 & 3 & 3 & 4 & 5 \\
			2 & 2 & 2 & 3 & 4 \\
			1 & 2 & 2 & 3 & 4 \\
			0 & 1 & 2 & 3 & 4
		\end{pmatrix}
		\qquad \qquad
		\begin{pmatrix}
			3 & 3 & 3 & 3 & 4 \\
			2 & 2 & 2 & 2 & 3 \\
			1 & 2 & 2 & 2 & 3 \\
			0 & 1 & 2 & 2 & 3
		\end{pmatrix}
	\]
	The one on the left is simple, the one on the right is not.
	
	For future use, we state the following proposition which implies that the set of values of the rank function in a matricube of rank $r$ is the interval $[r]$.
	
	\begin{prop} \label{prop:inequality_one}
		Let $\rk$ be a rank function on $\hcube{\undr}{}$. Let $i \in \{1, \dots, \grdeg\}$. For an element $\undx \in \hcube{\undr}{} $ such that $\undx + \unde_i \in \hcube{\undr}{}$, we have
		\[ \rk(\undx) \leq \rk(\undx + \unde_i) \leq \rk(\undx) + 1. \]
	\end{prop}
	
	\begin{proof}
		Let $\undy = (\ssx_i + 1) \, \unde_i$, and note that $\undx \vee \undy = \undx + \unde_i$ and $\undx \wedge \undy = \ssx_i \, \unde_i$. Applying the submodularity of $\rk$ to the vectors $\undx$ and $\undy$, and using \ref{def:axiom_rank_function_1} in Definition~\ref{def:matricube}, we get $\rk(\undx) + 1 \geq \rk(\undx + \unde_i)$. The first inequality follows from the non-decreasing property of $\rk$.
	\end{proof}
	
	\subsection{Uniform matricubes}
		
		Notation as in the previous section, let $\undr = (\ssr_1, \dots, \ssr_{\grdeg})$, and consider the corresponding hypercuboid $\hcube{\undr}{}$. Let $r \in [\ssr_1 + \dots + \ssr_{\grdeg}]$ be a non-negative integer. We define the uniform matricube $\unimatricube_{\undr, r}$ of width $\undr$ and rank $r$ as the matricube defined by the \emph{standard rank function} defined as follows
		\[ \rkstandard(\unda) \coloneqq \min(r, \ssa_1 + \cdots + \ssa_\grdeg) \qquad \textrm{ for } \unda = (\ssa_1, \dots, \ssa_\grdeg) \in \hcube{\undr}{}. \qedhere \]
		Notice that the uniform matricube $\unimatricube_{\undr, r}$ is simple if, and only if, $r \geq \max_i \ssr_i$.
		
		Below are the uniform matricubes $\unimatricube_{(4, 3), 3}$ and $\unimatricube_{(4, 3), 5}$.
		\[
			\begin{pmatrix}
				3 & 3 & 3 & 3 & 3 \\
				2 & 3 & 3 & 3 & 3 \\
				1 & 2 & 3 & 3 & 3 \\
				0 & 1 & 2 & 3 & 3
			\end{pmatrix}
			\qquad \qquad
			\begin{pmatrix}
				3 & 4 & 5 & 5 & 5 \\
				2 & 3 & 4 & 5 & 5 \\
				1 & 2 & 3 & 4 & 5 \\
				0 & 1 & 2 & 3 & 4
			\end{pmatrix}
		\]
		
		\begin{prop}
			The rank function $\rk$ of any matricube $\matricube$ of rank $r$ on the ground set $\hcube{\undr}{}$ is dominated by the rank function $\rkstandard$ of the uniform matricube $\unimatricube_{\undr, r}$. That is, for every $\undx \in \hcube{\undr}{}$, we have $\rk(\undx) \leq \rkstandard(\undx)$.
		\end{prop}
		
		\begin{proof}
			It follows directly from Proposition~\ref{prop:inequality_one} that we have $\rk(\undx) \leq \ssx_1 + \cdots + \ssx_\grdeg$. Combined with $\rk(\undx) \leq r$, we deduce the result.
		\end{proof}
	
	\subsection{Representable matricubes} \label{subsec:representable_matricubes}
		
		Let $\dvs$ be a non-negative integer, and let $\VS$ be a vector space of dimension $\dvs$ over some field $\k$. An \emph{initial (decreasing) flag} of $\VS$ of length $r$ consists of a chain of vector subspaces
		\[ \VS = \ssfilt^0 \supseteq \ssfilt^1 \supseteq \dots \supseteq \ssfilt^{r - 1} \supseteq \ssfilt^r \supsetneq (0), \]
		where for each positive $i \in [r]$, $\filt^i$ is a vector subspace of codimension 0 or 1 in $\filt^{i - 1}$. We say that $\filt^\bullet$ is \emph{simple} if each $\filt^{i}$ has codimension $i$ in $\VS$. A \emph{complete flag} is a simple initial flag of length $\dvs - 1$.
		
		\smallskip
		
		Let $\grdeg$ be a positive integer, and let $\mathcal A$ be a collection of $\grdeg$ initial flags $\ssfilt_1^\bullet, \dots, \ssfilt_\grdeg^{\bullet}$ of $\VS$ of lengths $\ssr_1, \dots, \ssr_\grdeg$, respectively. Define the function $\rk \colon \hcube{\undr}{} \rightarrow \Z$ by
		\begin{equation} \label{eq:def_rank_function}
			\rk(\undx) \coloneqq \codim_\k \mleft(\ssfilt_{1}^{\ssx_1} \cap \dots \cap \ssfilt_{\grdeg}^{\ssx_\grdeg}\mright) \quad \forall \, \undx = (\ssx_1, \dots, \ssx_\grdeg) \in \hcube{\undr}{}.
		\end{equation}
		
		\begin{prop}
			The hypercuboid $\hcube{\undr}{}$ endowed with the function $\rk$ defined in~\eqref{eq:def_rank_function} is a matricube. This matricube is simple if, and only if, all the initial flags are simple.
		\end{prop}
		
		\begin{proof}
			Let $\unda$ and $\undb$ be two points of $\hcube{\undr}{}$, and let $\undx \coloneqq \unda \wedge \undb$ and $\undy \coloneqq \unda \vee \undb$. We have an injection
			\[ \rquot{\bigl(\ssfilt_1^{\ssa_1} \cap \dots \cap \ssfilt_{\grdeg}^{\ssa_\grdeg}\bigr)}{\bigl(\ssfilt_1^{\ssy_1} \cap \dots \cap \ssfilt_\grdeg^{\ssy_\grdeg}\bigr)} \hookrightarrow \rquot{\bigl(\ssfilt_1^{\ssx_1} \cap \dots \cap \ssfilt_\grdeg^{\ssx_\grdeg}\bigr)}{\bigl(\ssfilt_1^{\ssb_1} \cap \dots \cap \ssfilt_\grdeg^{\ssb_\grdeg}\bigr)}, \]
			from which, comparing the dimensions, we get $\rk(\undb) - \rk(\undx) \geq \rk\mleft(\undy\mright) - \rk(\unda)$. This proves the submodularity of $\rk$. Properties~\ref{def:axiom_rank_function_1} and~\ref{def:axiom_rank_function_2} in Definition \ref{def:matricube} are trivially verified. This proves the first assertion. The matricube is simple if, and only if, each $\ssfilt^i_j$ has codimension $i$ in $\VS$, that is if, and only if, $\ssfilt^\bullet_j$ is simple, for $j = 1, \dots, \grdeg$.
		\end{proof}
		
		We denote by $\matricube_{\mathcal A}$ the matricube associated to $\mathcal A$.
		
		\begin{defi}[Representable matricube] \label{def:representable_matricube}
			A matricube $\matricube$ on ground set $\hcube{\undr}{}$ is called \emph{representable over a field $\k$} if it is the matricube associated to an arrangement of $\grdeg$ initial flags $\ssfilt_1^\bullet, \dots, \ssfilt_\grdeg^\bullet$ of lengths $\ssr_1, \dots, \ssr_{\grdeg}$, respectively, in a $\k$-vector space $\VS$.
		\end{defi}
		
		\begin{example}
			The matricube $\unimatricube_{\undr, r}$ is representable over every field of large enough cardinality. Indeed, it is the matricube associated to an arrangement of $\grdeg$ initial flags of lengths $\ssr_1, \dots, \ssr_\grdeg$ in $\VS$ of dimension $r$ which are in \emph{general relative position}, that is, whose intersection patterns have the smallest possible dimensions.
		\end{example}
		
		Let $\matricube$ be a matricube associated to an arrangement of $\grdeg$ initial (decreasing) flags $\ssfilt_1^\bullet, \dots, \ssfilt_\grdeg^\bullet$ inside $\VS$. For every $i \in \{1, \dots, \grdeg\}$, duality transforms the initial (decreasing) flag in $\VS$
		\[ \VS = \ssfilt_i^0 \supseteq \ssfilt_i^1 \supseteq \cdots \supseteq \ssfilt_i^{\ssr_i - 1} \supseteq \ssfilt_i^{\ssr_i} \supsetneq (0) \]
		into an \emph{initial (increasing) flag} in $\VS^*$
		\[ (0) = \ssfiltg^i_0 \subseteq \ssfiltg^i_1 \subseteq \cdots \subseteq \ssfiltg^i_{\ssr_i - 1} \subseteq \ssfiltg^i_{\ssr_i} \subsetneq \VS^*, \]
		where $\ssfiltg^i_j$ is the orthogonal to $\ssfilt_i^j$ for the duality pairing $\langle \cdot \, , \cdot \rangle \colon \VS \times \VS^* \to \k$, that is, $\ssfiltg_j^i \coloneqq \bigcap_{v \in \ssfilt_i^j} \ker \mleft(v \colon \VS^* \to \k\mright)$ and $\ssfilt_i^j \coloneqq \bigcap_{\ell \in \ssfiltg^i_j} \ker \mleft(\ell \colon \VS \to \k\mright)$, so that one filtration can be recovered from the other. Note that in the case of matroids, this duality corresponds to the one between arrangements of hyperplanes in $\VS$ and arrangements of vectors in the dual $\VS^*$.
		
		The rank function $\rk$ of $\matricube$, defined in Equation~\eqref{eq:def_rank_function} using intersections of the $\ssfilt_i^j$, can be alternatively described using the flags $\ssfiltg_j^i$ in the following way:
		\begin{align}
			\label{eq:rank_dual_space} \rk(\undx) = \dim_\k \mleft(\ssfiltg^1_{\ssx_1} + \cdots + \ssfiltg^\grdeg_{\ssx_\grdeg}\mright) \qquad \qquad \forall \, \, \undx = (\ssx_1, \dots, \ssx_\grdeg) \in \hcube{\undr}{}.
		\end{align}
		
		We will discuss the representability of matricubes further in Section~\ref{subsec:disc_minor}.
		
		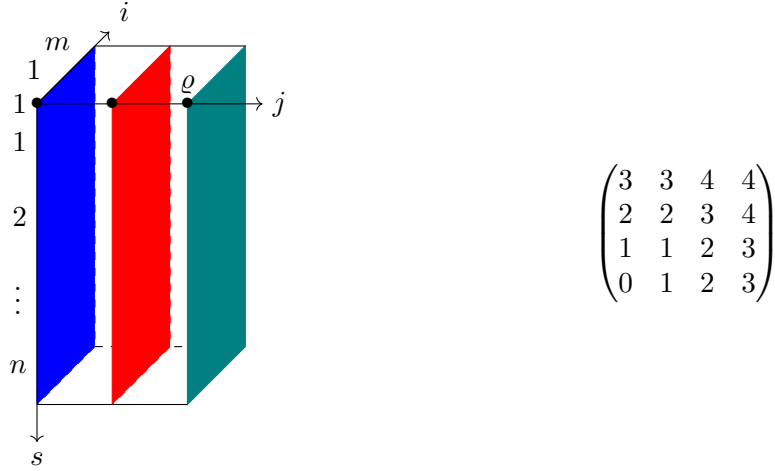
\begin{figure}[h!]
			\begin{minipage}{0.45\textwidth}
			\centering
			\begin{tikzpicture}[scale=1]
				% Background edge.
				\draw[dashed] (0,0,0) -- (2,0,0);
				% Blue layer.
				\begin{scope}[color=lightgray]
					\draw (0,3,2) -- (0,3,0);
					\draw (0,2,2) -- (0,2,0);
					\draw (0,1,2) -- (0,1,0);
				\end{scope}
				\begin{scope}[color=blue]
					\draw[dashed] (0,4,0) -- (0,0,0) -- (0,0,2);
					\draw (0,0,2) -- (0,4,2) -- (0,4,0);
					\draw[dashed,line width=0.05mm] (0,4,1) -- (0,0,1);
					\fill[opacity=0.1] (0,0,2) -- (0,4,2) -- (0,4,0) -- (0,0,0) -- cycle;
					\begin{scope}[canvas is yz plane at x=0,transform shape]
						\node[rotate=-90] at (3.5,1.5) {$1$};
						\node[rotate=-90] at (2.5,1.5) {$0$};
						\node[rotate=-90] at (1.5,1.5) {$0$};
						\node[rotate=-90] at (0.5,1.5) {$0$};
						\node[rotate=-90] at (3.5,0.5) {$\scaleto{\pi}{7pt}$};
						\node[rotate=-90] at (2.5,0.5) {$0$};
						\node[rotate=-90] at (1.5,0.5) {$0$};
						\node[rotate=-90] at (0.5,0.5) {$0$};
					\end{scope}
				\end{scope}
				% Red layer.
				\begin{scope}[color=lightgray]
					\draw (1,3,2) -- (1,3,0);
					\draw (1,2,2) -- (1,2,0);
					\draw (1,1,2) -- (1,1,0);
				\end{scope}
				\begin{scope}[color=red]
					\draw[dashed] (1,4,0) -- (1,0,0) -- (1,0,2);
					\draw (1,0,2) -- (1,4,2) -- (1,4,0);
					\draw[dashed,line width=0.05mm] (1,4,1) -- (1,0,1);
					\fill[opacity=0.1] (1,0,2) -- (1,4,2) -- (1,4,0) -- (1,0,0) -- cycle;
					\begin{scope}[canvas is yz plane at x=1,transform shape]
						\node[rotate=-90] at (3.5,1.5) {$0$};
						\node[rotate=-90] at (2.5,1.5) {$1$};
						\node[rotate=-90] at (1.5,1.5) {$0$};
						\node[rotate=-90] at (0.5,1.5) {$0$};
						\node[rotate=-90] at (3.5,0.5) {$0$};
						\node[rotate=-90] at (2.5,0.5) {$0$};
						\node[rotate=-90] at (1.5,0.5) {$0$};
						\node[rotate=-90] at (0.5,0.5) {$1$};
					\end{scope}
				\end{scope}
				% Green layer.
				\begin{scope}[color=lightgray]
					\draw (2,3,2) -- (2,3,0);
					\draw (2,2,2) -- (2,2,0);
					\draw (2,1,2) -- (2,1,0);
				\end{scope}
				\begin{scope}[color=teal]
					\draw (2,0,2) -- (2,4,2) -- (2,4,0) -- (2,0,0) -- cycle;
					\draw[line width=0.05mm] (2,4,1) -- (2,0,1);
					\fill[opacity=0.1] (2,0,2) -- (2,4,2) -- (2,4,0) -- (2,0,0) -- cycle;
					\begin{scope}[canvas is yz plane at x=2,transform shape]
						\node[rotate=-90] at (3.5,1.5) {$0$};
						\node[rotate=-90] at (2.5,1.5) {$0$};
						\node[rotate=-90] at (1.5,1.5) {$1$};
						\node[rotate=-90] at (0.5,1.5) {$0$};
						\node[rotate=-90] at (3.5,0.5) {$1$};
						\node[rotate=-90] at (2.5,0.5) {$0$};
						\node[rotate=-90] at (1.5,0.5) {$0.5$};
						\node[rotate=-90] at (0.5,0.5) {$0$};
					\end{scope}
				\end{scope}
				% Foreground edges.
				\draw (0,0,2) -- (2,0,2);
				\draw (0,4,0) -- (2,4,0);
				% Axes
				\draw[->] (0,4,2) -- (3,4,2);
				\draw[->] (0,4,2) -- (0,-0.5,2);
				\draw[->] (0,4,2) -- (0,4,-0.5);
				\draw (3,4,2) node[right]{$j$};
				\draw (0,-0.5,2) node[below]{$s$};
				\draw (0,4,-0.5) node[above right]{$i$};
				\draw (0,4,2) node[left]{$1$} node{$\bullet$};
				\draw (1,4,2) node{$\bullet$};
				\draw (2,4,2) node[above]{$\ssr$} node{$\bullet$};
				\draw (0,4,1.5) node[above left]{$1$};
				\draw (0,4,0.5) node[above left]{$m$};
				\draw (0,3.5,2) node[left]{$1$};
				\draw (0,2.5,2) node[left]{$2$};
				\draw (0,1.5,2) node[left]{$\vdots \,$};
				\draw (0,0.5,2) node[left]{$\dvs$};
			\end{tikzpicture}
			\end{minipage}
			\begin{minipage}{0.45\textwidth}
			\begin{center}
				$\begin{pmatrix}
					3 & 3 & 4 & 4 \\
					2 & 2 & 3 & 4 \\
					1 & 1 & 2 & 3 \\
					0 & 1 & 2 & 3
				\end{pmatrix}$
			\end{center}
			\end{minipage}
			\caption{The left figure represents a three-dimensional matrix $A = (A^i_{js})$ of size $\grdeg \times \ssr \times \dvs$ with $\grdeg = 2$, $\ssr = 3$ and $\dvs = 4$. The blue (resp. red, resp. green) layer contains vertically the coordinates of the vectors $\ssv^1_1$ and $\ssv^2_1$ (resp. $\ssv^1_2$ and $\ssv^2_2$, resp. $\ssv^1_3$ and $\ssv^2_3$). The associated matricube is given on the right.}
			\label{fig:matricube_matrix}
		\end{figure}
	
	\subsection{Matricube induced by a cubical matrix} \label{sec:cubical_matrices}
		
		Using the duality between initial (decreasing) flags in $\VS$ and initial (increasing) flags in $\VS^*$, we explain a procedure that associates a representable matricube to any three-dimensional matrix with coordinates in a field $\k$. This construction extends the representation of representable matroids by matrices. As in the case of matroids, this justifies the terminology \emph{matricube}, which encompasses both the idea of \emph{cubical matrix} (like ``matroid'', coming from ``matrix'') and \emph{hypercuboid} (a matricube is described by a hypercuboid of numbers, given by the rank function).
		
		Notation as in the previous section, first assume that $\ssr_1 = \cdots = \ssr_\grdeg = \ssr$ and let $i \in \{1, \dots, \grdeg\}$. We choose vectors $\ssv^i_1, \dots, \ssv^i_\ssr \in \VS^*$ such that for every $j \in \{1, \dots, \ssr\}$, $\ssfiltg^i_j = \mleft\langle \ssv^i_1, \dots, \ssv^i_j \mright\rangle$. This shows that a representable matricube $\matricube$ can be determined by the collection of vectors $\ssv^i_j$, for $i \in \{1, \dots, \grdeg\}$ and $j \in \{1, \dots, \ssr\}$. Said otherwise, choosing a basis of $\VS^*$, $\matricube$ is determined by a three-dimensional matrix $A = (A^i_{js})$ where $i \in \{1, \dots, \grdeg\}$, $j \in \{1, \dots, \ssr\}$ and $s \in \{1, \dots, \dvs\}$, $\dvs$ being the dimension of $\VS^*$. Inversely, using the definition of rank function given in \eqref{eq:rank_dual_space}, this procedure gives a way to associate to every three-dimensional matrix $A = (A^i_{js})$ of size $\grdeg \times \ssr \times \dvs$ with entries in a field $\k$ a matricube $\matricube_A$ on the hypercube $\hcube{\ssr}{\grdeg}$.
		
		In the general case, if not all $\ssr_i$ are equal, we set $\ssr \coloneqq \max_i \ssr_i$ and choose, for every $i \in \{1, \dots, \grdeg\}$, a family of vectors $\ssv^i_1, \dots, \ssv^i_{\ssr_i}, 0, \dots, 0$, with $\ssv^i_j$ for $j \leq \ssr_i$ as before, completed now with $\ssr - \ssr_i$ copies of the zero vector. This gives a matrix $A$ of size $\grdeg \times \ssr \times \dvs$. In the matricube $\matricube_A$ associated to $A$, we now delete, for every $i \in \{1, \dots, \grdeg\}$, $\ssr - \ssr_i$ times the element $i$ (we refer to Section~\ref{subsubsec:deletion} below for the definition of the operation of deletion). This gives the matricube associated to the original family of flags.	
		
		It follows from the construction above that every representable matricube is associated to some three-dimensional matrix, possibly after a few deletions corresponding to zero vectors.
		
		An example of a matricube associated to a three-dimensional real matrix with $\grdeg = 2$, $\ssr_1 = \ssr_2 = 3$ and $\dvs = 4$, is depicted in Figure~\ref	{fig:matricube_matrix}. The ground set of the corresponding matricube is $[3] \times [3]$.
	
	\subsection{Matricube induced by a collection of flag matroids} \label{subsec:flag_matroids}
		
		We show that a finite collection of initial flag matroids on the same ground set $E$ gives rise to a (simple) matricube. We refer to~\cite[Chapter 1]{BGW03} and \cite{CDMS17} for a nice introduction to flag matroids.
	 	
		Let $E$ be a finite set and $\ssr$ be a positive integer. An \emph{initial flag} of size $s$ is an increasing chain of subsets $\ssF_0 = \emptyset \subsetneq \ssF_1 \subsetneq \ssF_2 \subsetneq \dots \subsetneq \ssF_{s} \subseteq E$ with $|\ssF_j| = j$ for $j = 1, \dots, s$. Note that the data of an initial flag is equivalent to an ordered sequence $\ssvare_1, \dots, \ssvare_s$ of distinct elements of $E$, $\ssF_j$ consisting of the first $j$ elements $\ssvare_1, \dots, \ssvare_j$.
		
		A total order $<_{\order}$ on $E$ induces a partial, element-wise order on $E^s$. Through the bijection between initial flags of size $s$ and ordered sequences of size $s$ in $E$, $<_{\order}$ induces a partial order on initial flags of size $s$.
		
		An \emph{initial flag matroid} $\mat$ of rank $s$ is a collection $\mathscr F$ of initial flags of size $s$ as above such that for any total order $<_{\order}$ on $E$, there exists a unique flag in $\mathscr F$ maximal with respect to the induced partial order on initial flags of size $s$. In this case, the following properties hold:
		
		\begin{itemize}
			\item the collection consisting of the terms $\ssF_j$ of flags $\ssF_{\bullet}$ appearing in $\mathscr F$ forms the set of bases of a matroid $\mat_j$ of rank $j$ on the ground set $E$, for $j = 0, \dots, s$;
			
			\item the matroid $\mat_{j}$ is a quotient of the matroid $\mat_{j + 1}$;
			
			\item any sequence $\ssF_0 = \emptyset \subset \ssF_1 \subset \ssF_2 \subset \dots \subset \ssF_s$ with $\ssF_j$ a basis of $\mat_j$ is an element of $\mathscr F$.
		\end{itemize}
		
		(These properties are equivalent to $\mathscr F$ defining a flag matroid, see~\cite[Theorem 1.7.1]{BGW03}.) Elements of $\mathscr F$ are called \emph{bases}, and $\mat_j$ is called the \emph{$j$-th constituent} of $\mat$.
		
		Let now $\undr = (\ssr_1, \dots, \ssr_{\grdeg})$ be a vector with positive integer entries. Consider a collection $\mat^1, \dots, \mat^{\grdeg}$ of initial flag matroids of respective ranks $\ssr_1, \dots, \ssr_{\grdeg}$ on the ground set $E$. For each $i = 1, \dots, \grdeg$, and each $j \in [\ssr_i]$, denote by $\mat^i_j$ the $j$-th constituent of $\mat^i$. This is a matroid of rank $j$.
		
		For each $\undx \in \hcube{\undr}{}$, denote by $\mat_{\undx}$ the matroid union $\mat_{x_1}^1 \cup \dots \cup \mat_{x_\grdeg}^{\grdeg}$ of $\mat_{x_1}^1, \dots, \mat_{x_\grdeg}^{\grdeg}$. Recall that the independent sets of the matroid union $\mat_\undx$ are subsets of $E$ of the form $\ssI_1 \cup \dots \cup \ssI_\grdeg$ where each $\ssI_i$ is an independent of the matroid $\mat_{x_i}^i$ for $i = 1, \dots, \grdeg$.
		
		Consider the function $\rk \colon \hcube{\undr}{} \to \Z_{\geq 0}$ on the hypercuboid defined by
		\[
			\rk(\undx) \coloneqq r(\mat_{\undx}) \qquad \qquad \forall \undx \in \hcube{\undr}{},
		\]
		where $r(\mat_\undx)$ is the rank of the matroid $\mat_{\undx}$.
		
		\begin{thm}
			Notation as above, $\rk$ is the rank function of a simple matricube $\matricube = \matricube(\mat^1, \dots, \mat^\grdeg)$ with ground set the hypercuboid $\hcube{\undr}{}$.
		\end{thm}
		
		\begin{proof}
			We first note that $\rk(t \, \unde_i) = r(\mat^i_t) = t$ for all $t \in [\ssr_i]$. This shows that \ref{def:axiom_rank_function_1nd} is verified. The axiom~\ref{def:axiom_rank_function_2} is obviously verified by the definition of the matroid union. It thus remains to show \ref{def:axiom_rank_function_3}, i.e., that $\rk$ is submodular. By Theorem~\ref{thm:eq_submodularity_diamond_property}, it will be enough to show that $\rk$ verifies the diamond property, that is, for all $\undx \in \hcube{\undr}{}$ and distinct $1 \leq i, j \leq \grdeg$ with $\undx + \unde_i, \undx + \unde_j \in \hcube{\undr}{}$, we have
			\begin{align} \label{eq:flag_matroids}
				\rk(\undx + \unde_i) + \rk(\undx + \unde_j) \geq \rk(\undx + \unde_i + \unde_j) + \rk(\undx).
			\end{align}
			Removing an element from each independent set of $\mat^i_{x_i + 1}$ results in an independent set of $\mat^i_{x_i}$. This implies that $\rk(\undx + \unde_i) \leq \rk(\undx) + 1$. We thus have $\rk(\undx) \leq \rk(\undx + \unde_i + \unde_j) \leq \rk(\undx) + 2$.
			
			Let $\unda = \undx + \unde_i$, $\undb = \undx + \unde_j$, and $\undc = \undx + \unde_i + \unde_j$. Three cases can occur, depending on whether $\rk(\undc) = \rk(\undx)$, $\rk(\undx) + 1$, or $\rk(\undx) + 2$.
			
			\noindent $\bullet$ In the first case, $\rk(\undc) = \rk(\undx)$, inequality~\eqref{eq:flag_matroids} holds trivially.
			
			\noindent $\bullet$ Consider the third case $\rk(\undc) = \rk(\undx) + 2$. In this case, using the inequality $\rk\mleft(\undy + \unde_k\mright) \leq \rk\mleft(\undy\mright) + 1$ for all $\undy, \undy + \unde_k \in \hcube{\undr}{}$, we infer that $\rk(\unda) = \rk(\undb) =\rk(\undx) + 1$, and inequality~\eqref{eq:flag_matroids} holds again trivially.
			
			\noindent $\bullet$ It remains to treat the case $\rk(\undc) = \rk(\undx) + 1$. Let $I = \ssI_1 \cup \dots \cup \ssI_\grdeg$ be a basis of $\mat_{\undx}$ with $\ssI_k$ an independent of $\mat_{x_k}^k$ for $k = 1, \dots, \grdeg$. There exists a basis $J$ of $\mat_{\undc}$ which contains $I$ and an extra element $\ssvare$ of $E$. Write $J = \ssJ_1 \cup \dots \cup \ssJ_\grdeg$ with $\ssJ_k$ an independent of $\mat^k_{c_k}$, for $k = 1, \dots, \grdeg$. Since $J$ is not an independent of $\mat_{\undx}$, $\ssvare$ appears in either $\ssJ_i$ or $\ssJ_j$. Removing it if necessary from one of the two, we can suppose that $\ssvare$ appears in exactly one of the two sets $\ssJ_i$ or $\ssJ_j$, say, without loss of generality, in $\ssJ_i$. Then, $J$ will be an independent set of $\mat_{\unda}$, and so $\rk(\unda) = \rk(\undx) + 1$. This shows that inequality~\eqref{eq:flag_matroids} holds. The theorem follows.
		\end{proof}
		
		Note that if, in the definition of an initial flag matroid, we relax strict inclusions, the same construction as above gives rise as well to matricubes which are not necessarily simple.
	
	\subsection{Operations on matricubes} \label{subsec:operations_on_matricubes}
		
		Let $\matricube$ be a matricube with underlying ground set $\hcube{\undr}{}$, $\undr = (\ssr_1, \dots, \ssr_{\grdeg})$.
		
		\subsubsection{Deletion} \label{subsubsec:deletion}
			
			Let $i \in \{1, \dots, \grdeg\}$. We define the \emph{deletion} of $i$ in $\matricube$, denoted by $\matricube \mysetminus i$, as the matricube with ground set $\hcube{\undr'}{}$, $\undr' = (\ssr_1, \dots, \ssr_i - 1, \dots, \ssr_{\grdeg})$, defined as follows. We view $\hcube{\undr'}{}$ as the subset of $\hcube{\undr}{}$ consisting of all the points $\undx$ with $\ssx_i < \ssr_i$ and define the rank function $\rk'$ of $\matricube \mysetminus i$ to be the restriction of $\rk$ to $\hcube{\undr'}{}$. Obviously, $\rk'$ verifies the axiomatic system \ref{def:axiom_rank_function_1}-\ref{def:axiom_rank_function_2}-\ref{def:axiom_rank_function_3} of matricube rank functions. Furthermore, note that if $\matricube$ is simple, then so is $\matricube \mysetminus i$.
			
			As an example, here is a (simple) matricube $\matricube$ with $\undr = (4, 3)$ (left) and its deletion $\matricube \mysetminus 2$ in the vertical direction (right).
			\[
				\begin{pmatrix}
					3 & 3 & 3 & 4 & 5 \\
					2 & 2 & 2 & 3 & 4 \\
					1 & 2 & 2 & 3 & 4 \\
					0 & 1 & 2 & 3 & 4
				\end{pmatrix}
				\qquad \qquad
				\begin{pmatrix}
					2 & 2 & 2 & 3 & 4 \\
					1 & 2 & 2 & 3 & 4 \\
					0 & 1 & 2 & 3 & 4
				\end{pmatrix}
			\]
		
		\subsubsection{Contraction} \label{subsubsec:contraction}
			
			Let $i \in \{1, \dots, \grdeg\}$. We define the \emph{contraction} of $i$ in $\matricube$, denoted by $\matricube \mycontr i$, as the matricube with ground set $\hcube{\undr'}{}$, $\undr' = (\ssr_1, \dots, \ssr_i - 1, \dots, \ssr_{\grdeg})$, defined as follows. We define an embedding of $\hcube{\undr'}{}$ in $\hcube{\undr}{}$ by sending each point $\undx$ to $\undx + \unde_i$. We then define the rank function $\rk'$ of ${\matricube \mycontr i}$ by setting $\rk'(\undx) \coloneqq \rk(\undx + \unde_i) - \rk(\unde_i)$. The embedding of $\hcube{\undr'}{}$ in $\hcube{\undr}{}$ respects the two operations of $\wedge$ and $\vee$. It is easy to see that $\rk'$ verifies the axiomatic system \ref{def:axiom_rank_function_1}-\ref{def:axiom_rank_function_2}-\ref{def:axiom_rank_function_3} of matricube rank functions. Note that $\matricube \mycontr i$ is not necessarily simple, even if $\matricube$ is so.
			
			As an example, here is a (simple) matricube $\matricube$ with $\undr = (4, 3)$ (left) and its (non-simple) contraction $\matricube \mycontr 2$ in the vertical direction (right).
			\[
				\begin{pmatrix}
					3 & 3 & 3 & 4 & 5 \\
					2 & 2 & 2 & 3 & 4 \\
					1 & 2 & 2 & 3 & 4 \\
					0 & 1 & 2 & 3 & 4
				\end{pmatrix}
				\qquad \qquad
				\begin{pmatrix}
					2 & 2 & 2 & 3 & 4 \\
					1 & 1 & 1 & 2 & 3 \\
					0 & 1 & 1 & 2 & 3
				\end{pmatrix}
			\]
		
		\subsubsection{Minors} \label{subsubsec:minors}
			
			A matricube $\matricube'$ is a \emph{minor} of another matricube $\matricube$ if it can be obtained by a sequence of deletions and contractions from $\matricube$.
			
			Both operations of contraction and deletion respect the representability over a given field $\k$. It follows that if $\matricube$ is representable over $\k$, then all of its minors $\matricube'$ are also representable over $\k$. We will discuss the connection between representability and minors in Section~\ref{subsec:disc_minor}.
	
	\subsection{Duality} \label{subsec:duality}
		
		Let $\matricube$ be a matricube on the ground set $\hcube{\undr}{}$ with rank function $\rk$. The \emph{dual matricube $\matricube^*$} is the matricube on $\hcube{\undr}{}$ with rank function $\rk^*$ defined by
		\[ \rk^*(\undx) \coloneqq \ellone{\undx} + \rk(\undx^c) - \rk(\matricube) \qquad \forall \, \undx \in \hcube{\undr}{}. \]
		Here,
		\[
			\undx^c \coloneqq \undr - \undx = (\ssr_1 - \ssx_1, \dots, \ssr_\grdeg - \ssx_\grdeg)
		\] 
		is the \emph{complement} of $\undx$ and $\ellone{\undx} \coloneqq \ssx_1 + \dots + \ssx_\grdeg$ is the $\ell_1$-norm of $\undx = (\ssx_1, \dots, \ssx_\grdeg)$. A direct verification shows that $\rk^*$ verifies the axiomatic system \ref{def:axiom_rank_function_1}-\ref{def:axiom_rank_function_2}-\ref{def:axiom_rank_function_3} of matricube rank functions. Moreover, $\matricube^*$ has rank $\rk^*(\matricube^*) = \ellone{\undr} - \rk(\matricube)$, and we have $\ssub{(\matricube^*)}!^* = \matricube$. Note however that $\matricube$ can be simple without $\matricube^*$ being so, and vice-versa.
		
		Here is a (simple) matricube $\matricube$ with $\undr = (4, 3)$ (left) and its (non-simple) dual $\matricube^*$ (right).
		\[
			\begin{pmatrix}
				3 & 3 & 3 & 4 & 5 \\
				2 & 2 & 2 & 3 & 4 \\
				1 & 2 & 2 & 3 & 4 \\
				0 & 1 & 2 & 3 & 4
			\end{pmatrix}
			\qquad \qquad
			\begin{pmatrix}
				2 & 2 & 2 & 2 & 2 \\
				1 & 1 & 1 & 2 & 2 \\
				0 & 0 & 0 & 1 & 2 \\
				0 & 0 & 0 & 1 & 2
			\end{pmatrix}
		\]

%%%%%%%%%% Flats

\section{Flats} \label{sec:flats}
	
	In this section, we define flats of matricubes and provide an axiomatic system for them. This extends the axioms of flats in matroid theory.
	
	\subsection{Definition and basic properties}
		
		Let $\matricube = (\hcube{\undr}{}, \rk)$ be a matricube of rank $r$.
		
		\begin{defi}[Flats of a matricube] \label{def:flats}
			 A point $\unda \in \hcube{\undr}{}$ is called a \emph{flat} for $\rk$ if for every $1 \leq i \leq \grdeg$ such that $\unda + \unde_i$ belongs to $\hcube{\undr}{}$, we have $\rk(\unda + \unde_i) = \rk(\unda) + 1$. We denote by $\Fl = \Fl(\matricube) \subseteq \hcube{\undr}{}$ the set of flats of $\matricube$.
		\end{defi}
		
		Here are two (simple) matricubes with $\undr = (4, 3)$ for the first, and $\undr = (5, 4)$ for the second. The flats in each case are depicted in blue.
		\[
			\begin{pmatrix}
				3 & 3 & \color{blue} 3 & \color{blue} 4 & \color{blue} 5 \\
				2 & 2 & \color{blue} 2 & \color{blue} 3 & \color{blue} 4 \\
				\color{blue} 1 & 2 & 2 & 3 & 4 \\
				\color{blue} 0 & \color{blue} 1 & 2 & 3 & 4
			\end{pmatrix}
			\qquad \qquad
			\begin{pmatrix}
				4 & \color{blue} 4 & 5 & 5 & \color{blue} 5 & \color{blue} 6 \\
				\color{blue} 3 & 4 & 5 & 5 & 5 & 6 \\
				\color{blue} 2 & \color{blue} 3 & 4 & 4 & \color{blue} 4 & \color{blue} 5 \\
				\color{blue} 1 & \color{blue} 2 & 3 & \color{blue} 3 & 4 & 5 \\
				\color{blue} 0 & \color{blue} 1 & \color{blue} 2 & 3 & 4 & 5
			\end{pmatrix}
		\]
		
		\begin{prop}[Stability of flats under meet] \label{prop:flat_stability}
			The set $\Fl(\matricube)$ of flats of a matricube $\matricube$ is stable under $\wedge$.
		\end{prop}
		
		We need the following lemma.
		
		\begin{lemma} \label{lem:diamond}
			Let $i \in \{1, \dots, \grdeg\}$ and let $\undx, \undy, \undx + \unde_i, \undy + \unde_i$ be elements of $\hcube{\undr}{}$ with $\undx \subfaceeq \undy$ and $\ssx_i = \ssy_i$. If $\rk\mleft(\undy + \unde_i\mright) = \rk\mleft(\undy\mright) + 1$, then we have $\rk(\undx + \unde_i) = \rk(\undx) + 1$.
		\end{lemma}
		
		\begin{proof}
			This follows from the submodularity of $\rk$ applied to $\undx + \unde_i$ and $\undy$, and Proposition~\ref{prop:inequality_one}.
		\end{proof}
		
		\begin{proof}[Proof of Proposition~\ref{prop:flat_stability}]
			Let $\unda$ and $\undb$ be two flats and let $\undc = \unda \wedge \undb$. Let $i \in \{1, \dots, \grdeg\}$ be such that $\undc + \unde_i$ belongs to $\hcube{\undr}{}$. We have to show that $\rk(\undc + \unde_i) = \rk(\undc) + 1$. By symmetry, we can suppose that $\ssa_i \leq \ssb_i$, that is, $\ssc_i = \ssa_i$. Since $\unda$ is a flat, we have $\rk(\unda + \unde_i) = \rk(\unda) + 1$. Applying Lemma~\ref{lem:diamond} to $\undx = \undc$ and $\undy = \unda$, we conclude.
	 	\end{proof}
	 	
		The above result implies the following.
		
		\begin{thm}
			The set $\Fl(\matricube)$ of flats of a matricube endowed with the partial order $\subfaceeq$ is a graded lattice. The grading is induced by the rank function.
		\end{thm}
		
		\begin{proof}
			$\Fl(\matricube)$ has a minimum and a maximum element, and is stable under meet. It follows that it is a lattice, with the operation $\vee$ between two elements $\unda$ and $\undb$ in $\Fl(\matricube)$ defined as the meet of all the upper bounds $\undc$ for $\unda$ and $\undb$.
			
			Note that for $\unda \subface \undb$ two distinct and comparable flats of $\matricube$, we have $\rk(\unda) < \rk(\undb)$. The statement that $\Fl(\matricube)$ is graded is a consequence of Proposition~\ref{prop:grading} below.
		\end{proof}

		\begin{lemma} \label{lem:layer_withouh_dots}
			Let $\matricube$ be a matricube on the ground set $\hcube{\undr}{}$. Let $\undx$ be an element of $\hcube{\undr}{}$, and let $\undc$ be the minimum flat with $\undc \supfaceeq \undx$. Then, we have 
			\begin{enumerate}
					\item \label{lem:layer_without_dots_1} $\rk(\undc) = \rk(\undx)$. 
				
					\item \label{lem:layer_without_dots_2} Let $i \in \{1, \dots, \grdeg\}$ be such that $\undx + \unde_i \in \hcube{\undr}{}$. Then, 
					\begin{itemize}
						\item if $\ssc_i > \ssx_i$, then we have $\rk(\undx + \unde_i) = \rk(\undx)$.
						
						\item if $\ssc_i = \ssx_i$, then we have $\rk(\undx + \unde_i) = \rk(\undx) + 1$.
					\end{itemize}
			\end{enumerate}
		\end{lemma}
		
		\begin{proof}
			To prove~\ref{lem:layer_without_dots_1}, it will be enough to show there exists a flat $\undb \supfaceeq \undx$ with $\rk(\undb) = \rk(\undx)$. Then, since flats are closed under meet, $\undb$ will coincide with $\undc$ and (1) follows. We proceed by a reverse induction on the $\ell_1$-norm of $\undx$. If $\undx$ is a flat, in particular, if $\undx = \undr$, there is nothing to prove. Otherwise, there exists $i \in \{1, \dots, \grdeg\}$ with $\undy \coloneqq \undx + \unde_i \in \hcube{\undr}{}$ and $\rk(\undx + \unde_i) = \rk(\undx)$. By the induction hypothesis, there is a flat $\undb \supfaceeq \undy$ with $\rk(\undb) = \rk\mleft(\undy\mright) = \rk(\undx)$, and we conclude.
			
			We prove~\ref{lem:layer_without_dots_2}. If $\ssc_i \geq \ssx_i + 1$, then $\undx \subfaceeq \undx + \unde_i \subfaceeq \undc$, and thus $\rk(\undx) \leq \rk(\undx + \unde_i) \leq \rk(\undc) = \rk(\undx)$. We infer that $\rk(\undx + \unde_i) = \rk(\undx)$.
			
			If $\ssc_i = \ssx_i$, then using $\rk(\undc + \unde_i) = \rk(\undc) + 1$, we apply Lemma~\ref{lem:diamond} and deduce that $\rk(\undx + \unde_i) = \rk(\undx) + 1$, as required.
		\end{proof}
		
		We get the following corollary.
		
		\begin{prop} \label{prop:grading}
			Let $\unda \subface \undb$ be two distinct flats of $\matricube$ with $\rk(\undb) \geq \rk(\unda) + 2$. There exists a flat $\unda \subface \undc \subface \undb$ with $\rk(\undc) = \rk(\unda) + 1$.
		\end{prop}
		
		\begin{proof}
			Since $\unda\subface \undb$, there is an index $i$ such that $\ssa_i < \ssb_i$. Let $\undx \coloneqq \unda + \unde_i \subfaceeq \undb$. Since $\unda$ is a flat, $\rk(\undx) = \rk(\unda) + 1$. Let $\undc$ be the minimum flat with $\undc \supfaceeq \undx$. Obviously, $\unda \subface \undc$. Also $\undc \subfaceeq \undb$, as $\undb$ is a flat and $\undc$ is minimum. By Property~\ref{lem:layer_without_dots_1} in the previous lemma, $\rk(\undc) = \rk(\undx) = \rk(\unda) + 1$. We thus have strict inequality $\undc \subface \undb$, and the result follows.
		\end{proof}
	
	\subsection{Axiomatic system of flats} \label{subsec:axiomatic_flats}
		
		For a subset $\Fl$ of $\hcube{\undr}{}$, consider the following properties.
		
		\begin{enumerate}[label=(F\arabic*)]
			\item \label{def:axiom_flats_1} $\undr$ is in $\Fl$.
			
			\item \label{def:axiom_flats_2} $\Fl$ is closed under meet.
			
			\item \label{def:axiom_flats_3} If $\unda$ is an element of $\Fl$ and $i \in \{1, \dots, \grdeg\}$ is such that $\unda + \unde_i \in \hcube{\undr}{}$, then there exists an element $\undb$ in $\Fl$ such that $\undb \supfaceeq \unda + \unde_i$, and $\undb \ssupface \unda$ in $\Fl$.
		\end{enumerate}
		
		We recall that $\undb \ssupface \unda$ means that $\undb$ \emph{covers} $\unda$, i.e., $\undb \supface \unda$ in $\Fl$ and there is no element $\undc \in \Fl$ such that $\undb \supface \undc \supface \unda$.
		
		We also introduce the following non-degeneracy property.

		\begin{enumerate}[label=(F$^{\ast}$)]
			\item \label{def:axiom_flats_star} Each layer $\ssL^i_t$, $i = 1, \dots, \grdeg$, $t \in [\ssr_i]$, contains an element of $\Fl$.
		\end{enumerate}
		
		We prove the following result.
		
		\begin{thm} \label{thm:axiomatic_flats}
			The set of flats $\Fl$ of a matricube $\matricube$ with ground set $\hcube{\undr}{}$ verifies \ref{def:axiom_flats_1}-\ref{def:axiom_flats_2}-\ref{def:axiom_flats_3}.
			Conversely, let $\Fl \subseteq \hcube{\undr}{}$ be a subset verifying \ref{def:axiom_flats_1}-\ref{def:axiom_flats_2}-\ref{def:axiom_flats_3}. Then, $\Fl$ is the set of flats of a matricube $\matricube$ with underling ground set $\hcube{\undr}{}$.
			
			Moreover, the matricube $\matricube$ is simple if, and only if, Property~\ref{def:axiom_flats_star} holds.
		\end{thm}
		
		\begin{remark}
			The first three conditions \ref{def:axiom_flats_1}-\ref{def:axiom_flats_2}-\ref{def:axiom_flats_3} are the matricube analogs of the three axioms that define the set of flats of a matroid. Property \ref{def:axiom_flats_star} requires that the matricube does not contain any \emph{loop}: if a layer does not verify this condition, deleting it, we get a smaller hypercuboid with the same collection $\Fl$ that verifies the same axioms.
		\end{remark}
		
		\begin{remark} \label{rem:zero_automatically_independent}
			The axioms~\ref{def:axiom_flats_2} and~\ref{def:axiom_flats_star} together imply that $\underline 0$ is in $\Fl$. Indeed, \ref{def:axiom_flats_star} yields a flat in $\ssL^i_0$ for each $i = 1, \dots, \grdeg$. The meet of these flats is $\underline 0$, and by~\ref{def:axiom_flats_2} belongs to $\Fl$.
		\end{remark}
		
		\begin{remark} \label{rem:element_covering_automatically_unique}
			It is easy to see that if $\Fl$ verifies~\ref{def:axiom_flats_2}, then the element $\undb$ in~\ref{def:axiom_flats_3} is unique.
		\end{remark}
	
	\subsection{Flats of a matricube verify the axioms}
		
		Let $\matricube$ be a matricube $\hcube{\undr}{}$ and let $\Fl = \Fl(\matricube)$ be the set of flats of $\matricube$. We prove that $\Fl$ verifies properties \ref{def:axiom_flats_1}-\ref{def:axiom_flats_2}-\ref{def:axiom_flats_3}. Moreover, if $\matricube$ is simple, then we show that~\ref{def:axiom_flats_star} holds.
		
		\begin{proof}[Proof of the first part of Theorem~\ref{thm:axiomatic_flats}]
			 Property \ref{def:axiom_flats_1} follows from the definition of flats of a matricube. We already proved property \ref{def:axiom_flats_2} in Proposition~\ref{prop:flat_stability}. It remains to show \ref{def:axiom_flats_3}. If $\rk\mleft(\undr\mright) = \rk(\unda) + 1$, then $\undb = \undr$ satisfies \ref{def:axiom_flats_3}. Otherwise, we have $\rk\mleft(\undr\mright) \geq \rk(\unda) + 2$, and by Proposition~\ref{prop:grading}, we have an element $\undb$ in $\Fl$ of rank $\rk(\unda) + 1$ with $\undb \supfaceeq \unda + \unde_i$. Again, $\undb \ssupface \unda$, as required.
			 
			 Now suppose that $\matricube$ is simple. Property~\ref{def:axiom_flats_star} is a consequence of Lemma~\ref{lem:layer_withouh_dots} above applied to $t \, \unde_i$. Let $\undc$ be the minimum flat with $\undc \supfaceeq t \, \unde_i$. If Property~\ref{def:axiom_flats_star} does not hold for the layer $\ssL^i_t$, then necessary $\ssc_i > t$ and thus $\rk((t + 1) \, \unde_i) = \rk(t \, \unde_i)$, contradicting the simpleness of $\matricube$.
		\end{proof}
		
		In the rest of the section, we prove the second part of the theorem.
	
	\subsection{Diamond property}
		
		We recall the following definition.
		
		\begin{defi}[Diamond property for lattices] \label{def:diamond_property_lattices}
			Let $(L, \subfaceeq)$ be a lattice with the meet and join operations $\wedge$ and $\vee$, respectively. We say that $L$ has the \emph{diamond property} if for every triple of elements $a, b, c \in L$ such that $b \neq c$ and $b$ and $c$ both cover $a$, the join $b \vee c$ covers both $b$ and $c$.
		\end{defi}
		
		\begin{lemma} \label{lem:diamond_implies_gradedness}
			Let $L$ be a lattice that satisfies the diamond property. Then, it admits a grading, i.e., all its maximal chains have the same length.
		\end{lemma}
		
		\begin{proof}
			This is well-known. We give a rather informal proof. We apply the diamond property multiple times to show that two maximal chains $C$ and $C'$ in $L$ have the same length. This is done by induction on the elements of $C$ and $C'$; the ``diamonds'' drawn in the Hasse diagram by repeated application of the diamond property provide a finite sequence of chains of constant length between $C$ and $C'$, ultimately proving that $C$ and $C'$ have the same length.
		\end{proof}

	\subsection{The axiomatic system of flats implies the diamond property}
			
		\begin{lemma} \label{lem:axioms_imply_diamond}
			Let $\Fl \subseteq \hcube{\undr}{}$ be a subset verifying axioms \ref{def:axiom_flats_1} and \ref{def:axiom_flats_2}. Then, $\Fl$ is a lattice. If additionally $\Fl$ verifies axioms \ref{def:axiom_flats_3}, then it verifies the diamond property.
		\end{lemma}
		
		\begin{proof}
			Since $\Fl$ is closed under meet and it has a maximum element, it is a lattice. The join of two elements $\unda$ and $\undb$ in $\Fl$ is the meet of all the elements $\undc \in \Fl$ that verify $\undc \supfaceeq \unda$ and $\undc \supfaceeq \undb$.
			
			Now let $\unda, \undb, \undc \in \Fl$ be such that $\undb \neq \undc$, and $\undb$ and $\undc$ both cover $\unda$. By assumption, $\undb$ and $\undc$ are not comparable, so there exist $j \neq k \in \{1, \dots, \grdeg\}$ such that $\ssb_j > \ssc_j$ and $\ssc_k > \ssb_k$. In fact, $\ssc_j = \ssa_j$, because otherwise, we would have $\unda \lneq \undb \wedge \undc \lneq \undb$, which is impossible because $\undb \ssupface \unda$.
			
			Now, applying~\ref{def:axiom_flats_3} yields $\undx \in \Fl$ such that $\undx \supfaceeq \undc + \unde_j$ and $\undx \ssupface \undc$. We show that $\undx \gneq \undb$. Indeed, first, if $\undx \not \supfaceeq \undb$, then we would have $\unda \lneq \undx \wedge \undb \lneq \undb$, the strictness of the first inequality coming from the inequalities $\ssx_j \geq \ssa_j + 1$ (because $\undx \supfaceeq \undc + \unde_j$ and $\ssc_j = \ssa_j$) and $\ssb_j > \ssa_j$. This would be in contradiction with $\undb \ssupface \unda$. Second, $\undx = \undb$ is not possible because $\undx \supfaceeq \undc + \unde_j \supfaceeq \undc$ and $\undb$ and $\undc$ are not comparable. Therefore, $\undx \gneq \undb$.
			
			Symmetrically, \ref{def:axiom_flats_3} provides an element $\undy \in \Fl$ such that $\undy \supfaceeq \undb + \unde_k$, $\undy \ssupface \undb$ and $\undy \gneq \undc$.
			
			Let $\undu \coloneqq \undx \wedge \undy$ and notice that $\undu$ verifies the chains of inequalities $\undb \subfaceeq \undu \subfaceeq \undy$ and $\undc \subfaceeq \undu \subfaceeq \undx$. In other words, $\undu$ belongs to the interval $\mleft[\undb, \undy\mright]$, defined as the set of all elements $\undz \in \Fl$ such that $\undb \subfaceeq \undz \subfaceeq \undy$. This interval is equal to $\mleft\{\undb, \undy\mright\}$ since $\undy \ssupface \undb$. Likewise, $\undu \in [\undc, \undx] = \{\undc, \undx\}$. Since $\undb$ and $\undc$ are not comparable, $\undy \ssupface \undb$ and $\undx \ssupface \undc$, the only possibility for the sets $\mleft\{\undb, \undy\mright\}$ and $\{\undc, \undx\}$ to have the element $\undu$ in common is that $\undu = \undx = \undy$. This shows that $\undu$ covers $\undb$ and $\undc$. Then, $\undu$ is necessarily equal to $\undb \vee \undc$. We conclude that $\Fl$ verifies the diamond property.
		\end{proof}
		
		Applying Lemma~\ref{lem:diamond_implies_gradedness}, we infer the following.
		
		\begin{prop} \label{prop:flats_axioms_graded}
			Let $\Fl \subseteq \hcube{\undr}{}$ be a subset verifying axioms \ref{def:axiom_flats_1}-\ref{def:axiom_flats_2}-\ref{def:axiom_flats_3}. Then, $\Fl$ is a graded lattice.
		\end{prop}
		
		We denote by $\rk \colon \Fl \to \N \cup \{0\}$ the corresponding grading. The following properties hold.
		
		\begin{enumerate}[label=(\alph*)]				
			\item \label{fact:rank_function_increasing} The function $\rk$ is increasing on $\Fl$, in the following sense: if $\unda \lneq \undb \in \Fl$, then $\rk(\unda) < \rk(\undb)$.
			
			\item \label{fact:rank_function_covering} If $\unda, \undb \in \Fl$ are such that $\undb \ssupface \unda$, then $\rk(\undb) = \rk(\unda) + 1$.
		\end{enumerate} 
	
	\subsection{Proof of the second part of Theorem~\ref{thm:axiomatic_flats}} 
		
		Let $\Fl \subseteq \hcube{\undr}{}$ be a subset verifying the axioms~\ref{def:axiom_flats_1}-\ref{def:axiom_flats_2}-\ref{def:axiom_flats_3}. We define a map 
		\[ \varphi \colon \hcube{\undr}{} \to \Fl \subseteq \hcube{\undr}{} \]
		as follows. For each $\undx \in \hcube{\undr}{}$, we define $\varphi(\undx)$ to be the minimum flat $\undb \in \Fl$ such that $\undb \supfaceeq \undx$.
		
		\begin{lemma} \label{lem:gradedness_implies_facts}
			Notation as above, the map $\varphi$ is well-defined, and $\varphi$ and $\Fl$ have the following properties.
			
			\begin{enumerate}[label=(\roman*)]
				\item \label{property:covering_flat_non_decreasing} The map $\varphi \colon \hcube{\undr}{} \to \Fl \subseteq \hcube{\undr}{}$ is non-decreasing.
				
				\item \label{property:possible_variation_covering_flat} Let $\undx \in \hcube{\undr}{}$ and $i \in \{1, \dots, \grdeg\}$ such that $\undx + \unde_i \in \hcube{\undr}{}$. Then, either $\varphi(\undx + \unde_i) = \varphi(\undx)$, or $\varphi(\undx + \unde_i) \ssupface \varphi(\undx)$.
			\end{enumerate}
		\end{lemma}
		
		\begin{proof}
			The first part is immediate by definition. We prove the second one.
			If we assume that $\varphi(\undx + \unde_i) \neq \varphi(\undx)$, then we must have $\varphi(\undx + \unde_i) \supfaceeq \varphi(\undx) + \unde_i$. Indeed, otherwise, $\varphi(\undx + \unde_i) \supface \varphi(\undx) \supfaceeq \undx + \unde_i$, contradicting the minimality of $\varphi(\undx + \unde_i)$.
			
			Now, \ref{def:axiom_flats_3} yields the existence of $\undb \supfaceeq \varphi(\undx) + \unde_i$ in $\Fl$ such that $\undb \ssupface \varphi(\undx)$. Then, using Property~\ref{property:covering_flat_non_decreasing} above and the definition of $\varphi$, we get $\undb \supfaceeq \varphi(\varphi(\undx) + \unde_i) \supfaceeq \varphi(\undx + \unde_i)$. The chain of inequalities $\undb \supfaceeq \varphi(\undx + \unde_i) \supfaceeq \varphi(\undx)$ and the fact that $\undb \ssupface \undx$ imply $\undb = \varphi(\undx + \unde_i)$. Therefore, we have $\varphi(\undx + \unde_i) \ssupface \varphi(\undx)$, as required.
		\end{proof}
		
		We can now complete the proof of Theorem~\ref{thm:axiomatic_flats}.
		
		\begin{proof}[Proof of the second part of Theorem~\ref{thm:axiomatic_flats}]
			Let $\Fl$ be a subset of $\hcube{\undr}{}$ verifying axioms \ref{def:axiom_flats_1}-\ref{def:axiom_flats_2}-\ref{def:axiom_flats_3}. As we have shown already in Proposition~\ref{prop:flats_axioms_graded}, $\Fl$ is a graded lattice.
			
			We first extend the function $\rk \colon \Fl \to \N\cup\{0\}$ to a function $\rk \colon \hcube{\undr}{} \to \N\cup\{0\}$ by setting, for each $\undx \in \hcube{\undr}{}$, $\rk(\undx) \coloneqq \rk(\varphi(\undx))$. We claim that $\rk$ is the rank function of a matricube. We will use the properties proven in Lemma~\ref{lem:gradedness_implies_facts}.
			
			Part~\ref{property:possible_variation_covering_flat} of the lemma implies directly that $\rk(\undx + \unde_i) \leq \rk(\undx) + 1$ for each $\undx \in \hcube{\undr}{}$, proving~\ref{def:axiom_rank_function_1}.

			The fact that $\rk$ is non-decreasing on $\hcube{\undr}{}$ is a consequence of Property~\ref{property:covering_flat_non_decreasing} in the lemma and Fact~\ref{fact:rank_function_increasing} stipulating that $\rk$ is increasing on elements of $\Fl$.
			
			We show that $\rk$ is submodular on $\hcube{\undr}{}$. By Theorem~\ref{thm:eq_submodularity_diamond_property} that we will prove in Section~\ref{sec:diamond_property}, it is sufficient to prove the diamond property for functions on hypercuboids (see Section~\ref{sec:diamond_property}). Let $\undx \in \hcube{\undr}{}$ and $i \neq j \in \{1, \dots, \grdeg\}$ be such that $\undx + \unde_i + \unde_j \in \hcube{\undr}{}$. Using Property~\ref{property:possible_variation_covering_flat}, we may assume that $\rk(\undx + \unde_i + \unde_j) = \rk(\undx + \unde_j) + 1$, and then need to prove that $\rk(\undx + \unde_i) = \rk(\undx) + 1$. The equality $\rk(\undx + \unde_i + \unde_j) = \rk(\undx + \unde_j) + 1$ means that $\varphi(\undx + \unde_i + \unde_j) \supface \varphi(\undx + \unde_j)$. This implies that $\varphi(\undx + \unde_j) \in \ssL^i_{\ssx_i}$ (as otherwise, we would get $\varphi(\undx + \unde_j)_i \geq \ssx_i + 1$, that is, $\varphi(\undx + \unde_j) \geq \undx + \unde_i + \unde_j$, and so we would have $\varphi(\undx + \unde_i + \unde_j) = \varphi(\undx + \unde_j)$). Since $\varphi(\undx + \unde_j) \supfaceeq \varphi(\undx)$, this in turn implies that $\varphi(\undx) \in \ssL^i_{x_i}$. However, $\varphi(\undx + \unde_i)_i \geq \ssx_i + 1$. We infer that $\varphi(\undx + \unde_i) \supface \varphi(\undx)$, and thus $\rk(\undx + \unde_i) = \rk(\undx) + 1$, as desired.
			
			We have shown that $\rk$ is the rank function of a matricube with ground set $\hcube{\undr}{}$. The fact that the set of flats of $\rk$ is exactly $\Fl$ is immediate by the definition of $\rk$ and the map $\varphi$, and the fact that $\rk$ is increasing on $\Fl$, see~\ref{fact:rank_function_increasing}.
			
			It remains to show that if \ref{def:axiom_flats_star} holds, then $\matricube$ is simple. Let $i \in \{1, \dots, \grdeg\}$. We show by induction that for every $0 \leq t \leq \ssr_i$, we have $\rk(t \, \unde_i) = t$. The base case $t = 0$ holds by definition. We now suppose that $\rk(t \, \unde_i) = t$ with $0 \leq t < \ssr_i$ and show that $\rk((t + 1) \, \unde_i) = t + 1$. \ref{def:axiom_flats_star} implies that $\varphi(t \, \unde_i) \in \ssL^i_t$ and $\varphi((t + 1) \, \unde_i) \in \ssL^i_{t + 1}$. In particular, $\varphi((t + 1) \, \unde_i) \neq \varphi(t \, \unde_i)$. Again, Lemma~\ref{lem:gradedness_implies_facts} implies that $\varphi((t + 1) \, \unde_i) \ssupface \varphi(t \, \unde_i)$ which, using that $\rk$ is increasing on $\Fl$ implies $\rk((t + 1) \, \unde_i) = \rk(\varphi((t + 1) \, \unde_i)) = \rk(\varphi(t \, \unde_i)) + 1 = \rk(t \, \unde_i) + 1 = t + 1$, as desired.
		\end{proof}
		
%%%% 

\section{Circuits} \label{sec:circuits}
	
	We define circuits in matricubes and provide an intrinsic axiomatic system for them.
	
	\subsection{Duality and circuits} \label{subsec:circuits_and_duality}
		
		Let $\matricube$ be a matricube on the ground set $\hcube{\undr}{}$, and denote by $\matricube^*$ its dual. Denote by $\Fl(\matricube^*)$ the set of flats of the dual matricube, and consider
		\[
			\CCir \coloneqq \mleft\{\unda^c \, \st \, \unda \in \Fl(\matricube^*) \mright\} \subseteq \hcube{\undr}{},
		\]
		where, we recall, $\unda^c = \undr - \unda$. Since $\Fl(\matricube)$ is closed under meet, $\CCir$ will be closed under the join operation. 
		
		Given a subset $\mathcal A \subset \hcube{\undr}{}$, we say that an element $\unda$ of $\mathcal A$ is \emph{join-irreducible in $\mathcal A$} if it is not the join of any set of elements of $\mathcal A \setminus \{\unda\}$.
		
		\begin{defi}[Circuits] \label{def:circuits}
			The collection of circuits of $\matricube$, denoted by $\Cir$, is defined as the set of nonzero join-irreducible elements of $\CCir$.
		\end{defi}
		
		Here is a (simple) matricube $\matricube$ with $\undr = (5, 4)$. On the left, $\matricube$ is represented by its rank function $\rk$, with its circuits in red and the join-reducible elements of $\CCir$ in blue. On the right, the dual $\matricube^*$ of $\matricube$, which is not simple, is represented by its rank function $\rk^*$, with its flats in teal.
		\begin{equation} \label{ex:circuits_matricube}
			\begin{pmatrix}
				4 & 4 & \color{red} 4 & \color{blue} 4 & 5 & 6 \\
				3 & 3 & 4 & 4 & 5 & 6 \\
				2 & 2 & 3 & \color{red} 3 & 4 & 5 \\
				1 & \color{red} 1 & 2 & 3 & 4 & 5 \\
				0 & 1 & 2 & 3 & 4 & 5
			\end{pmatrix}
			\qquad \qquad 
			\begin{pmatrix}
				3 & 3 & 3 & 3 & 3 & \color{teal} 3 \\
				2 & 2 & 2 & 2 & \color{teal} 2 & 3 \\
				1 & 1 & \color{teal} 1 & 2 & 2 & 3 \\
				1 & 1 & 1 & 2 & 2 & 3 \\
				0 & 0 & \color{teal} 0 & \color{teal} 1 & 2 & 3
			\end{pmatrix}
		\end{equation}
		
		Obviously, by definition, $\Cir$ determines $\CCir$, and therefore, gives the set of flats of the dual matroid $\matricube^*$. By Theorem~\ref{thm:axiomatic_flats}, this implies that $\Cir$ determines $\matricube$.
	
	\subsection{Axiomatic system of circuits}
		
		For a subset $\Cir$ of $\hcube{\undr}{}$, denote by $\CCir$ the join-closure of $\Cir$, obtained by taking the join of any set of elements of $\Cir$. Consider the following set of properties:
		
		\begin{enumerate}[label=(C\arabic*)]
			\item \label{def:axiom_circuits_1} $\underline 0$ is \emph{not} in $\Fl$.
			
			\item \label{def:axiom_circuits_2} All elements of $\Cir$ are join-irreducible in $\Cir$.
			
			\item \label{def:axiom_circuits_3} If $\unda \in \CCir$ and $i \in \{1, \dots, \grdeg\}$ is such that $\unda - \unde_i \in \hcube{\undr}{}$, then there exists an element $\undb \subfaceeq \unda - \unde_i$ in $\CCir \cup \{\underline 0\}$ such that $\undb \ssubface \unda$ in $\CCir \cup \{\underline 0\}$.
		\end{enumerate}
		
		We also introduce the following simpleness property.
		
		\begin{enumerate}[label=(C$^{\ast}$)]
			\item \label{def:axiom_circuits_star} For every $i \in \{1, \dots, \grdeg\}$ and $t \in [\ssr_i]$, $t \, \unde_i$ is \emph{not} in $\Cir$.
		\end{enumerate}
		
		We prove the following result.
		
		\begin{thm} \label{thm:axiomatic_circuits}
			The set of circuits $\Cir$ of a matricube $\matricube$ with ground set $\hcube{\undr}{}$ verifies \ref{def:axiom_circuits_1}-\ref{def:axiom_circuits_2}-\ref{def:axiom_circuits_3}. Conversely, let $\Cir \subseteq \hcube{\undr}{}$ be a subset verifying \ref{def:axiom_circuits_1}-\ref{def:axiom_circuits_2}-\ref{def:axiom_circuits_3}. Then, $\Cir$ is the set of circuits of a matricube $\matricube$ with underlying ground set $\hcube{\undr}{}$.
			
			Moreover, the matricube $\matricube$ is simple if, and only if, Property~\ref{def:axiom_circuits_star} holds.
		\end{thm}
		
		\begin{proof}
			$(\Longrightarrow)$ Let $\Cir$ be the set of circuits of a matricube $\matricube$ with ground set $\hcube{\undr}{}$. Properties~\ref{def:axiom_circuits_1} and~\ref{def:axiom_circuits_2} hold by definition of $\Cir$ (see Definition~\ref{def:circuits}). As for Property~\ref{def:axiom_circuits_3}, it is a translation through duality of Property~\ref{def:axiom_flats_3} which holds for the set of flats of the dual matricube $\matricube^*$.
			
			Assume moreover that $\matricube$ is simple. Let $i \in \{1, \dots, \grdeg\}$ and $t \in [\ssr_i]$. We have $\rk(t \, \unde_i) = t$ and therefore, denoting by $\rk^*$ the rank function on $\matricube^*$, we have $\rk^*\mleft(\undr - t \, \unde_i\mright) = \rk^*\mleft(\undr\mright)$. This implies that for every $i$ and $t$, $\undr - t \, \unde_i$ is not a flat in $\matricube^*$, which means that $t \, \unde_i$ is not in $\CCir$. As a consequence, $t \, \unde_i \notin \Cir$.
			
			$(\Longleftarrow)$ Let $\Cir \subseteq \hcube{\undr}{}$ be a subset verifying \ref{def:axiom_circuits_1}-\ref{def:axiom_circuits_2}-\ref{def:axiom_circuits_3}, and $\CCir$ the join-closure of $\Cir$. Define
			\[ \Fl \coloneqq \mleft\{\unda^c \st \unda \in \CCir\mright\} \cup \mleft\{\undr\mright\}. \]
			We claim that $\Fl$ is the set of flats of a matricube. We need to show that it satisfies \ref{def:axiom_flats_1}-\ref{def:axiom_flats_2}-\ref{def:axiom_flats_3}. By construction, \ref{def:axiom_flats_1} and~\ref{def:axiom_flats_2} hold. Property~\ref{def:axiom_flats_3} is a translation through duality of Property~\ref{def:axiom_circuits_3} which holds for $\Cir$. Therefore, $\Fl$ is the set of flats of a matricube. We denote the dual of this matricube by $\matricube$, so that $\Fl = \Fl(\matricube^*)$. It is immediate by construction that $\Cir$ is its set of circuits of $\matricube$.
	
			Let $\rk$ be the rank function of $\matricube$ and $\rk^*$ that of $\matricube^*$. 
			
			To prove the last assertion, assume that for every $i \in \{1, \dots, \grdeg\}$ and $t \in [\ssr_i]$, $t \, \unde_i$ is not in $\Cir$. Then, for every $i$ and $t$, $\undr - t \, \unde_i$ is not a flat of $\matricube^*$. By a simple induction, this implies that for every $i$ and $t$, $\rk^*\mleft(\undr - t \, \unde_i\mright) = \rk^*\mleft(\undr\mright)$, and consequently, $\rk(t \, \unde_i) = t$. We infer that $\matricube$ is simple.
		\end{proof}

\section{Independents} \label{sec:independents}
	
	In this section, we define the independents of a matricube and study their properties. As in the case of flats and circuits, we give the axiomatic system of independents of a matricube.
	
	\subsection{Definition and basic properties}
		
		Let $\matricube$ be a matricube on the ground set $\hcube{\undr}{}$.
		
		\begin{defi}[Independents of a matricube] \label{def:independents_matricube}
			We say that a point $\undp$ of $\hcube{\undr}{}$ is called an \emph{independent of $\matricube$} if for each $i = 1, \dots, \grdeg$ such that $\undp - \unde_i \in \hcube{\undr}{}$, we have $\rk\mleft(\undp - \unde_i\mright) = \rk\mleft(\undp\mright) - 1$. We denote by $\Ind(\matricube) \subseteq \hcube{\undr}{}$ the set of independents of $\matricube$.
		\end{defi}
		Here are two matricubes with $\undr = (4, 3)$ for the first, and $\undr = (5, 4)$ for the second. The independents in each case are depicted in blue.	
		\[
			\begin{pmatrix}
				\color{blue} 3 & 3 & 3 & \color{blue} 4 & \color{blue} 5 \\
				\color{blue} 2 & 2 & 2 & 3 & 4 \\
				\color{blue} 1 & \color{blue} 2 & 2 & 3 & 4 \\
				\color{blue} 0 & \color{blue} 1 & \color{blue} 2 & \color{blue} 3 & \color{blue} 4
			\end{pmatrix}
			\qquad \qquad
			\begin{pmatrix}
				\color{blue} 4 & 4 & 5 & 5 & 5 & 6 \\
				\color{blue} 3 & \color{blue} 4 & \color{blue} 5 & 5 & 5 & \color{blue} 6 \\
				\color{blue} 2 & \color{blue} 3 & \color{blue} 4 & 4 & 4 & 5 \\
				\color{blue} 1 & \color{blue} 2 & \color{blue} 3 & 3 & 4 & 5 \\
				\color{blue} 0 & \color{blue} 1 & \color{blue} 2 & \color{blue} 3 & \color{blue} 4 & \color{blue} 5
			\end{pmatrix}
		\]
		
		The following proposition provides a list of properties of independent sets in a matricube.
		
		\begin{prop} \label{prop:properties_independents_matricube}
			Let $\Ind(\matricube)$ be the set of independents of a matricube $\matricube$. The following properties hold.
			
			\begin{itemize} 
				\item $\Ind(\matricube)$ is non-empty and closed under meet.
				
				\item For every independent $\undp \in \Ind(\matricube)$ and every distinct elements $\ssi_1, \dots, \ssi_k \in \{1, \dots, \grdeg\}$ with $\ssp_{\ssi_j} \neq 0$, $j = 1, \dots, k$, we have $\rk\mleft(\undp - \unde_{\ssi_1} - \dots - \unde_{\ssi_k}\mright) = \rk(\unda) - k$.
			\end{itemize}
		\end{prop}
		
		\begin{proof}
			Both statements follow from Lemma~\ref{lem:diamond}, as in the proof of Proposition~\ref{prop:flat_stability}.
		\end{proof}
		
		Note that $\Ind(\matricube)$ in general does not have a maximum, and lacks the existence of a join.
		
		In order to study more refined properties of independents, we will associate a notion of \emph{size} to each independent element in $\matricube$ by defining a \emph{removal} operation on elements of $\Ind(\matricube)$.
	
	\subsection{Removal and size} \label{subsec:removal_and_size}
		
		Let $\Jind$ be a subset of $\hcube{\undr}{}$. Let $\unda$ be an element of $\Jind$ and $i \in \{1, \dots, \grdeg\}$ such that $\ssa_i \neq 0$. If there is at least one element $\undb \subface \unda$ in $\Jind$ that differs from $\unda$ only in the $i$-th component, we define $\unda \setminus i$ to be such an element in $\Jind$ with the largest $i$-th coordinate. In this case, we say that $\unda \setminus i$ is the \emph{removal of $i$ in $\unda$ in $\Jind$.}
		
		\begin{defi} \label{def:orderable_subset}
			Let $\Jind$ be a subset of $\hcube{\undr}{}$.
				
			\begin{enumerate}[label=(\alph*)]
				\item We say that \emph{removals exist in $\Jind$} if for every $\unda \in \Jind$ and $i \in \{1, \dots, \grdeg\}$, if $\ssa_i \geq 1$, the removal $\unda \setminus i$ exists in $\Jind$.
			\end{enumerate}
			(If removals exist in $\Jind$, then necessarily, we have $\underline 0 \in \Jind$. Moreover, for every element $\unda \in \Jind$, there exists a sequence of removals in $\Jind$, that reduces $\unda$ to $\underline 0$.)
			
			\begin{enumerate}[label=(\alph*),resume]	
				\item We say that $\Jind$ is \emph{orderable} if removals exist in $\Jind$ and for every $\unda \in \Jind$, all the sequences of removals in $\Jind$ that bring $\unda$ to $\underline 0$ have the same length.
			\end{enumerate}
			If $\Jind$ is orderable, we define the \emph{size} of each element $\unda \in \Jind$ denoted by $|\unda|$ as the number of removals needed to reduce $\unda$ to $\underline 0$.
		\end{defi}
		
		In Lemma~\ref{lem:second_axiom_implies_orderability} below, we formulate a simple orderability criterion.
	
	\subsection{Axiomatic system of independents}
		
		We first make the following definition, which turns out to be useful in the proof of the main theorem of this section.
		
		\begin{defi} \label{def:difference_and_energy}
			Let $\unda, \undb \in \hcube{\undr}{}$ be two elements. We define:
			\[ 
				\Delta(\unda, \undb) \coloneqq \mleft \{k = 1, \dots, \grdeg \st \ssa_k < \ssb_k\mright \} \qquad \textrm{and} \qquad E(\unda, \undb) \coloneqq \sum_{k \in D(\unda, \undb)} (\ssb_k - \ssa_k). \qedhere
			\]
		\end{defi}
		
		For a subset $\Ind$ of $\hcube{\undr}{}$, consider the following property:
		
		\begin{enumerate}[label=(I\arabic*)]
			\item \label{def:axiom_independents_1} Removals exist in $\Ind$ and the following holds. For all $\undp \in \Ind$ and removals $\undp \setminus i$ and $\undp \setminus j$, with $i, j \in \{1, \dots, \grdeg\}$, the meet $\undq \coloneqq \mleft(\undp \setminus i\mright) \wedge \mleft(\undp \setminus j\mright)$ belongs to $\Ind$ and, moreover, the two intervals $\mleft[\undq, \undp \setminus i \mright]$ and $\mleft[\undq, \undp \setminus j\mright]$ in $\Ind$ have the same size.
		\end{enumerate}
		(The interval $[\unda, \undb]$ in $\Ind$ is defined as the set of all $\undc \in \Ind$ such that $\unda \subfaceeq \undc \subfaceeq \undb$.) 
		
		\smallskip
		
		It follows from Lemma~\ref{lem:second_axiom_implies_orderability}, proved in Section~\ref{sec:proof_orerability}, that a subset $\Ind \subseteq \hcube{\undr}{}$ that verifies \ref{def:axiom_independents_1} is orderable. We can thus define the size $|\unda|$ of each element $\unda \in \Ind$. This enables us to formulate the second property of interest:
		
		\begin{enumerate}[label=(I\arabic*)]
			\setcounter{enumi}{1}
			
			\item \label{def:axiom_independents_2} $|\cdot|$ is increasing on independents, i.e., for all $\unda, \undb \in \Ind$ such that $\unda \subface \undb$, we have $|\unda| < |\undb|$. Moreover, let $\unda$ and $\undb$ be two elements of $\Ind$ such that $|\unda| < |\undb|$ and $D(\unda, \undb)$ contains at least two elements. Then, there exists $\undc \in \Ind$ that verifies:
			\begin{itemize}
				\item $\undc \subfaceeq \unda \vee \undb$,
				
				\item $|\undc| > |\unda|$.
				
				\item There exists $i \in D(\unda,\undb)$ such that $\ssc_i < \ssb_i$.
			\end{itemize}
		\end{enumerate}
		
		We also introduce the following notion of simpleness.
		
		\begin{enumerate}[label=(I$^{\ast}$)]
			\item \label{def:axiom_independents_star} For $i = 1, \dots, \grdeg$, the points $t \, \unde_i$ for $t \in [\ssr_i]$ are all in $\Ind$.
		\end{enumerate}
		
		This is the main result of this section.
		
		\begin{thm} \label{thm:axiomatic_independents}
			The set of independents $\Ind(\matricube)$ of a matricube $\matricube$ with ground set $\hcube{\undr}{}$ verifies \ref{def:axiom_independents_1}-\ref{def:axiom_independents_2}. Conversely, let $\Ind \subseteq \hcube{\undr}{}$ be a subset that verifies \ref{def:axiom_independents_1}-\ref{def:axiom_independents_2}. Then, $\Ind$ is the set of independents of a matricube $\matricube$ with underling ground set $\hcube{\undr}{}$.
			
			Moreover, the matricube $\matricube$ is simple if, and only if, $\Ind$ verifies~\ref{def:axiom_independents_star}.
		\end{thm}
		
		\begin{remark}
			Axiom~\ref{def:axiom_independents_1} is an analog of the hereditary property for independents of matroids. It also implies that $\underline 0 \in \Ind$, analog of the first axiom of independents in matroids. Axiom~\ref{def:axiom_independents_2} plays the role of the augmentation property for independents. These axioms take into account the more singular nature of independents in the context of matricubes: for example, in matroids, all maximal (for inclusion) independents have the same cardinality, whereas in matricubes, maximal independents (for the partial order $\subfaceeq$) can have different sizes as one of the two examples above shows (see Section~\ref{sec:special_features_independents} for further discussion). Axiom~\ref{def:axiom_independents_star} requires that there is no ``loop''.
		\end{remark}
	
	\subsection{Independents of a matricube verify the axioms}
		
		Let $\matricube$ be a matricube $\hcube{\undr}{}$. We prove that $\Ind(\matricube)$ verifies properties \ref{def:axiom_independents_1}-\ref{def:axiom_independents_2}. Moreover, if $\matricube$ is simple, then we prove that~\ref{def:axiom_independents_star} holds.
		
		The proof shows that the size function on $\Ind(\matricube)$ coincides with the rank function.
		
		\begin{proof}[Proof of the first part of Theorem~\ref{thm:axiomatic_independents}]
			
			We start by proving~\ref{def:axiom_independents_1}. Let $\undp \in \Ind(\matricube)$ and let $i \in \{1, \dots, \grdeg\}$ be such that $\ssp_i \neq 0$. We claim that $\unda \coloneqq \undp - \ssp_i \unde_i \in \Ind(\matricube)$. By design, $\ssa_i = 0$. Now, let $j \in \{1, \dots, \grdeg\}$ be different from $i$ such that $\ssa_j \neq 0$. We have $\ssp_j = \ssa_j \neq 0$. We have $\rk\mleft(\undp - \unde_j\mright) = \rk\mleft(\undp\mright) - 1$. Applying Lemma~\ref{lem:diamond} with $\undx = \unda - \unde_j$ and $\undy = \undp - \unde_j$, we get $\rk(\unda - \unde_j) = \rk(\unda) - 1$. This shows that $\unda \in \Ind(\matricube)$. Therefore, removals exist in $\Ind(\matricube)$.
			
			We next show that $\rk\mleft(\undp \setminus i\mright) = \rk\mleft(\undp\mright) - 1$. For the sake of a contradiction, suppose this not being the case, that is, $\rk\mleft(\undp \setminus i\mright) \leq \rk\mleft(\undp\mright) - 2$. Then, there would exist $\undp \setminus i\subface \undb \subface \undp$ such that $\rk(\undb) = \rk\mleft(\undp\mright) - 1$ and $\rk(\undb - \unde_i) = \rk(\undb) - 1$. Note that $\ssb_j = \ssp_j$ for all $j\neq i$. Applying again Lemma~\ref{lem:diamond} as above, we infer that $\undb$ belongs to $\Ind(\matricube)$. This would be a contradiction to the definition of the removal.
			
			This implies that a sequence of removals bringing $\undp \in \Ind(\matricube)$ to $\underline 0$ has size precisely $\rk\mleft(\undp\mright)$.
			
			Now, for distinct $i, j \in \{1, \dots, \grdeg\}$, we consider the removals $\undp \setminus i$ and $\undp \setminus j$ in $\Ind(\matricube)$, as well as $\undq \coloneqq \mleft(\undp \setminus i\mright) \wedge \mleft(\undp \setminus j\mright)$. By Proposition~\ref{prop:properties_independents_matricube}, $\undq \in \Ind(\matricube)$. Note that we have $\rk\mleft(\undp \setminus j\mright) = \rk\mleft(\undp\mright) - 1 = \rk\mleft(\undp \setminus i\mright)$.
			
			The element $\undq$ differs from $\undp\setminus i$ only in the $j$-th component, and therefore can be obtained from it by a sequence of removals of $j$. It follows that $\mleft[\undq, \undp \setminus i\mright]$ has cardinality $\rk\mleft(\undp \setminus i\mright) - \rk\mleft(\undq\mright) + 1$. Similarly, $\mleft[\undq, \undp \setminus j\mright]$ has cardinality $\rk\mleft(\undp \setminus j\mright) - \rk(\mleft(\undq\mright)) + 1$.
			
			We conclude that the two intervals $\mleft[\undq, \undp \setminus i\mright]$ and $\mleft[\undq, \undp \setminus j\mright]$ in $\Ind(\matricube)$ have the same cardinality, and~\ref{def:axiom_independents_1} follows. We thus get a well-defined size function $|\cdot|$ on $\Ind(\matricube)$. As the proof shows, we have $|\unda| = \rk(\unda)$ for every $\unda \in \Ind$.
			
			\smallskip
			
			The first half of Property~\ref{def:axiom_independents_2} results from the fact that if $\unda \subface \undb$ are two independents, then $\Delta(\unda, \undb) \neq \varnothing$. Then, taking $k \in \Delta(\unda, \undb)$, we get $\rk(\unda) \leq \rk(\undb - \unde_k) = \rk(\undb) - 1$.
			
			For the second half of Property~\ref{def:axiom_independents_2}, let $\unda, \undb$ be two independents such that $|\unda| < |\undb|$ and $|\Delta(\unda, \undb)| \geq 2$. We consider two cases depending on whether $|\unda| \leq |\undb| - 2$ or $|\unda| = |\undb| - 1$.
			
			First, consider the case $|\unda| \leq |\undb| - 2$. Let $i \in \Delta(\unda, \undb)$. Since $\ssb_i > \ssa_i \geq 0$, we can define $\undc \coloneqq \undb \setminus i$. Note that $|\undc| = |\undb| - 1 > |\unda|$. Furthermore, by construction, $\undc \subfaceeq \unda \vee \undb$ and $\ssc_i < \ssb_i$. This shows that $\undc$ is suitable.
			
			We now consider the case $|\unda| = |\undb| - 1$. Let $\undd \coloneqq \unda \vee \undb$. Let $\undy \subfaceeq \undd$ be an element of $\hcube{\undr}{}$ minimal for $\subfaceeq$ under the constraints that $\rk\mleft(\undy\mright) = \rk(\undd)$ and for all $k \notin \Delta(\unda, \undb)$, $\ssy_k = \ssd_k$ (note that for all those $k$, we have $\ssd_k = \ssa_k$). For all $k\in \Delta(\unda, \undb)$ with $\undy - \unde_k \in \hcube{\undr}{}$, we thus have $\rk\mleft(\undy - \unde_k\mright) = \rk\mleft(\undy\mright) - 1$.
			
			Next, let $\undx \subfaceeq \undy$ be an element of $\hcube{\undr}{}$ minimal for $\subfaceeq$ under the constraint that $\rk(\undx) = \rk\mleft(\undy\mright) = \rk(\undd)$, and for all $k \in \Delta(\unda, \undb)$, $\ssx_k = \ssy_k$. Since for all $k \in \Delta(\unda, \undb)$ with $\ssx_k = \ssy_k > 0$, we have $\rk\mleft(\undy - \unde_k\mright) = \rk\mleft(\undy\mright) - 1$, Lemma~\ref{lem:diamond} implies that $\rk(\undx - \unde_k) = \rk(\undx) - 1$. Moreover, by the choice of $\undx$, we also have, for all $i \notin \Delta(\unda, \undb)$, $\rk(\undx - \unde_i) = \rk(\undx) - 1$ provided that $\undx - \unde_i$ belongs to $\hcube{\undr}{}$. Therefore, combining all this, we conclude that $\undx \in \Ind(\matricube)$.
			
			If $\undy \subface \undd$, let $\undc \coloneqq \undx \in \Ind(\matricube)$. There exists then $k \in \Delta(\unda, \undb)$ such that $\ssy_k < \ssd_k = \ssb_k$. As a consequence, $\ssc_k < \ssb_k$. Moreover, $\rk(\undc) = \rk(\undd) \geq \rk(\undb) = \rk(\unda) + 1$ and therefore $|\undc| > |\unda|$. Since by construction $\undc \subface \undd$, $\undc$ is suitable.
			
			It remains to consider the case $\undy = \undd$. This means that for every $k \in \Delta(\unda, \undb)$, $\rk(\undd - \unde_k) < \rk(\undd)$. We claim that in this case, the strict inequality $\rk(\undd) > \rk(\undb)$ holds. Indeed, for the sake of a contradiction, suppose $\rk(\undd) = \rk(\undb)$. Since $|\Delta(\unda, \undb)| \geq 2$, there are two distinct elements $i, j \in \Delta(\unda, \undb)$, and for these $i, j$, we would have $\rk(\undd - \unde_i) = \rk(\undd - \unde_j) = \rk(\undd) - 1$. By submodularity (see Lemma~\ref{lem:diamond}), we would get $\rk(\undd - \unde_i - \unde_j) = \rk(\undd) - 2$. Since $i, j \in \Delta(\unda, \undb)$, we have $\unda \subfaceeq \undd - \unde_i - \unde_j$, and therefore, we would have $\rk(\unda) \leq \rk(\undd) - 2 = \rk(\undb) - 2$, contradicting the assumption that $|\unda| = |\undb| - 1$. This proves the claim that $\rk(\undd) > \rk(\undb)$.
			
			Let now $\undc \coloneqq \undx \setminus k$ for an element $k\in \Delta(\unda,\undb)$. We have $|\undc| = \rk(\undc) = \rk(\undx) - 1 = \rk(\undd) - 1 \geq \rk(\undb) = |\undb| > |\unda|$, and therefore $|\undc| > |\unda|$. Besides, $\ssc_k < \ssx_k = \ssd_k = \ssb_k$, and obviously $\undc \subfaceeq \undd$. This shows that $\undc$ is suitable in this last case. This ends the proof of~\ref{def:axiom_independents_2}. We have proved that $\Ind(\matricube)$ verifies~\ref{def:axiom_independents_1} and~\ref{def:axiom_independents_2}.

			To finish the proof, note that if $\matricube$ is simple, then~\ref{def:axiom_rank_function_1nd} immediately implies~\ref{def:axiom_independents_star}.
		\end{proof}
	
	\subsection{Orderability lemma}	\label{sec:proof_orerability}
		
		Before going to the proof of the second part of Theorem~\ref{thm:axiomatic_independents}, we show the following criterion for orderability.
		
		\begin{lemma} \label{lem:second_axiom_implies_orderability}
			Let $\Jind$ be a subset of $\hcube{\undr}{}$. The following are equivalent:
			
			\begin{enumerate}[label=(\arabic*)]
				\item \label{ord1} $\Jind$ satisfies \ref{def:axiom_independents_1}.
				
				\item \label{ord2} $\Jind$ is orderable in the sense of Definition~\ref{def:orderable_subset}.
			\end{enumerate}
			
		\end{lemma}
		
		\begin{proof}
			We first prove \ref{ord1} $\Rightarrow$ \ref{ord2}. Assume~\ref{def:axiom_independents_1} holds. Proceeding by induction under the partial order $\subfaceeq$, we show that for every $\unda \in \Jind$, the following property holds:
			\[ P(\unda): \quad \text{All the sequences of removals that bring } \unda \text{ to } \underline 0 \text{ have the same length.} \]
			Obviously, $P(\underline 0)$ holds. Let $\unda \in \Jind$ be an element such that for every $\undb \in \Jind$ with $\undb \subface \unda$, $P(\undb)$ holds. We prove that $P(\unda)$ is true. Let
			\[ \unda = \undb_0 \supface \undb_1 \supface \cdots \supface \undb_k = \underline 0 \quad \text{and} \quad \unda = \undc_0 \supface \undc_1 \supface \cdots \supface \undc_\ell = \underline 0 \]
			be two sequences of removals bringing $\unda$ to $\underline 0$. We need to prove that $k = \ell$. Let $i, j \in \{1, \dots, \grdeg\}$ be such that $\undb_1 = \unda \setminus i$ and $\undc_1 = \unda \setminus j$, and let $\undq \coloneqq \undb_1 \wedge \undc_1 \in \Jind$. Let
			\[ \undq = \undx_0 \supface \undx_1 \supface \cdots \supface \undx_m = \underline 0 \]
			be any sequence of removals bringing $\undq$ to $\underline 0$.
			
			By~\ref{def:axiom_independents_1}, the linear intervals $\mleft[\undq, \unda_1\mright]$ and $\mleft[\undq, \undb_1\mright]$ have the same size, that we denote by $s + 1$, with $s \geq 0$. Let
			\[ \mleft[\undq, \undb_1\mright] = \mleft\{\undb_1 = \undy_0 \supface \undy_1 \supface \cdots \supface \undy_s = \undq\mright\} \text{ and } \mleft[\undq, \undc_1\mright] = \Bigl\{\undc_1 = \undz_0 \supface \undz_1 \supface \cdots \supface \undz_s = \undq\Bigr\}. \]
			Then, $\undb_1 = \undy_0 \supface \undy_1 \supface \cdots \supface \ \undy_s = \undq = \undx_0 \supface \undx_1 \supface \cdots \supface \undx_m = \underline 0$
			is a sequence of removals that brings $\undb_1$ to $\underline 0$. Property $P(\undb_1)$ therefore implies that the length of this sequence is equal to the length of the sequence $\undb_1 \supface \cdots \supface \undb_k = \underline 0,$ that is, $k - 1 = s + m$. The same argument applied to $\undc_1$ yields $\ell - 1 = s + m$. We conclude that $k = \ell$.

			The implication \ref{ord2} $\Rightarrow$ \ref{ord1} follows from the identities 
			\[ |\undq| + |\mleft[\undq, \undp \setminus i\mright]| = |\undp| = |\undq| + |\mleft[\undq, \undp \setminus j\mright]|, \]
			using a sequence of removals of $j$ (resp. $i$) that brings $p \setminus i$ (resp. $p \setminus j$) to $q$.
		\end{proof}
	
	\subsection{Proof of the second part of Theorem~\ref{thm:axiomatic_independents}}
		
		We need the following lemma.
		
		\begin{lemma} \label{lem:increasing_and_axiom_imply_other_axiom}
			Let $\Jind$ be a subset of $\hcube{\undr}{}$ that verifies~\ref{def:axiom_independents_2}. Then, for two elements $\unda \subface \undb$ of $\Jind$ such that $D(\unda, \undb)$ has at least two elements, we have $|\unda| \leq |\undb| - 2$.
		\end{lemma}
		
		\begin{proof}
			The first part of Property~\ref{def:axiom_independents_2} ensures that $|\unda| < |\undb|$. The second part of Property~\ref{def:axiom_independents_2} now implies that there exists an element $\undc \in \Jind$ and $i \in D(a,b)$ such that $\undc \subfaceeq \unda \vee \undb = \undb$, $|\undc| > |\unda|$, and $\ssc_i < \ssb_i$. Combining the latter with $\undc \subfaceeq \undb$ yields that $\undc \subface \undb$. Applying~\ref{def:axiom_independents_2} again, we get $|\undc| < |\undb|$. All in all, we get $\unda \subface \undc \subface \undb$, and the inequality $|\unda| \leq |\undb| - 2$ follows.
		\end{proof}
		
		We now prove the second part of the main theorem.
		
		\begin{proof}[Proof of the second part of Theorem~\ref{thm:axiomatic_independents}]
			Notation as in the statement of the theorem, by Lemma~\ref{lem:second_axiom_implies_orderability}, we have a well-defined size function $|\cdot|$ on $\Ind$. We define a function $\rk$ on $\hcube{\undr}{}$ by setting
			\[ \rk(\undx) \coloneqq \max_{\substack{\unda \in \Ind \textrm{ with } \unda \subfaceeq \undx}}	\, |\unda| \qquad \forall \,\, \undx \in \hcube{\undr}{}, \]
			and show that $\rk$ is the rank function of a matricube. Note that by~\ref{def:axiom_independents_2}, $\rk(\unda) = |\unda|$ for $\unda \in \Ind$.
			
			Obviously, $\rk(\underline 0) = 0$. Moreover, by orderability of $\Ind$, for $1 \leq i \leq \grdeg $ and $1 \leq t \leq \ssr_i$, we have either $\rk(t \, \unde_i) - \rk((t - 1) \, \unde_i) = 0$ or $\rk(t \unde_i) - \rk((t - 1) \unde_i) = 1$. Therefore~\ref{def:axiom_rank_function_1} holds.
			
			Since $|\cdot|$ is increasing by~\ref{def:axiom_independents_2}, $\rk$ is non-decreasing on $\hcube{\undr}{}$. That is, Property~\ref{def:axiom_rank_function_2} holds.
			
			We show $\rk$ is submodular. Using Theorem~\ref{thm:eq_submodularity_diamond_property} proved in Section~\ref{sec:diamond_property}, it will be enough to show that $\rk$ verifies the diamond property for functions on the hypercuboid.
			
			We first observe that, by orderability of $\Ind$, for every $\undx \in \hcube{\undr}{}$ and $i \in \{1, \dots, \grdeg\}$, we have $\rk(\undx + \unde_i) - \rk(\undx) \leq 1$, provided that $\undx + \unde_i \in \hcube{\undr}{}$.
			
			Now let $\undx \in \hcube{\undr}{}$ and let $i \neq j$ be elements of $ \in \{1, \dots, \grdeg\}$ such that $\undx + \unde_i + \unde_j \in \hcube{\undr}{}$. Let $\undy \coloneqq \undx + \unde_i$, $\undz \coloneqq \undx + \unde_j$ and $\undw \coloneqq \undx + \unde_i + \unde_j$. Proving the diamond property for $\undx, \undy, \undz, \undw$ reduces to showing that the situation where $\rk(\undx) = \rk\mleft(\undy\mright) = \rk(\undz)$ and $\rk(\undw) = \rk(\undx) + 1$ never happens. For the sake of a contradiction, assume that we are in the situation where the above equalities hold. This implies in particular that $\undy, \undz \notin \Ind$. The rest of the argument is a case-by-case analysis. We first treat the case $\undw \in \Ind$, then, $\undw \notin \Ind$ but $\undx \in \Ind$, and then generalize the argument to treat the remaining case $\undw, \undx \notin \Ind$.
			
			First consider the case where $\undw \in \Ind$. Let $\unda \subfaceeq \undx$ be an element of $\Ind$ such that $\rk(\unda) = \rk(\undx)$. Since $\unda \subface \undw$ and $|D(\unda, \undw)| \geq |D(\undx, \undw)| = 2$, applying Lemma~\ref{lem:increasing_and_axiom_imply_other_axiom}, we get the inequality $\rk(\undx) \leq \rk(\undw) - 2$, which is a contradiction. This implies that $\undw \not\in \Ind$.
			
			At this point, we have deduced $\undy, \undz,\undw \notin \Ind$.
			Now consider the case where $\undx \in \Ind$. Let $\undb \subface \undw$ be an element of $\Ind$ such that $\rk(\undw) = \rk(\undb)$. Notice that $|\undb| = \rk(\undw) > \rk(\undx) = |\undx|$. Moreover, $\ssb_i = \ssw_i = \ssx_i + 1$ because otherwise we would have $\undb \subfaceeq \undy$ and $\rk(\undb) > \rk\mleft(\undy\mright)$, which would be impossible since $\rk$ is non-decreasing. Likewise, we have $\ssb_j = \ssw_j = \ssx_j + 1$. Since $\undb \subface \undw = \undx + \unde_i + \unde_j$, this shows that $D(\undx, \undb) = \{i, j \}$. By~\ref{def:axiom_independents_2}, there exists an independent $\undc \in \Ind$ such that $\undc \subfaceeq \undx \vee \undb = \undw$, $|\undc| > |\undx|$ and $\ssc_k < \ssb_k$ for some $k \in \{i, j\}$. But if $k = i$, then $\undc \subfaceeq \undy$ and therefore $|\undc| \leq \rk\mleft(\undy\mright) = \rk(\undx) = |\undx|$, a contradiction; we conclude similarly if $k = j$. We have shown that $\undx \notin \Ind$.
			
			We now treat the remaining case. We define a finite procedure by applying repeatedly an analogue of the preceding construction, as follows. Let $\unda \subface \undx$, $\undb \subface \undw$ be elements of $\Ind$ such that $\rk(\undx) = |\unda|$ and $\rk(\undw) = |\undb|$. We have $|\undb| > |\unda|$. We claim $\ssb_i = \ssw_i > \ssx_i\geq \ssa_i$. Indeed, otherwise, we would have $\undb \subfaceeq \undy$, impossible by the inequality $|\undb| = \rk(\undw) > \rk\mleft(\undy\mright)$. Likewise, we have $\ssb_j = \ssw_j > \ssa_j$. Consequently, we have $D(\unda, \undb) \supseteq \{i, j\}$. By~\ref{def:axiom_independents_2}, there exists an independent $\undc^1 \in \Ind$ such that 
			\[
				\undc^1 \subfaceeq \unda \vee \undb \subfaceeq \undw, \qquad |\undc^1| > |\unda|, \qquad \textrm{and} \qquad \ssc_{k_1}^1 < \ssb_{k_1} \quad \textrm{for some } \ssk_1 \in D(\unda, \undb).
			\]
		 	Next, if $D(\unda, \undc^1)$ contains itself at least two elements, since we have $|\undc^1| > |\unda|$, we can apply~\ref{def:axiom_independents_2}, and the same procedure as above, replacing the pair $\unda, \undb$ by the pair $\unda, \undc^1$, yields an element $\undc^2 \in \Ind$ such that
			\[
				\undc^2 \subfaceeq \unda \vee \undc^1 \subfaceeq \undw, \qquad |\undc^2| > |\unda|, \qquad \textrm{ and } \qquad \ssc_{k_2}^2 < \ssc^1_{k_2} \quad \textrm{ for some } \ssk_2 \in D(\unda, \undc^1).
			 \] 
			 
			 Repeating the procedure while it is possible, we get a sequence $\undc^1, \undc^2, \ldots, \undc^j, \ldots $ of elements of 
			 $\Ind$, satisfying for every $j \geq 1$,
			 \[
				\undc^j \subfaceeq \unda \vee \undc^{j - 1} \subfaceeq \undw, \qquad |\undc^j| > |\unda|, \qquad \textrm{ and } \qquad \ssc_{k_j}^j < \ssc^{j - 1}_{k_j} \quad \textrm{ for some } \ssk_j \in D(\unda, \undc^{j - 1}),
			\]
			with $\undc^0 = \undb$. We claim that this sequence is necessarily finite. Indeed, we observe that for every $j \geq 1$, we have by construction:
			\[ D(\unda, \undc^{j - 1}) \supseteq D(\unda, \undc^{j}) \quad \text{and} \quad E(\unde, \undc^{j - 1}) > E(\unde, \undc^{j}). \]
			Since the integers $E(\unda, \undc^j)$ are all non-negative, we infer that the sequence $\undc^\bullet$ ends at some integer $j > 0$. This means that the condition $|D(\unda, \undc^{j})| \geq 2$ fails, and thus it is impossible to have both $i$ and $j$ included in $D(\unda, \undc^{j})$.
			
			Without loss of generality, assume $i \notin D(\unda, \undc^{j})$. This implies that $\ssc^j_i \leq \ssa_i$, and so, we have $\unda^{j} \subfaceeq \undy$. We infer that $|\unda| < |\undc^{j}| \leq \rk\mleft(\undy\mright) = \rk(\undx) = |\unda|$, which is a contradiction.
			
			At this point, we have shown the diamond property, and therefore we conclude that $\rk$ is submodular, and~\ref{def:axiom_rank_function_3} follows.
			
			It follows that $\rk$ is the rank function of a matricube $\matricube$. Moreover, by definition of the rank function, $\Ind$ coincides with the set of independents of $\matricube$.
			
			Finally, by definition of $\rk$, Property~\ref{def:axiom_rank_function_1nd} is seen to be equivalent to~\ref{def:axiom_independents_star}, and so $\matricube$ is simple if, and only if,~\ref{def:axiom_independents_star} holds.
		\end{proof}

%%%Diamond property	

\section{Diamond property for functions} \label{sec:diamond_property}
 	
	The aim of this section is to generalize to the setting of matricubes the well-known result in matroid theory that the submodularity of the rank function of a matroid is equivalent to the diamond property for its graded lattice of flats~\cite[Proposition 3.3.2]{Stanley}. To this end, we here introduce a weaker version of submodularity.
	
	\begin{defi}[Diamond property for functions on hypercuboids] \label{def:diamond_property}
		We say an integer-valued function $\rk$ on $\hcube{\undr}{}$ satisfies the \emph{diamond property} if the following holds. For every point $\undx \in \hcube{\undr}{}$ and distinct integers $i \neq j \in \{1, \dots, \grdeg\}$ such that $\undx + \unde_i, \undx + \unde_j \in \hcube{\undr}{}$, we have
		\begin{equation} \label{eq:diamond}
			\rk(\undx + \unde_i) - \rk(\undx) \geq \rk(\undx + \unde_i + \unde_j) - \rk(\undx + \unde_j). \qedhere
		\end{equation}
	\end{defi}
	
	The following theorem shows that the above property is equivalent to submodularity.
	
	\begin{thm}[Equivalence of submodularity and the diamond property] \label{thm:eq_submodularity_diamond_property}
		Let $\rk$ be an integer-valued function on $\hcube{\undr}{}$. The following properties are equivalent:
		\begin{enumerate}[label=(\roman*)]
			\item \label{thm-item:weak_submodularity_simplified1} $\rk$ is submodular.
			
			\item \label{thm-item:weak_submodularity_simplified} $\rk$ verifies the diamond property.	
		\end{enumerate}
	\end{thm}
	
	In preparation for the proof, we provide generalizations of the diamond property, that allow to proceed by induction. We say an integer-valued function $\rk$ on $\hcube{\undr}{}$ satisfies the \emph{unidirectional submodularity at distance one} if the following holds. For all $i \in \{1, \dots, \grdeg\}$ and for all points $\undx \subfaceeq \undy \in \hcube{\undr}{}$ such that $\ssx_i = \ssy_i$ and $\undx + \unde_i \in \hcube{\undr}{}$, the inequality
	\begin{equation} \label{eq:weak_submodularity}
		\rk(\undx + \unde_i) - \rk(\undx) \geq \rk\mleft(\undy + \unde_i\mright) - \rk\mleft(\undy\mright)
	\end{equation}
	holds. More generally, we have the following generalization of Property~\eqref{eq:weak_submodularity} in \emph{several directions} and at \emph{higher distance}.
	
	\begin{defi}[Multidirectional submodularity at a given distance] \label{def:gen_weak_submodularity_property}
		For positive integers $k$ and $n$, we define the \emph{$k$-directional submodularity at distance up to $n$}, denoted $(*)^n_k$, as follows.
 		
		\smallskip
		
		\noindent $(*)^n_k$: \, Pick any integer $1 \leq s \leq k$, any integers $1 \leq \ssi_1 < \cdots < \ssi_s \leq \grdeg$ and $0 \leq \ssn_{\ssi_1}, \dots, \ssn_{\ssi_s} \leq n$.
		
		\smallskip
		
		\noindent Then, for every pair of elements $\undx \subfaceeq \undy \in \hcube{\undr}{}$ such that $\ssx_{\ssi_j} = \ssy_{\ssi_j}$ for all $1 \leq j \leq s$, and $\undx + \sum_{1 \leq j \leq s} \ssn_j \, \unde_{\ssi_j} \in \hcube{\undr}{}$, we have
		\[ \rk\Bigl(\undx + \sum_{j = 1}^s \ssn_{\ssi_j} \, \unde_{\ssi_j}\Bigr) - \rk (\undx) \geq \rk\Bigl(\undy + \sum_{j = 1}^s \ssn_{\ssi_j} \, \unde_{\ssi_j}\Bigr) - \rk\mleft(\undy\mright). \qedhere \]
	\end{defi}
	
	Notice that the property stated in~\eqref{eq:weak_submodularity} is exactly $(*)^1_1$ as defined above, and the terminologies are consistent. Moreover, any $(*)^n_k$ with $k, n \geq 1$ implies $(*)^1_1$.
	
	\begin{remark}[Alternative description of $(*)^n_k$]
		Using the notation of Definition~\ref{def:gen_weak_submodularity_property}, after the change of variables $\unda \coloneqq \undx, \; \undb \coloneqq \undy - \sum_{j = 1}^s \ssn_{\ssi_j} \, \unde_{\ssi_j}$, property $(*)^n_k$ can be rewritten as follows.
		
		For all elements $\unda$ and $\undb \in \hcube{\undr}{}$, we have the submodularity inequality
		\[ \rk(\unda) + \rk(\undb) \geq \rk(\unda \vee \undb) + \rk(\unda \wedge \undb) \]
		as long as there exist an integer $1 \leq s \leq k$ and integers $1 \leq \ssi_1 < \cdots < \ssi_s \leq \grdeg$ such that $\undb + \sum_{1 \leq j \leq s} (\ssa_{\ssi_j} - \ssb_{\ssi_j}) \, \unde_{\ssi_j}$ is an element of $\hcube{\undr}{}$ greater than or equal to $\unda$ and, such that, for all $1 \leq j \leq s$, we have $0 \leq \ssa_{\ssi_j} - \ssb_{\ssi_j} \leq n$.
		
		This parametrization using $\unda$ and $\undb$ enables to see instantaneously that the submodularity property of $\rk$ in the hypercuboid implies all the properties $(*)^n_k$. The other parametrization, using $\undx$ and $\undy$, will be useful to prove Theorem~\ref{thm:eq_submodularity_diamond_property} below, in that it behaves linearly (contrary to formulas involving the symbols $\wedge$ and $\vee$).
	\end{remark}
	
	\begin{proof}[Proof of Theorem~\ref{thm:eq_submodularity_diamond_property}]
		Obviously, \ref{thm-item:weak_submodularity_simplified1} implies \ref{thm-item:weak_submodularity_simplified}.
		
		\smallskip
		
		We explain how to deduce \eqref{eq:weak_submodularity}, that is $(*)^1_1$, from \ref{thm-item:weak_submodularity_simplified}. Let $\undx$ and $\undy$ be as in Definition~\ref{def:diamond_property}. The fact that $\ssx_i = \ssy_i$ implies that $\undy$ can be written as $\undy = \undx + \sum_{j \neq i} \ssn_j \, \unde_j$ with $\ssn_j \geq 0$, and we can sum inequalities of the form (\ref{eq:diamond}) to get the inequality~(\ref{eq:weak_submodularity}).
		
		\smallskip
		
		We then explain how to deduce \ref{thm-item:weak_submodularity_simplified1} from $(*)^1_1$. Proceeding by induction, we show that the property $(*)^1_1$ implies $(*)^n_k$ for all $k, n \geq 1$. We first show that $(*)^1_1$ implies $(*)^n_1$ for all $n \geq 1$. Let $i \in \{1, \dots, \grdeg\}$ and $0 \leq \ssn_i \leq n$, and let $\undx \subfaceeq \undy$ be elements of $\hcube{\undr}{}$ such that $\undx + \ssn_i \, \unde_i \in \hcube{\undr}{}$ and $\ssx_i = \ssy_i$. For all $0 \leq t < n_i$, the pair $\mleft(\undx + t \, \unde_i, \undy + t \, \unde_i\mright)$ satisfies the hypotheses needed to apply $(*)^1_1$ in direction $i$, so we know that
		\[ \rk(\undx + (t + 1) \, \unde_i) - \rk(\undx + t \, \unde_i) \geq \rk\mleft(\undy + (t + 1) \, \unde_i\mright) - \rk\mleft(\undy + t \, \unde_i\mright). \]
		Summing all these inequalities for $0 \leq t < \ssn_i$, and canceling out the terms which appear on both sides, yields
		\[ \rk(\undx + \ssn_i \, \unde_i) - \rk(\undx) \geq \rk\mleft(\undy + \ssn_i \, \unde_i\mright) - \rk\mleft(\undy\mright), \]
		which gives $(*)^n_1$.
		
		We now show that properties $(*)^n_1$ for $n \geq 1$ imply properties $(*)^n_2$. Let $i, j \in \{1, \dots, \grdeg\}$ and $0 \leq \ssn_i, \ssn_j \leq n$, and let $\undx \subfaceeq \undy$ be elements of $\hcube{\undr}{}$ such that $\undx + \ssn_i \, \unde_i + \ssn_j \, \unde_j \in \hcube{\undr}{}$, $\ssx_i = \ssy_i$ and $\ssx_j = \ssy_j$. We apply $(*)^n_1$ to the pair $\mleft(\undx, \undy\mright)$ in direction $i$ and get
		\[ \rk(\undx + \ssn_i \, \unde_i) - \rk(\undx) \geq \rk\mleft(\undy + \ssn_i \, \unde_i\mright) - \rk\mleft(\undy\mright). \]
		The pair $\mleft(\undx + \ssn_i \, \unde_i, \undy + \ssn_i \, \unde_i\mright)$ satisfies the hypotheses required for applying $(*)^n_1$ again, but this time in direction $j$. This yields
		\[ \rk(\undx + \ssn_i \, \unde_i + \ssn_j \, \unde_j) - \rk(\undx + \ssn_i \, \unde_i) \geq \rk\mleft(\undy + \ssn_i \, \unde_i + \ssn_j \, \unde_j\mright) - \rk\mleft(\undy + \ssn_i \, \unde_i\mright). \]
		Summing up these two inequalities shows that $\rk$ satisfies $(*)^n_2$. The same procedure inductively proves that $\rk$ satisfies all $(*)^n_k$, i.e., $\rk$ is submodular.
	\end{proof}
	
	\begin{remark}[Discrete partial derivatives and transverse local convexity]
		For $i \in \{1, \dots, \grdeg\}$, we can define the \emph{discrete partial derivative of $\rk$ in the direction $i$} as the function $\partial_i \rk$ defined by
		\[ \partial_i \rk(\undx) \coloneqq \rk(\undx + \unde_i) - \rk(\undx), \qquad \forall \undx \in \hcube{\undr}{} \text{ such that } \ssx_i < \ssr_i. \]
		We notice that property $(*)^1_1$ is equivalent to the fact that for all $i \in \{1, \dots, \grdeg\}$ and for all $0 \leq t < \ssr_i$, $\partial_i \rk\rest{\ssL_t^i}$ is non-increasing. This is why $(*)^1_1$ may be alternatively called \emph{transverse local concavity}. Submodularity is thus equivalent to transverse local concavity.
		
		In other contexts, submodularity is sometimes referred to as the discrete analogue of concavity: see, for example, \cite[Theorem~44.1]{schrijver2003combinatorial}. While this is fully relevant for a function $\rk$ defined on the collection $\mathcal{P}(S)$ of all subsets of a given set $S$, that is on the hypercube $\hcube{1}{\grdeg}$, it is not exactly true for supermodular functions on $\hcube{\undr}{}$ for larger values of $\ssr_1, \dots, \ssr_{\grdeg}$. This is because the functions $\partial_i \rk$ are non-decreasing only in directions different from $i$. For the function $\rk$ defined on $\hcube{(2,2)}{}$ by $ \begin{pmatrix}
			\color{red} 2 & \color{red} 2 & \color{red} 3 \\
			1 & 2 & 3 \\
			0 & 1 & 2 \\
		\end{pmatrix}$, we have $\rk((2, 1)) - \rk((2, 0)) \ngeq \rk((2, 2)) - \rk((2, 1))$, i.e., $\ssub\partial!_2 \rk$ is not non-increasing in direction $2$.
	\end{remark}

\section{Permutation arrays} \label{sec:permutation_arrays}
	
	The aim of this section is to study simple matricubes with ground set an actual hypercube ($\ssr_1 = \ssr_2 = \cdots = \ssr_\grdeg = r$) of minimum possible rank $r$ or $r + 1$. In the representable case, this corresponds to a collection of $\grdeg$ complete flags in a vector space of dimension $r + 1$. Theorem~\ref{thm:equiv_perm_supermodular} establishes a one-to-one correspondendance between these matricubes and permutation arrays introduced by Eriksson-Linusson ~\cite{eriksson2000combinatorial, eriksson2000decomposition}.	
	
	\subsection{Permutation arrays}
		
		First, we recall some terminology from~\cite{eriksson2000combinatorial}. Our presentation differs slightly from the original setting as our indexing of flags is by codimension while in their work, Eriksson and Linusson use an indexing by dimension. (Concretely, this amounts to having lower blocks in~\cite{eriksson2000combinatorial, eriksson2000decomposition} replaced here by upper blocks.)
		
		Let $\ssr_1, \dots, \ssr_\grdeg$ be $\grdeg$ non-negative integers. An $\grdeg$-\emph{dimensional dot array} $P$ is an $\grdeg$-dimensional array of type $[\ssr_1] \times \dots \times [\ssr_\grdeg]$ where some of the entries are dotted.
		
		For a dot array $P$, and $\undx \in \hcube{\undr}{}$, we denote by $P[\undx]$ the \emph{upper principal subarray} of $P$, which consists of all $\undy$ with $\undy \supfaceeq \undx$. It is naturally a dot array itself.
		
		To be precise, for $P[\undx]$ to become a dot array, we must coordinate-wise subtract the point $(\ssx_1, \dots, \ssx_\grdeg)$ to all its elements. In the following, we will use both parametrization conventions freely for the sake of convenience.
		
		For a dot array $P$ and $j \in \{1, \dots, \grdeg\}$, the \emph{rank along the $j$-axis}, denoted by $\rank_j(P)$, is the total number of $0 \leq t \leq \ssr_j$ such that there is at least one dot in some position whose $j$-th index is equal to $t$, i.e., there is at least one dot in the layer $\ssL_t^j$ of $P$. A dot array $P$ is called \emph{rankable} if we have $\rank_j(P) = \rank_{i}(P)$ for all $i, j \in \{1, \dots, \grdeg\}$. If $P$ is rankable, then we call $\rank_j(P)$ the \emph{rank} of $P$ for any $j \in \{1, \dots, \grdeg\}$.
		
		A dot array $P$ is called \emph{totally rankable} if every upper principal subarray of $P$ is rankable.
		
		We recall that in the terminology of \cite{eriksson2000combinatorial} and \cite{eriksson2000decomposition}, a position $\undx$ is \emph{redundant} if there exist dot positions $\undy_1, \dots, \undy_m \neq \undx$, for some $m \geq 2$, such that each $\undy_i$ has at least one coordinate in common with $\undx$, and such that $\undx = \bigwedge_{i = 1}^m \undy_i$. The set of redundant positions of $P$ is denoted by $R(P)$. A \emph{redundant dot} is a redundant position that is dotted. The reason for the term ``redundant'' is that placing or removing a redundant dot does not change the rank of any upper principal subarray of $P$. (In the language of lattices, a non-redundant position is meet-irreducible in the set of dotted positions.)
		
		If $A$ is a subset of $\hcube{\undr}{}$, then $P \cup A$ (resp. $P \setminus A$) denotes the dot array based on $P$ where, for every $\undx \in A$, we dot (resp. undot) the position $\undx$ in $P$, if necessary.
		
		A \emph{permutation array} of width $r$ and dimension $\grdeg$ is a totally rankable dot array $P$ of shape $\hcube{r}{\grdeg} = \hcube{\undr}{} = [r]^\grdeg$, $\undr = (r, \dots, r)$, of rank $r + 1$, and with no redundant dots.
	
	\subsection{Equivalence of permutation arrays with simple matricubes of rank $r$ or $r + 1$ on $\hcube{r}{\grdeg}$}
		
		Our next theorem establishes an equivalence between permutation arrays and simple matricubes of rank $r$ or $r + 1$ on the hypercube $\hcube{r}{\grdeg}$.
		
		\begin{thm} \label{thm:equiv_perm_supermodular}
			Let $P$ be a permutation array of width $r$ and dimension $\grdeg$. The function $\rk_P$ defined by $\rk_P(\unda) \coloneqq r + 1 - \rank(P[\unda])$ for every $\unda \in \hcube{r}{\grdeg}$ is the rank function of a simple matricube $\matricube_{P}$ with ground set $\hcube{r}{\grdeg}$. This matricube is of rank $r$ or $r + 1$ depending on whether the position $\undr$ in $\hcube{r}{\grdeg}$ is dotted or not. The set of flats of $\rk_P$ is precisely the union of the set of dot positions in $P$ with $R(P)$, and $\undr$.
			
			Conversely, the rank function $\rk$ of every simple matricube $\matricube$ of rank $r$ or $r + 1$ on the hypercube $\hcube{r}{\grdeg}$ defines a dot array $\ssP{\matricube}$ on $\hcube{r}{\grdeg} = [r]^\grdeg$ with dots positioned on the set of flats $\unda \neq \undr$ of $\matricube$, and also a dot positioned on $\undr$ if $\rk(\matricube) = r$. Then, $P \coloneqq \ssP{\matricube} \setminus R(\ssP{\matricube})$ is a permutation array.
		\end{thm}
		
		The proof of this theorem is given in the next section.
	
	\subsection{Proof of Theorem~\ref{thm:equiv_perm_supermodular}} \label{subsec:equiv_perm_supermodular}
 		
 		We start by proving the first part of the theorem. Let $P$ be a permutation array on $\hcube{r}{\grdeg} = [r]^\grdeg$. We claim that the function
		\[ \ssrho_P(\undx) \coloneqq r + 1 - \rank(P[\undx]), \qquad \forall \undx \in \hcube{r}{\grdeg}, \]
		is the rank function of a simple matricube on the ground set $\hcube{r}{\grdeg}$. We need to show properties \ref{def:axiom_rank_function_1nd}-\ref{def:axiom_rank_function_2}-\ref{def:axiom_rank_function_3}.
		
		Since $\undx \subfaceeq \undy$ implies $P[\undx] \supseteq P\mleft[\undy\mright]$, we deduce that $\rk_P$ is non-decreasing, which shows \ref{def:axiom_rank_function_2}.
		
		We now prove \ref{def:axiom_rank_function_1}. Let $i \in \{1, \dots, \grdeg\}$. We have to show that $\rk_P(t \, \unde_i) = t$ for $t \in [r]$. By definition, $\rank(P[t \, \unde_i]) \leq r + 1 - t$, which implies $\rk_P(t \, \unde_i) \geq t$.
		The reverse inequality is shown by induction on $t$. The case $t = 0$ is true by the definition of permutation arrays, which requires $\rank(P) = r + 1$. Assuming that $\rank(P[t \, \unde_i]) \geq r + 1 - t$, we show that $\rank(P[(t + 1) \, \unde_i]) \geq r - t$. This follows from the inequality $\rank(P[\unda + \unde_i]) = \rank_i(P[\unda + \unde_i]) \geq \rank_i(P[\unda]) - 1 = \rank(P[\unda]) - 1$, valid for every $\unda \in \hcube{r}{\grdeg}$ such that $\unda + \unde_i \in \hcube{r}{\grdeg}$.
		
		It remains to show that $\rk_P$ is submodular. Thanks to Theorem~\ref{thm:eq_submodularity_diamond_property}, it is sufficient to show that $\rk_P$ satisfies the diamond property for functions. We thus take two distinct integers $i \neq j \in \{1, \dots, \grdeg\}$ and an element $\undx \in \hcube{r}{\grdeg}$ such that $\undx + \unde_i, \undx + \unde_j \in \hcube{r}{\grdeg}$. We assume that $\rk(\undx + \unde_i + \unde_j) - \rk(\undx + \unde_j) = 1$ and show that $\rk(\undx + \unde_i) - \rk(\undx) = 1$. The hypothesis implies that the layer $\ssL_{\ssx_i}^i$ in the dot array $P[\undx + \unde_j]$ contains a dotted point. This point will be counted in the difference $\rk(\undx + \unde_i) - \rk(\undx)$, which proves the result.
		
		The matricube $\matricube_P$ is of rank $r$ or $r + 1$ depending on whether $\rank\mleft(P\mleft[\undr\mright]\mright) = 1$ or $0$, that is, whether $\undr$ is dotted or not.
		
		Finally, to see the statement about the flats, consider $\undx \neq \undr$ and assume first that $\undx$ is dotted. Then, for each direction $\unde_i$ with $\undx + \unde_i \in \hcube{r}{\grdeg}$, we get $\rank_i(P[\undx]) - \rank_i(P[\undx + \unde_i]) = 1$. This shows that $\undx$ is a flat.
		
		Next, assume that $\undx$ is not dotted. Since flats are closed under meet, if $\undx$ is a redundant position, then it is a flat. It remains to consider the case where $\undx$ is neither dotted nor a redundant point. This means there is an $i \in \{1, \dots, \grdeg\}$ such that the layer $\ssL^i_{x_i}$ in $P[\undx]$ does not contain any dot. Two cases can happen:
		
		\begin{itemize}
			\item If $\undx + \unde_i \in \hcube{r}{\grdeg}$, then 
			\[ \rank(P[\undx]) = \rank_i(P[\undx]) = \rank_i(P[\undx + \unde_i]) = \rank(P[\undx + \unde_i]), \]
			and thus $\rk(\undx) = \rk(\undx + \unde_i)$, and $\undx$ is not a flat of $\matricube_P$.
			
			\item Otherwise, $\ssx_i = r$, and so $\rank_i(P[\undx]) = 0$, that is, $\rk(\undx) = r + 1$. This implies that $\matricube$ is of rank $r + 1$, and since $\undx \neq \undr$, then, again $\undx$ is not a flat.
		\end{itemize}
 		This finishes the proof of the first direction.
		
		\smallskip
 		
		We now show the other direction. Suppose that $\matricube$ is a simple matricube of rank $r + 1$ or $r$ on $\hcube{r}{\grdeg}$ with rank function $\rk$. Let $\ssP{\matricube}$ be the corresponding dot array where a dot is positioned on every flat $\unda$ of $\matricube$ different from $\undr$, and if the rank of $\matricube$ is $r$, then a dot is also positioned on $\undr$. Let $P = \ssP{\matricube} \setminus R(\ssP{\matricube})$. We show that $P$ is a permutation array.
		
		By construction, $P$ has no redundant dots. We thus need to show that $P$ is totally rankable and has rank $r + 1$. We have to prove that for every $\undx \in \ssP{\matricube}$ and $i, j \in \{1, \dots, \grdeg\}$, $\rank_i(\ssP{\matricube}[\undx]) = \rank_j(\ssP{\matricube}[\undx])$. This is a direct consequence of Proposition~\ref{prop:link_rank_rho} below, which also shows that the rank of $\ssP{\matricube}$ is $r + 1$. We conclude that $P$ is a permutation array. \qed
		
		\begin{prop} \label{prop:link_rank_rho}
			Suppose that $\matricube$ is a simple matricube of rank $r$ or $r + 1$ on the ground set $\hcube{r}{\grdeg}$ and denote by $\rk$ its rank function. Let $\ssP{\matricube}$ be the corresponding dot array with a dot positioned at each flat $\unda \neq \undr$, and also a dot positioned at $\undr$ in the case $\rk(\matricube) = r$.
			
			Let $\undx$ be an element of the dot array $\ssP{\matricube}$ and $1 \leq i \leq \grdeg$. Then, $\mathrm{rank}_i(\ssP{\matricube}[\undx]) = r + 1 - \rk(\undx)$.
		\end{prop}
		
		\begin{proof}
			We proceed by reverse induction in the lattice $\hcube{r}{\grdeg}$, starting from $\undr$. For $\undr$, we have $\rank_i\mleft(\ssP{\matricube}\mleft[\undr\mright]\mright) = 0$ or $1$ depending on whether $\rk(\matricube) = r + 1$ or $r$, respectively, for each $i \in \{1, \dots, \grdeg\}$, as required.
			Assume $\undx \neq \undr$. We suppose the following equalities hold:
			\[ \rank_i\mleft(\ssP{\matricube}\mleft[\undy\mright]\mright) = r + 1 - \rk\mleft(\undy\mright) \quad \forall \, \, \undy \gneq \undx \qquad \textrm{and} \qquad \forall \, \, i \in \{1, \dots, \grdeg\}. \] 
			We prove that the equalities hold as well for $\undx$.
			
			\smallskip
			
			Suppose first that $\undx$ is a flat. Two cases can occur.
			
			\begin{itemize}
				\item[(I.1)] If $\undx$ has rank $r$, then for each $i \in \{1, \dots, \grdeg\}$ with $\undx + \unde_i \in \hcube{r}{\grdeg}$, we have $\rk(\undx + \unde_i) = r + 1$. By induction, $\rank_i(\ssP{\matricube}[\undx + \unde_i]) = r + 1 - \rk(\undx + \unde_i) = 0$. It follows that $\rank_i(\ssP{\matricube}[\undx]) = \rank_i(\ssP{\matricube}[\undx + \unde_i]) + 1 = 1$, as required.
				For the other values of $i$, we have $\undx + \unde_i \notin \hcube{r}{\grdeg}$, that is, $\ssx_i = r$, and, in this case, we have as well $\rank_i(\ssP{\matricube}[\undx]) = 1$.
				
				\item[(I.2)] If $\rk(\undx) < r$, then using the inequality $\ssx_i = \rk(\ssx_i \unde_i) \leq \rk(\undx)$, we get $\ssx_i < r$. This implies that $\undx + \unde_i \in \hcube{r}{\grdeg}$ for all $i$. We get
				\[ \rank_i(\ssP{\matricube}[\undx]) = \rank_i(\ssP{\matricube}[\undx + \sse_i]) + 1 = r + 1 - \rk(\undx + \sse_i) + 1 = r + 1 - \rk(\undx), \]
				as required.
			\end{itemize}
			
			Now suppose that $\undx$ is not a flat of $\matricube$. Again, two cases can occur.
			
			\begin{itemize}
				\item[(II.1)] If $\ssx_i < r$ for all $i \in \{1, \dots, \grdeg\}$, then let $\unda$ be the minimum flat with $\unda \supfaceeq \undx$ and let $i \in \{1, \dots, \grdeg\}$. Applying the induction hypothesis to $\undy = \undx + \unde_i$, we get $\rank_i\mleft(\ssP{\matricube}\mleft[\undy\mright]\mright) = r + 1 - \rk\mleft(\undy\mright)$. By Lemma~\ref{lem:layer_withouh_dots}, we deduce that
				\[
					\rk\mleft(\undy\mright) = \begin{cases}
						\rk(\undx) + 1 & \textrm{ if } \ssa_i = \ssx_i, \\
						\rk(\undx) & \textrm{ if } \ssa_i >\ssx_i.
					\end{cases}	
				\]
				In the first case, when $\ssa_i = \ssx_i$, there is a dot in the layer $\ssL^i_{\ssx_i}$ of $\ssP{\matricube}[\undx]$, and thus $\rank_i(\ssP{\matricube}[\undx]) = \rank\mleft(\ssP{\matricube}\mleft[\undy\mright]\mright) + 1$. In the second case, when $\ssa_i > \ssx_i$, there is no dot in the layer $\ssL^i_{\ssx_i}$ of $\ssP{\matricube}[\undx]$, and therefore $\rank_i(\ssP{\matricube}[\undx]) = \rank\mleft(\ssP{\matricube}\mleft[\undy\mright]\mright)$. We infer that $\rank_i(\ssP{\matricube}[\undx]) = r + 1 - \rk(\undx)$, as required.
				
				\item[(II.2)] We now treat the remaining case where $\ssx_i = r$ for some indices $i$ among $1, \dots, \grdeg$. In this case, $\rk(\undx)$ is either $r$ or $r + 1$. We treat each of these possibilities separately.
				
				\begin{enumerate}
					\item Firstly, suppose that $\rk(\undx) = r$. Consider an index $i$ with $\ssx_i = r$. Then, if $\rk(\matricube) = r$, there is a dot positioned at $\undr$, and thus $\rank_i(\ssP{\matricube}[\undx]) = 1 = r + 1 - \rk(\undx)$. If $\rk(\matricube) = r + 1$, the minimum flat $\unda \supfaceeq \undx$ has rank $r$, and lives in the layer $\ssL^i_r$ of $\ssP{\matricube}[\undx]$. Again, we get $\rank_i(\ssP{\matricube}[\undx]) = 1$.
					
					Now consider an index $j$ with $\ssx_j < r$, so that $\undx + \unde_j \in \hcube{r}{\grdeg}$. Let $\unda$ be the minimum flat dominating $\undx$.
					A reasoning similar to (II.1), based on the use of Lemma~\ref{lem:layer_withouh_dots}, shows that $\rank_i(\ssP{\matricube}[\undx]) = r + 1 - \rk(\undx)$, as required.
					
					\item Secondly, suppose that $\rk(\undx) = r + 1$. The unique flat $\unda$ that dominates $\undx$ is $\undr$, which is not dotted. Therefore, $\rank_i(\ssP{\matricube}[\undx]) = 0 = r + 1 - \rk(\undx)$, as required. \qedhere
				\end{enumerate}
			\end{itemize}
		\end{proof}

\section{Local matroids} \label{sec:local_matroids}
	
	In this section, we define local matroids of matricubes, and formulate a local obstruction to their representability. We then turn this into an equivalent characterization of matricubes.
	
	\subsection{Local matroids of a matricube}
		
		In the following, for all $\unda \in \hcube{\undr}{}$, we define $\ssI_{\unda}$ as the set of all $i \in \{1, \dots, \grdeg\}$ such that $\unda + \unde_i \in \hcube{\undr}{}$.
		
		\smallskip
		
		To motivate the definition, first consider a representable matricube, given by $\grdeg$ initial flags $\ssfilt_1^\bullet, \dots, \ssfilt_\grdeg^\bullet$ of length $\ssr_1, \dots, \ssr_\grdeg$, respectively, in a $\k$-vector space $\VS$. Let $\ssfilt^{\unda} \coloneqq \bigcap_{1 \leq i \leq \grdeg} \ssfilt_i^{\ssa_i}$. The arrangement of subspaces $\ssfilt^{\unda + \unde_i}$ in $\ssfilt^{\unda}$ given by $i\in \ssI_{\unda}$ defines a matroid $\mat_{\unda}$ on the ground set $\ssI_{\unda}$. The rank one elements of this matroid correspond to those $i$ with $\ssfilt^{\unda + \unde_i}$ a proper vector subspace of $\ssfilt^{\unda}$; all the other elements are loops.
		
		This picture generalizes to any matricube, and defines local matroids associated to elements of the hypercuboid $\hcube{\undr}{}$. Let $\rk$ be a rank function on $\hcube{\undr}{}$. Let $\unda \in \hcube{\undr}{}$, and define a function $\ssrho_{\unda} \colon 2^{\ssI_{\unda}} \rightarrow \Z_{\geq 0}$ as follows. For every subset $X \subseteq \ssI_{\unda}$, set
		\[ \ssrho_{\unda}(X) \coloneqq \rk\Bigl(\unda + \sum_{i \in X} \unde_i\Bigr) - \rk(\unda). \]
		
		\begin{prop} \label{prop:local_matroids}
			The pair $(\ssI_{\unda}, \ssrho_{\unda})$ defines a matroid $\mat_{\unda}$ on the set of elements $\ssI_{\unda}$.
		\end{prop}
		
		\begin{proof}
			By Proposition~\ref{prop:inequality_one}, $\ssrho_{\unda}$ takes values in the set $\mleft\{0, \dots, |\ssI_{\unda}|\mright\}$, and $\ssrho_{\unda}(X) \leq |X|$. Since $\rk$ is non-decreasing, we also have $\ssrho_{\unda}(Y) \leq \ssrho_{\unda}(X)$ for $Y \subseteq X$. Therefore, it is enough to show that $\ssrho_a$ is submodular, that is,
			\[ \forall \, X, Y \subseteq \ssI_{\unda} \, \quad \ssrho_{\unda}(X) + \ssrho_{\unda}(Y) \geq \ssrho_{\unda}(X \cup Y) + \ssrho_{\unda}(X \cap Y). \]
			This follows from the submodularity of $\rk$ applied to $\unda + \sum_{i \in X} \unde_i$ and $\unda + \sum_{i \in Y} \unde_i$ in $\hcube{\undr}{}$.
		\end{proof}
		
		\begin{prop} \label{prop:local_obstructions}
			A necessary condition for the representability of a matricube $\matricube$ with ground set $\hcube{\undr}{}$ on a field $\k$ is the representability of all the matroids $\mat_{\unda}$, $\unda \in \hcube{\undr}{}$, on $\k$.
		\end{prop}
		
		\begin{proof}
			This follows directly from the above discussions.
		\end{proof}
	
	\subsection{The case of permutation arrays}
		
		Let $\k$ be a field. By Theorem~\ref{thm:equiv_perm_supermodular}, the representability of a simple matricube $\matricube$ of rank $r$ or $r + 1$ with the ground set the hypercube $\hcube{r}{\grdeg}$ is equivalent to the representability of the corresponding permutation array in the terminology of~\cite{eriksson2000decomposition}. Billey and Vakil~\cite{billey2008intersections} provide several examples of permutation arrays which are non-representable. Proposition~\ref{prop:local_obstructions} above provides a conceptual explanation of the examples treated in~\cite{billey2008intersections}. Theorem~\ref{thm:representability_matricube_matroid}, combined with Theorem~\ref{thm:equiv_perm_supermodular}, shows that over an infinite field, the representability for permutation arrays reduce to the representability of matroids.
	
	\subsection{Coherent complexes of matroids and matricubes} \label{sec:matroidal-characterization}
		We show that the data of a matricube on the ground set $\hcube{\undr}{}$ is equivalent to the data of a set of matroids indexed by $\hcube{\undr}{}$ satisfying compatibility properties~\ref{def:axiom_coherent_complex_1} and~\ref{def:axiom_coherent_complex_2} listed below.
		
		Notation as in the previous section, for the sake of convenience, if $i \in \ssI_{\unda}$, we write $\ssrho_{\unda}(i)$ instead of $\ssrho_{\unda}(\{i\})$. We start with the definitions below.
		
		\begin{defi}[Increasing path] \label{def:increasing_path}
			Let $\unda$ and $\undb$ be points of $\hcube{\undr}{}$ with $\unda \subfaceeq \undb$. We define an \emph{increasing path} from $\unda$ to $\undb$ to be any finite sequence
			\[ \unda = \undc_0, \undc_1, \dots, \undc_k = \undb \]
			such that for every $0 \leq j < k$, we have $\undc_{j + 1} = \undc_j + \unde_\ell$ for some $\ell \in \{1, \dots, \grdeg\}$.
		\end{defi}
		Note that the integer $k$ is equal to $\sum_{i = 1}^\grdeg (\ssb_i - \ssa_i)$.
		
		\begin{defi}[Coherent complex of matroids] \label{def:coherent_complex_matroids}
			Let $\mleft(\mat_{\unda}\mright)_{\unda \in \hcube{\undr}{}}$ be a set of matroids indexed by $\hcube{\undr}{}$, with $\mat_{\unda}$ a matroid on the set $\ssI_{\unda}$ and with rank function $\rk_{\unda}$. We say $\mleft(\mat_{\unda}\mright)$ form a \emph{coherent complex of matroids} if the following two properties are satisfied:
			
			\begin{enumerate}[label=(CC\arabic*)]
				\item \label{def:axiom_coherent_complex_1} For all $i \in \{1, \dots, \grdeg\}$ and $0 \leq t < \ssr_i$, we have $\ssrho_{t \, \unde_i}(i) \leq 1$.
				
				\item \label{def:axiom_coherent_complex_2} The matroids satisfy the following relation.
				\[ \mat_{\unda + \unde_i} =
				\begin{cases}
					\mat_{\unda} \mycontr i &\textrm{if } \ssa_i = \ssr_i - 1 \\
					\mat_{\unda} \mycontr i \sqcup \{i\} &\textrm{else}
				\end{cases}. \]
			\end{enumerate}
			
			Here, $\mat \mycontr e$ denotes the contraction of a matroid $\mat$ by its element $e$, and $i$ is the element of the matroid set corresponding to the direction $i$. Moreover, $\mat_{\unda} \mycontr i \sqcup \{i\}$ denotes an extension of $\mat_{\unda} \mycontr i$ by a single element denoted $i$.
		\end{defi}
		
		In the following, we denote by $\ssrho_{\unda \mycontr i}$ the rank function on $\ssI_{\unda} \setminus \{i\}$ that defines the matroid $\mat_{\unda} \mycontr i$. It is explicitely given by the following equation, in terms of the rank function $\ssrho_{\unda}$:
		\[ \ssrho_{\unda \mycontr i}(X) = \ssrho_{\unda}(X \cup \{i\}) - \ssrho_{\unda}(i) \qquad \textrm{ for all } X \subseteq \ssI_{\unda} \setminus \{i\}. \]
		
		\begin{remark} \label{rem:loop_coherent_complex}
			Property \ref{def:axiom_coherent_complex_2} above implies the following: let $\undx \subfaceeq \undy$ be two points of $\hcube{\undr}{}$ and $i \in \{1, \dots, \grdeg\}$ such that $\ssx_i = \ssy_i$ and $i \in \ssI_{\undx}$. Then, $i$ being a loop in $\mat_{\undx}$ implies that $i$ is a loop in $\mat_{\undy}$. Indeed, $\mat_{\undy}$ is obtained from $\mat_{\undx}$ through a sequence of operations consisting of either the contraction of an element different from $i$ or an extension. These operations do not change the property of $i$ being a loop.
		\end{remark}
		
		\begin{thm}[Matroidal characterization of matricubes] \label{thm:matroidal_characterization}
			There is a one-to-one correspondence between coherent complexes of matroids indexed by the hypercuboid $\hcube{\undr}{}$ and matricubes with ground set $\hcube{\undr}{}$.
		\end{thm}
		
		\begin{proof}
			$(\Longleftarrow)$ If we start with a matricube $\matricube$ of rank function $\rk$ on $\hcube{\undr}{}$, the collection of matroids $\mat_{\unda}$ defined above forms a coherent complex of matroids. Indeed, Property~\ref{def:axiom_coherent_complex_1} is trivially satisfied because of Property~\ref{def:axiom_rank_function_1} in Definition~\ref{def:matricube}. We check Property~\ref{def:axiom_coherent_complex_2}. Let $\unda \in \hcube{\undr}{}$ and $i \in \ssI_{\unda}$. If $\ssa_i = \ssr_i - 1$, then $\ssI_{\unda + \unde_i} = \ssI_{\unda} \setminus \{i\}$, which is the ground set of the matroid $\mat_{\unda} \mycontr i$. If $\ssa_i < \ssr_i - 1$, then $\ssI_{\unda + \unde_i} = \ssI_{\unda}$, which is the ground set of the matroid $\mat_{\unda} \mycontr i \sqcup \{i\}$. We now check the equality of the rank functions on subsets of $\ssI_{\unda} \setminus \{i\}$. Consider $X \subseteq \ssI_{\unda + \unde_i}$ not containing the element $i$. We need to show that
			\[ \rk_{\unda + \unde_i}(X) = \rk_{\unda}(X \cup \{i\}) - \rk_{\unda}(i). \]
			The left-hand side is by definition $\rk\bigl(\unda + \unde_i + \sum_{j \in X} \unde_j\bigr) - \rk(\unda + \unde_i)$, and the right-hand side is
			\[ \rk\Bigl(\unda + \unde_i + \sum_{j \in X} \unde_j\Bigr) - \rk(\unda) - \rk(\unda + \unde_i) + \rk(\unda). \]
			Both sides are therefore equal.
			
			\smallskip
			
			$(\Longrightarrow)$ The other way around, we consider a coherent complex of matroids $(\mat_{\unda})_{\unda}$ and associate a matricube $\matricube$ on $\hcube{\undr}{}$ by specifying its rank function $\rk$. Let $\unda \in \hcube{\undr}{}$. We take any increasing path $\underline 0 = \undb_0, \undb_1, \dots, \undb_k = \unda$ from $\underline 0$ to $\unda$, and define
			\[ \rk(\unda) \coloneqq \sum_{j = 0}^{k - 1} \ssrho_{\undb_j} \mleft(\undb_{j + 1} - \undb_j\mright). \]
			
			We first prove that $\rk$ is well-defined, which amounts to showing that $\rk(\unda)$ does not depend on the choice of the increasing path $\mleft(\undb_j\mright)$. Two different such paths can be linked by a finite sequence of increasing paths such that between two consecutive increasing paths in the sequence, the only change is an inversion between two consecutive elementary moves $\unde_i$ and $\unde_j$, $i \neq j$. We thus have to check that, for every $\unda \in \hcube{\undr}{}$ and $i, j \in \ssI_{\unda}$ with $i \neq j$, we have
			\[ \ssrho_{\unda}(\unde_i) + \ssrho_{\unda + \unde_i}(\unde_j) = \ssrho_{\unda}(\unde_j) + \ssrho_{\unda + \unde_j}(\unde_i). \]
			But by \ref{def:axiom_coherent_complex_2}, we have $\ssrho_{\unda + \unde_i}(\unde_j) = \ssrho_{\unda \mycontr i}(\unde_j) = \ssrho_{\unda}(\unde_j + \unde_i) - \ssrho_{\unda}(\unde_i)$ in the left-hand part and $\ssrho_{\unda + \unde_j}(\unde_i) = \ssrho_{\unda \mycontr j}(\unde_i) = \ssrho_{\unda}(\unde_i + \unde_j) - \ssrho_{\unda}(\unde_j)$ in the right-hand part, so the desired equality holds.
			
			We now check~\ref{def:axiom_rank_function_1}-\ref{def:axiom_rank_function_2}-\ref{def:axiom_rank_function_3}. It is obvious by construction that $\rk$ takes integer values, is non-decreasing and that $\rk(t \, \unde_i) - \rk((t - 1) \, \unde_i) = 0$ or $1$ for all $i \in \{1, \dots, \grdeg\}$ and $0 \leq t \leq \ssr_i$. These imply ~\ref{def:axiom_rank_function_1} and~\ref{def:axiom_rank_function_2}.
			
			It remains to show that $\rk$ is submodular. By Theorem~\ref{thm:eq_submodularity_diamond_property}, it is sufficient to check the diamond property. We show unidirectional submodularity at distance one stated in \eqref{eq:weak_submodularity}. Let thus $i \in \{1, \dots, \grdeg\}$ and $\undx, \undy \in \hcube{\undr}{}$ such that $\undx \subfaceeq \undy$, $i \in \ssI_{\undx}$ and $\ssx_i = \ssy_i$. We assume that $\rk_{\undx}(\unde_i) = 0$ and show that $\rk_{\undy}(\unde_i) = 0$. But this has been shown to be the case in Remark~\ref{rem:loop_coherent_complex}.
			
			We have defined two maps linking coherent complexes of matroids and matricubes. It is straightforward to check that they are inverse of each other.
		\end{proof}

\section{Further discussions} \label{sec:further_discussions}
	
	In this final section, we discuss further related results and questions.
		
	\subsection{Bases of matricubes and special features of independents} \label{sec:special_features_independents}

		We do not know how to define a good notion of bases for matricubes. We review some natural attempts in this section. For each definition, we show with an example that the data of the set of bases according to that definition does not determine the matricube in a unique way. Below, $\matricube$ refers to a matricube of rank $r$ on the ground set $\hcube{\undr}{}$ and $\unda$ is an element of $\hcube{\undr}{}$.
		
		First, consider the idea closest to that of matroids.
		
		\begin{enumerate}[label=(\alph*)]
			\item \label{def:tentative_definition_bases_a} \emph{A basis of $\matricube$ is an independent $\unda \in \Ind$ which is maximal for the partial order $\subfaceeq$.}
		\end{enumerate}
		
		This does not carry enough information. The two matricubes below, with $\undr = (2, 2)$, have the same set of bases according to Definition~\ref{def:tentative_definition_bases_a}, but not the same sets of independents. The elements verifying~\ref{def:tentative_definition_bases_a} are highlighted in red, the other independents in blue.
		\begin{equation} \label{ex:bases_a}
			\begin{pmatrix}
				\color{blue} 2 & \color{blue} 3 & \color{red} 4 \\
				\color{blue} 1 & \color{blue} 2 & \color{blue} 3 \\
				\color{blue} 0 & \color{blue} 1 & \color{blue} 2
			\end{pmatrix}
			\qquad \qquad
			\begin{pmatrix}
				\color{blue} 2 & 2 & \color{red} 3 \\
				\color{blue} 1 & 1 & 2 \\
				\color{blue} 0 & \color{blue} 1 & \color{blue} 2
			\end{pmatrix}
		\end{equation}
		
		We note that, unlike matroids, a maximal independent of a matricube is not necessarily of maximal rank. This is not visible in Example~\eqref{ex:bases_a}, but the following matricube provides such an example. Two maximal independents (in red) have distinct ranks.
		\begin{equation} \label{ex:bases_a_not_same_rank}
			\begin{pmatrix}
				\color{red} 2 & 2 & 3 \\
				\color{blue} 1 & \color{blue} 2 & \color{red} 3 \\
				\color{blue} 0 & \color{blue} 1 & \color{blue} 2
			\end{pmatrix}
		\end{equation}
		Consider the following alternative to~\ref{def:tentative_definition_bases_a}.
		
		\begin{enumerate}[label=(\alph*)]
			\setcounter{enumi}{1}
			
			\item \label{def:tentative_definition_bases_b} \emph{A basis of $\matricube$ is an independent $\unda \in \Ind$ of maximal rank $r$.}
		\end{enumerate}
		The same examples given in~\eqref{ex:bases_a} show that this does not work neither.
				
		Definitions~\ref{def:tentative_definition_bases_a} and~\ref{def:tentative_definition_bases_b} are global. Seeking for local counterparts, similar to flats, circuits and independents, treated in the previous sections, leads to the following candidates.
		
		\begin{enumerate}[label=(\alph*)]
			\setcounter{enumi}{2}
			
			\item \label{def:tentative_definition_bases_c} \emph{A basis of $\matricube$ is an independent $\unda \in \Ind$ which is locally maximal, in the sense that for every $i \in \{1, \dots, \grdeg\}$ with $\unda + \unde_i \in \hcube{\undr}{}$, we have $\unda + \unde_i \notin \Ind$.}
			
			\item \label{def:tentative_definition_bases_d} \emph{A basis of $\matricube$ is an independent $\unda \in \Ind$ which is locally of maximal rank, in the sense that for every $i \in \{1, \dots, \grdeg\}$ with $\unda + \unde_i \in \hcube{\undr}{}$, we have $\rk(\unda + \unde_i) = \rk(\unda)$.}
		\end{enumerate}
		
		It is immediate that~\ref{def:tentative_definition_bases_d} implies~\ref{def:tentative_definition_bases_c}. Definitions~\ref{def:tentative_definition_bases_c} and~\ref{def:tentative_definition_bases_d} are in fact equivalent (the proof is omitted). The two matricubes below, with $\undr = (5, 4)$, have the same set of bases according to Definition~\ref{def:tentative_definition_bases_c}, but not the same sets of independents. The elements verifying~\ref{def:tentative_definition_bases_c} are depicted in red, and the other independents are in blue.
		\begin{equation} \label{ex:bases_c}
			\begin{pmatrix}
				\color{red} 4 & 4 & 4 & 4 & 5 & 6 \\
				\color{blue} 3 & 3 & \color{red} 4 & 4 & \color{blue} 5 & \color{red} 6 \\
				\color{blue} 2 & 2 & \color{blue} 3 & 3 & 4 & 5 \\
				\color{blue} 1 & 1 & 2 & 3 & 4 & 5 \\
				\color{blue} 0 & \color{blue} 1 & \color{blue} 2 & \color{blue} 3 & \color{blue} 4 & \color{red} 5
			\end{pmatrix}
			\qquad \qquad
			\begin{pmatrix}
				\color{red} 4 & 4 & 5 & 5 & 5 & 6 \\
				\color{blue} 3 & \color{blue} 4 & \color{red} 5 & 5 & 5 & \color{red} 6 \\
				\color{blue} 2 & \color{blue} 3 & \color{blue} 4 & 4 & 4 & 5 \\
				\color{blue} 1 & \color{blue} 2 & \color{blue} 3 & 3 & 4 & 5 \\
				\color{blue} 0 & \color{blue} 1 & \color{blue} 2 & \color{blue} 3 & \color{blue} 4 & \color{red} 5
			\end{pmatrix}
		\end{equation}
		
		The removal operation defined in Section~\ref{subsec:removal_and_size} takes an independent and produces smaller independents. Consider the following candidate.
		
		\begin{enumerate}[label=(\alph*)]
			\setcounter{enumi}{4}
			
			\item \label{def:tentative_definition_bases_e}\emph{A basis of $\matricube$ is an independent $\unda \in \Ind$ which is not the removal $\undb \setminus i$ for some $\undb \in \Ind$ and $i \in \{1, \dots, \grdeg\}$.}
		\end{enumerate}
		
		This is genuinely new, but, if we consider the matricubes in Example~\eqref{ex:bases_a}, the bases given by~\ref{def:tentative_definition_bases_e} and~\ref{def:tentative_definition_bases_a} are the same. 
		
		In a matroid with ground set $E$ and rank function $\rk$, every subset $S \subseteq E$ satisfies $\rk(S) + \rk^*(S^c) \leq |E|$, where $\rk^*$ is the dual rank function and $S^c$ is the complement of $S$, with equality if, and only if, $S$ is a basis. This leads to the following candidate.
		
		\begin{enumerate}[label=(\alph*)]
			\setcounter{enumi}{5}
			
			\item \label{def:tentative_definition_bases_f} \emph{A basis of $\matricube$ is an independent $\unda \in \Ind$ such that $\rk(\unda) + \rk^*(\unda^c) = \ellone{\undr}$, where $\unda^c = \undr - \unda$ is the complement of $\unda$ and $\ellone{\undr} = \sum_i \ssr_i$.}
		\end{enumerate}
		
		The inequality $\rk(\undx) + \rk^*(\undx^c) \leq \ellone{\undr}$ does hold for every element $\undx \in \hcube{\undr}{}$. However, some matricubes have no bases at all according to this definition. This is for instance the case for both matricubes in Example~\eqref{ex:bases_c}.
		
		The question of finding a good notion of bases in matricubes therefore remains open.
	
	\subsection{The natural polymatroid and the natural matroid associated to a matricube} \label{subsec:matricube_polymatroid}
		
		We refer to~\cite{HH02} for the definition and basic properties of polymatroids. Let $\polymat$ be an integer polymatroid on the ground set $E$ with rank function $\rho \colon 2^E \to \Z_{\geq 0}$. Replace each element $e$ of $E$ with $\rho(e)$ elements, and let $\widehat E$ be the resulting set. For each subset $S\subseteq E$, let $\widehat S \subseteq \widehat E$ be the union of all the elements associated to each $e \in S$. Define $\widehat \rho \colon 2^{\widehat E} \to \Z_{\geq 0}$ by the formula
		\[
			\widehat \rho(Y) \coloneqq \min_{S \subseteq E} \mleft(\rho(S) + |Y \setminus \widehat S|\mright).
		\]
		This defines a matroid on the ground set $\widehat E$, called the \emph{natural matroid} of $\polymat$, which is symmetric with respect to the permutation of the $\rho(e)$ elements associated to each $e$.
		
		To a given matricube $\matricube$ on the ground set $\hcube{\undr}{}$, we can associate an integer polymatroid $\polymat$ on the ground set $E \coloneqq [\ssr_1] \sqcup \dots \sqcup [\ssr_\grdeg]$, that we name the \emph{natural polymatroid} of $\matricube$. The rank function $\rho$ of $\polymat$ is given by associating to any subset $S \subseteq E$ the integer value $\rho(S) \coloneqq \rk\bigl(\bigvee_{\unda \in S} \unda\bigr)$. The join is taken in $\hcube{\undr}{}$ and $\rk$ is the rank function of $\matricube$.
		
		There is, moreover, a natural way to send elements of $\matricube$ to subsets of $E$, by mapping every $\undx = (\ssx_1, \dots, \ssx_\grdeg) \in \hcube{\undr}{}$ to the subset $\psi(x) \coloneqq [0, \ssx_1] \sqcup \cdots \sqcup [0, \ssx_\grdeg] \subset E$, where each interval $[0, \ssx_i]$ is taken inside $[0, \ssr_i]$. Note that the polymatroid $\polymat$ can be viewed as a polytope in the vector space $\R^E$, but there does not seem to be a natural way of associating vectors in this polytope to elements of the matricube.
		
		Proceeding as above, we can thus associate to a given matricube $\matricube$ a natural matroid $\mat$ on the ground set $\widehat E$. 
		Note that $\mat$ has $\sum_{i = 1}^{\grdeg} \sum_{t = 1}^{\ssr_i} \rk(t \, \unde_i)$ elements.
		
		Using this construction, it seems natural to transfer the notions of flats, independents, bases, and circuits from the matroid $\mat$ to the matricube, in the spirit of the work~\cite{BCF23} on integer polymatroids. This however gives a story complementary to the theory exposed in this paper.
		
		In the case of flats, for example, our definition coincides with the definition of flats in the corresponding polymatroid, in the sense that the map $\psi$ described above establishes a bijection between the flats of $\matricube$ and the flats of $\polymat$. However, we are not aware of any intrinsic axiomatic system for flats in polymatroids, and the one for flats in matricubes given in the present paper does not seem to be directly related to the one for cyclic flats in polymatroids, due to Csirmaz~\cite{csirmaz2020cyclic}; see as well~\cite[Section~5]{BCF23}.
		
		When it comes to independents of matricubes, our definition differs entirely from that of independents in a polymatroid~\cite[Section~3]{BCF23}. As we observed in the previous section, maximal independents in matricubes can be of various ranks, whereas maximal independents in polymatroids all have the same rank. The independents in polymatroids are only defined as vectors in the corresponding polytope, not set-theoretically. Besides, as we mentioned previously, we do not have yet a good notion of bases for matricubes.
		
		The same situation holds for circuits: the definition and axiomatic system we give in the present paper differ from the ones given for polymatroids in ~\cite[Section~4]{BCF23}. Like independents, circuits in polymatroids can only be defined as vectors. Our definition of circuits does not rely on independents, and yields a different story. Note in particular that in matricubes, we can have comparable circuits for the partial order $\subfaceeq$ (see Example~\eqref{ex:circuits_matricube} in Section~\ref{subsec:circuits_and_duality}), whereas two distinct circuits in a polymatroid are never comparable.
		
	\subsection{Representability and minors} \label{subsec:disc_minor}
		
		We do not know whether the local obstructions given by Proposition~\ref{prop:local_obstructions} are the only obstructions for the representability of a matricube. On the other hand, the representability of a matricube over an infinite field is equivalent to the representability of the corresponding natural matroid, as we show in the next section.
		
		\begin{thm} \label{thm:representability_matricube_matroid}
			A matricube $\matricube$ is representable over an infinite field $\k$ if, and only if, the corresponding natural matroid $\mat$ is representable over $\k$.
		\end{thm}
		
		We will deduce this from the observation, firstly, that the representability of a matricube by a flag arrangement over \emph{any} field is equivalent to the representability of the associated natural polymatroid by a subspace arrangement over the same field. And secondly, we show that the representability of an integer polymatroid by a subspace arrangement over an \emph{infinite} field is equivalent to the representability of the natural matroid associated to the polymatroid over that field. As the proof shows, the latter statement, as well as the theorem, remains valid for a \emph{finite} field of large enough cardinality.
		
		\begin{proof}
			We first prove the equivalence between the representability of a matricube and the representability of the associated natural polymatroid. We take the dual point of view described in Section~\ref{subsec:representable_matricubes}, using increasing initial flags. Let $\matricube$ be a representable matricube. The associated natural polymatroid $\polymat$ is obviously representable, remembering only the data of the subspace arrangement coming from the flag arrangement that represents $\matricube$. In the other direction, we assume that the natural polymatroid $\polymat$ associated to a matricube $\matricube$ is represented by a subspace arrangement $\mleft\{\ssfiltg_j^i, \, i \in \{1, \dots, \grdeg\}, j \in [\ssr_i]\mright\}$. Let $i \in \{1, \dots, \grdeg\}$ and $t \in [\ssr_i - 1]$. By construction of $\polymat$, we have $\rho([j + 1]) = \rk((j + 1) \, \unde_i) \geq \rk(j \, \unde_i) = \rho([j])$, where $[j]$ and $[j + 1]$ are included in the interval $[\ssr_i]$ in the ground set of $\polymat$. This shows that for every $i$ and $j$, we have $\ssfiltg_j^i \subseteq \ssfiltg_{j + 1}^i$, and therefore that the subspace arrangement can be arranged into a flag arrangement in a compatible way. It is immediate to see that this flag arrangement is a representation of $\matricube$.
			
			We now prove the equivalence between the representability of an integer polymatroid (whether or not associated to a matricube) and the representability of its natural matroid, using vector representations. Let first $\mat$ be a matroid on a ground set $\widehat E$, with rank function $\widehat \rho$, and consider a partition $\widehat E = \bigsqcup_{e \in E} A_e$ of $\widehat E$, indexed by a set $E$. We assume that $\mat$ is represented by a configuration of vectors $\mleft\{v_x \in \VS, \, x \in \widehat E\mright\}$ in a $\k$-vector space $\VS$. We define, for every $S \subseteq E$, $\rho(S) \coloneqq \widehat \rho(\widehat S) = \widehat \rho\mleft(\bigcup_{e \in S} A_e\mright)$. It is easy to see that $\rho$ is submodular on $\ssub{2}!^E$, and therefore it is the rank function of a polymatroid $\polymat(\mat)$ on the ground set $E$. Moreover, $\polymat(\mat)$ is represented by the subspace arrangement $\{\ssfiltg_e \subseteq \VS, \, e \in E\}$, where, for every $e \in E$, $\ssfiltg_e \coloneqq \langle v_x, \, x \in A_e \rangle$ is the vector subspace generated by the $v_x$, $x \in A_e$. Finally, if $\mat$ is the natural matroid of some integer polymatroid $\polymat$, then $\polymat(\mat)$ is in fact the polymatroid $\polymat$, which concludes.
			
			In the other direction, we use the notation of Section~\ref{subsec:matricube_polymatroid}. Let $\polymat$ be an integer polymatroid on a ground set $E$, represented by a subspace arrangement $\{\ssfiltg_e \subseteq \VS, \, e \in E\}$. Let $\mat$ be the natural matroid of $\polymat$, on the ground set $\widehat E$. For every $e \in E$, let $\mathcal B_e$ be a generic vector basis of the subspace $\ssfiltg_e$, i.e., let $\mathcal B_e$ be chosen in a Zariski dense open subset (to be specified afterwards) of the variety of bases of $\ssfiltg_e$. The (disjoint) union $\mathcal B$ of the bases $\mathcal B_e$ is indexed by $\widehat E$ in the natural way, say $\mathcal B = \mleft\{v_x, \, x \in \widehat E\mright\}$. Let $\overbar \rho$ be the rank function on $\widehat E$ defined by $\overbar \rho(Y) \coloneqq \dim_\k \langle v_x, \, x \in Y \rangle$ for every $Y \subseteq \widehat E$. We show the natural matroid $\mat$ to be representable by proving that $\overbar \rho = \widehat \rho$. We therefore fix any $Y \subseteq \widehat E$. For every $S \subseteq E$, the inequality $\overbar \rho(Y) \leq \rho(S) + |Y \setminus \widehat S|$ is immediate. In the rest of the proof, we show the reverse inequality.
			
			Let $S \subseteq E$ be the set of all $e \in E$ such that $\ssfiltg_e \subseteq \langle v_x, \, x \in Y \rangle$. Let $\ssfiltg_0$ be the subspace of $\VS$ defined by $\ssfiltg_0 \coloneqq \sum_{e \in S} \ssfiltg_e$, and let $\ssfiltg_1, \dots, \ssfiltg_k$ be the other subspaces $\ssfiltg_e$ for $e \in E \setminus S$. Here, $k = |E| - |S|$. Rearrange the ground set $\widehat E$ so that the first $k$ disjoint intervals $[\ssr_1], \dots, [\ssr_k]$ that make it up are those indexed by $E \setminus S$. For every $j \in \{1, \dots, k\}$, let $\ssY_j$ be the family of vectors $v_x$ with $x \in Y \cap [\ssr_j]$ in the ground set $\widehat E$ of $\mat$. We then admit, for now, that for every $j \in \{1, \dots, k\}$, the following inequality holds:
			\begin{equation} \label{eq:dimension_bases_polymatroid}
				\dim(\ssfiltg_0 + \ssfiltg_1 + \cdots + \ssfiltg_j) > \dim(\ssfiltg_0) + |\ssY_1| + \cdots + |\ssY_j|.
			\end{equation}
			We show by induction that if Inequality~\eqref{eq:dimension_bases_polymatroid} is true for every $j \in \{1, \dots, k\}$, then, for every $j \in \{0, \dots, k\}$, the configuration of vectors $\ssY_1 \sqcup \cdots \sqcup \ssY_j$ is linearly independent in the quotient space $\rquot{\bigl(\ssfiltg_0 + \ssfiltg_1 + \cdots + \ssfiltg_j\bigr)}{\ssfiltg_0}$. This is obviously true for $j = 0$. If now this is true for some $j < k$, Inequality~\eqref{eq:dimension_bases_polymatroid} for $j + 1$ reads
			\[ |\ssY_{j + 1}| < \dim(\ssfiltg_0 + \ssfiltg_1 + \cdots + \ssfiltg_{j + 1}) - (\dim(\ssfiltg_0) + |\ssY_1| + \cdots + |\ssY_j|). \]
			The induction hypothesis then implies that the quantity on the right is the dimension of the quotient space
			\[ \rquot{\bigl(\ssfiltg_0 + \ssfiltg_1 + \cdots + \ssfiltg_{j + 1}\bigr)}{\bigl(\ssfiltg_0 + \langle \ssY_1, \dots, \ssY_j \rangle \bigr)}. \]
			The linear projection map
			\[ \rquot{\bigl(\ssfiltg_0 + \ssfiltg_1 + \cdots + \ssfiltg_j + \ssfiltg_{j + 1} \bigr)}{\ssfiltg_0} \twoheadrightarrow \rquot{\bigl(\ssfiltg_0 + \ssfiltg_1 + \cdots + \ssfiltg_j + \ssfiltg_{j + 1}\bigr)}{\bigl(\ssfiltg_0 + \langle \ssY_1, \dots, \ssY_j \rangle\bigr)} \]
			now enables to view the family of vectors $\ssY_{j + 1}$ in the quotient space $\rquot{\bigl(\ssfiltg_0 + \ssfiltg_1 + \cdots + \ssfiltg_{j + 1}\bigr)}{\bigl(\ssfiltg_0 + \langle \ssY_1, \dots, \ssY_j \rangle\bigr)}$. Since the number of vectors in $\ssY_{j + 1}$ is less than the dimension of this quotient space, and since all the vectors $v_x$ are generic, then $\ssY_{j + 1}$ is linearly independent in this quotient space. This means exactly that $\ssY_1 \sqcup \cdots \sqcup \ssY_j \sqcup \ssY_{j + 1}$ is linearly independent in the quotient space $\rquot{\bigl(\ssfiltg_0 + \ssfiltg_1 + \cdots + \ssfiltg_j + \ssfiltg_{j + 1}\bigr)}{\ssfiltg_0}$, which concludes the induction. Specializing the independence property to $j = k$, we get the desired equality:
			\[ \overbar \rho(Y) = \rho(S) + |Y \setminus \widehat S|. \]
			
			To finish the proof, before turning to Inequality~\eqref{eq:dimension_bases_polymatroid}, notice that the choice of the vectors $v_x$ determines the subset $S \subseteq E$. But since there is a finite number of such $S$, there is still a non-empty Zariski open set of choices of vectors $v_x$ for which the same set $S$ is associated to all these choices. Now, we have the above equality for a fixed $Y$. Since there is a finite number of such subsets $Y \subseteq \widehat E$, there is still a non-empty Zariski open set of choices of the vectors $v_x$ for which the above equality holds for \emph{every} $Y$.
			
			We now explain why Inequality~\eqref{eq:dimension_bases_polymatroid} is true for every $j$. If by contradiction it was not, let $j_0 \geq 1$ be the smallest $j$ such that it does not hold, i.e.,
			\[ |\ssY_{j_0}| \geq \dim(\rquot{\ssfiltg_0 + \ssfiltg_1 + \cdots + \ssfiltg_{j_0}}{\ssfiltg_0 + \langle \ssY_1, \dots, \ssY_{j_0 - 1} \rangle}). \]
			The preceding argument shows that the vectors in $\ssY_j$ for $j < j_0$ are independent in the quotient space $\rquot{\bigl(\ssfiltg_0 + \ssfiltg_1 + \cdots + \ssfiltg_{j_0 - 1}\bigr)}{\ssfiltg_0}$. Since the vectors in $\ssY_{j_0}$ are generic, the above inequality on the size of $\ssY_{j_0}$ implies that they generate the quotient space $\rquot{\bigl(\ssfiltg_0 + \ssfiltg_1 + \cdots + \ssfiltg_{j_0}\bigr)}{\bigl(\ssfiltg_0 + \langle \ssY_1, \dots, \ssY_{j_0 - 1} \rangle\bigr)}$. This leads to the inclusion $\ssfiltg_{j_0} \subseteq \langle v_x, \, x \in Y \rangle$, contradicting the definition of $S$.
		\end{proof}
		
		We call a matricube $\matricube$ \emph{regular} if it is representable over every field.
		For matroids, Seymour's theorem describes regular matroids in terms of sums of graphic, cographic, and an exceptional regular matroid on 10 elements~\cite{Sey80}. A theorem by Tutte characterizes regular matroids as those representable over the fields $\mathbb F_2$ and $\mathbb F_3$ with two and three elements, respectively. Another result by Tutte characterizes regular matroids as those that contain as a minor neither the Fano matroid $\mathrm F_7$ nor its dual $\mathrm F_7^*$. We refer to~\cite[Chap. 9]{Tru92} for a presentation of these results.
		
		\begin{question}
			Provide a characterization of regular matricubes.
		\end{question}
		
		In analogy with Rota's conjecture on the characterization of the representability of matroids over finite fields using a finite set of forbidden minors, we formulate the following question.
		
		\begin{question}[Rota's conjecture for matricubes]
			Let $\k$ be a finite field. Does there exist a finite collection of matricubes such that a matricube is representable over $\k$ if, and only if, it does not contain any of the matricubes in the collection as a minor?
		\end{question}
		
		We note that as it was recently shown by Oxley, Semple, and Whittle~\cite{OSW16}, the analogue of Rota's conjecture for 2-polymatroids fails in general. This does not exclude a positive answer to the above question, as matricubes behave more like matroids than polymatroids.
	
	\subsection{Stratification of products of flag varieties} \label{sec:stratification}
		
		Let $\k$ be a field and $n$ be a positive integer. Let $\VS$ be a $\k$-vector space of dimension $n$. For each positive integer $\ssr \leq n$, denote by $F(\ssr, n)$ the flag variety parametrizing initial flags of vector subspaces $\ssfiltg_1 \subset \ssfiltg_2 \subset \dots \subset \ssfiltg_{\ssr}$ of dimensions $\dim(\ssfiltg_j) = j$, for $j = 1, \dots, \ssr$. Given a vector $\undr = (\ssr_1, \dots, \ssr_\grdeg)$ with $\grdeg$ positive integers, consider the product variety $F\mleft(\undr, n\mright) \coloneqq F(\ssr_1, n) \times \dots \times F(\ssr_\grdeg, n)$. We get a natural stratification of $F\mleft(\undr, n\mright)$ by matricubes as follows. Given a simple matricube $\matricube$ with ground set $\hcube{\undr}{}$, the cell $\cell_{\matricube}$ parametrizes those collections of $\grdeg$ flags $\ssfiltg_{\bullet}^i \in F(\ssr_i, n)$, $i = 1, \dots, \grdeg$, whose associated matricube, through the constructions of Section~\ref{subsec:representable_matricubes}, coincides with $\matricube$. This stratification is analogous to that induced by matroids for Grassmannians. Theorem~\ref{thm:representability_matricube_matroid} provides a correspondence between strata of given rank $r$ in $F\mleft(\undr, n\mright)$ and strata of the Grassmannian $\Gr(r, N)$, for $N = \frac 12 \sum_{i = 1}^\grdeg (\ssr_i^2 + \ssr_i)$, see Section~\ref{subsec:matricube_polymatroid}. It would be interesting to study the combinatorics of this stratification, and the geometric meaning of this correspondence.
	
	\subsection{Polymatricubes} \label{sec:polymatricubes}
		
		A polymatricube is a function $f \colon \hcube{\undr}{} \to \R$ with $f(\underline 0) = 0$ which is non-decreasing and submodular, that is, it verifies~\ref{def:axiom_rank_function_2} and~\ref{def:axiom_rank_function_3}. Examples of polymatricubes are the representable ones which, by definition, are those associated to a collection of \emph{arbitrary} (instead of initial) flags in a vector space. Generalizing the discussion of Section~\ref{subsec:matricube_polymatroid}, we can associate a natural polymatroid and a natural matroid to any polymatricube. Theorem~\ref{thm:representability_matricube_matroid} extends to this setting: a polymatricube is representable over an infinite field if, and only if, the corresponding natural matroid is representable over the same field. In particular, the discussion of Section~\ref{sec:stratification} can be extended to arbitrary collections of flag varieties.
	
	\subsection{Tutte polynomial} \label{subsec:Tutte_polynomial}
		
		An important algebraic invariant associated to a matroid is its Tutte polynomial. This is a two-variable polynomial that specializes to the characteristic polynomial of the matroid. The Tutte polynomial of a matroid $\mat$ on the ground set $E$ is the unique polynomial $\ssT_\mat(\x, \y)$ that verifies the following recursive equation for every $e \in E$:
		\begin{align*}
			\ssT_\mat(\x, \y) = \begin{cases} 
				\x \, \ssT_{\mat \mysetminus e}(\x, \y) & \textrm{ if } e \textrm{ is a coloop} \\
				\y \, \ssT_{\mat \mycontr e}(\x, \y) & \textrm{ if } e \textrm{ is a loop} \\
				\ssT_{\mat \mysetminus e}(\x, \y) + \ssT_{\mat \mycontr e}(\x, \y) & \textrm{ if } e \textrm{ is neither a loop nor a coloop},
			\end{cases}
		\end{align*}
		and is defined for the matroid $\emptyset$ with empty ground set by $\ssT_{\emptyset} \equiv 1$.
		
		We can define the notion of loop and coloop in matricubes. We say $i \in \{1, \dots, \grdeg\}$ is a \emph{loop} of $\matricube$ if $\ssr_i > 0$ and $\rk(\unde_i) = 0$. We say $i$ is a \emph{coloop} of $\matricube$ if $i$ is a loop of the dual matricube $\matricube^*$. This is equivalent to having $\rk(\matricube \mysetminus i) = \rk(\matricube) - 1$. The recursive equation above, however, does not lead to an invariant of matricubes.
		
		The Tutte polynomial of a matroid $\mat$ on the ground set $E$ can be defined directly by the following formula:
		\[
			\ssT_\mat(\x, \y) = \sum_{S \subseteq E} (\x - 1)^{r - \rk(S)} (\y - 1)^{|S| - \rk(S)}.
		\]
		
		This definition naturally extends to any matricube.
		
		\begin{defi} \label{def:Tutte_polynomial_matricubes}
			Let $\matricube$ be a matricube of rank $r$ on the ground set $\hcube{\undr}{}$. The \emph{Tutte polynomial} of $\matricube$ is the two-variable polynomial
			\[
				\ssT_{\matricube}(\x, \y) \coloneqq \sum_{\undx \in \hcube{\undr}{}} (\x - 1)^{r - \rk(\undx)} (\y - 1)^{\ellone{\undx} - \rk(\undx)},
			\]
			where $\ellone{\undx}$ is the $\ell_1$-norm of $\undx$.
		\end{defi}
		
		Tutte polynomials of matricubes verify the following properties:
		
		\begin{itemize}	
			\item Let $\matricube^*$ denote the dual of the matricube $\matricube$. Then, we have $\ssT_{\matricube^*}(\x, \y) = \ssT_{\matricube}(\y, \x)$.
			
			\item For two matricubes $\matricube$ and $\matricube'$, we have $\ssT_{\matricube \oplus \matricube'} = \ssT_{\matricube} \cdot \ssT_{\matricube'}$.
		\end{itemize}
		
		Here, the direct sum $\matricube \oplus \matricube'$ of matricubes $\matricube$ and $\matricube'$ on hypercuboids $\hcube{\undr}{}$ and $\hcube{\undr}{}$, respectively, has ground set $\hcube{\undr}{} \times \hcube{\undr}{}$ and rank function defined by
		\[ \rk_{\matricube \oplus \matricube'}(\undx \oplus \undx') = \rk_{\matricube}(\undx) + \rk_{\matricube'}(\undx'). \]
		
		\begin{question}
			Does there exist a recursive identity which defines $\ssT_{\matricube}$, in terms of deletions and contractions?
		\end{question}
		
		There is a version of the Tutte polynomial for polymatroids defined by Cameron and Fink~\cite{CF22}. This polynomial satisfies a relation involving elementary operations reminiscent of deletion and contraction in polymatroids, and specializes to the Tutte polynomial for matroids. We do not know if there is any relation between the Tutte polynomial of a matricube and the Tutte polynomial of the corresponding natural polymatroid.
	
	\subsection{Matricubes arising from linear series on curves} \label{subsec:geometric_rank_functions}
		
		As pointed out in Section~\ref{subsec:further_related_questions}, matricubes naturally arise in our work on tropical degenerations of linear series on algebraic curves. We provide a brief hint to this by explaining how a finite collection of points and a finite dimensional space of rational functions on an algebraic curve gives rises to a matricube.
		
		Let $\k$ be an algebraically closed field, and let $C$ be a smooth proper curve over $\k$. Let $r$ be a non-negative integer, and let $p$ be a $\k$-point on $C$. Let $\k(C)$ be the function field of $C$, and let $\VS \subseteq \k(C)$ be a vector subspace of rational functions of dimension $r + 1$ over $\k$. The point $p$ leads to a complete flag $\ssfilt_p^\bullet$ of $\VS$ by looking at the orders of vanishing at $p$ of functions in $\VS$, as follows. Define the set $\ssS_{p} \coloneqq \mleft\{\ord_{p}(f) \, \st \, f \in \VS \setminus \{0\}\mright\}$. The set $\ssS_p$ is finite of cardinality $r + 1$. Denote by $\sss_0^p < \dots < \sss_r^p$ the elements of $\ssS_p$, enumerated in increasing order. The flag $\ssfilt_p^\bullet$ is defined by setting, for $j = 0, \dots, r$,
		\[ \ssfilt_p^j \coloneqq \mleft\{f \in \VS \setminus \{0\} \, \st \, \ord_{p}(f) \geq \sss_j^p \mright\} \cup \{0\}. \]
		It follows that $\ssfilt_p^j$ has codimension $j$ in $\VS$.
		
		Let now $\grdeg$ be a natural number, and let $A = \{\ssp_1, \dots, \ssp_\grdeg\}$ be a collection of $\grdeg$ distinct $\k$-points on $C$. By the construction above, each point $\ssp_i$ leads to a complete flag $\ssfilt_{i}^\bullet$. Denoting $\ssS_{i} \coloneqq \{\ord_{\ssp_i}(f) \, \st \, f \in H \setminus \{0\}\}$, and enumerating the elements of $\ssS_i$ in increasing order $\sss_0^i < \dots < \sss_r^i$, the flag $\ssfilt_{i}^\bullet$ is defined by setting
		\[ \ssfilt_i^j \coloneqq \mleft\{f \in H \setminus \{0\} \, \st \, \ord_{\ssp_i}(f) \geq \sss_j^i\mright\} \cup\{0\}. \]
		The data of $C, \VS, \ssp_1, \dots, \ssp_\grdeg$ defines a matricube on the ground set $\hcube{r}{\grdeg} = [r]^\grdeg$.
		
		We may call \emph{geometric} a matricube with ground set $\hcube{r}{\grdeg}$ that arises from the above construction for a curve $C$ over an algebraically closed field $\k$. By construction, geometric matricubes are all representable over the field $\k$ over which the curve $C$ is defined.
		
		\begin{question}
			Is it true that all representable matricubes on $\hcube{r}{\grdeg}$ of rank $r$ or $r + 1$ are geometric? What is the smallest possible genus of a curve representing a geometric matricube?
		\end{question}

\bibliographystyle{alpha}
\bibliography{Bibliography}

\end{document}